\documentclass[%
]{mpi2015-cscpreprint}
\usepackage[american]{babel}
\usepackage{graphicx}
\usepackage{algorithm}
\usepackage{algorithmic}
\usepackage{amssymb}
\usepackage{amsthm,amsmath}
\usepackage{ulem}
\usepackage{todonotes}
\usepackage{mathtools}

\newtheorem{theorem}{Theorem}[section]
\newtheorem{lemma}{Lemma}[section]
\newtheorem{definition}{Definition}[section]

\newtheorem{remark}{Remark}[section]

\def\IR{{\mathbb R}}
\def\IC{{\mathbb C}}
\def\IZ{{\mathbb Z}}

\def\IL{{\mathbb L}}
\def\sIL{{\mathbb L}_s}

\def\IV{{\mathbb V}}
\def\IW{{\mathbb W}}

\newcommand{\bA}{{\textbf A}}
\newcommand{\bB}{{\textbf B}}
\newcommand{\bC}{{\textbf C}}
\newcommand{\bD}{{\textbf D}}
\newcommand{\bE}{{\textbf E}}
\newcommand{\bF}{{\textbf F}}
\newcommand{\bG}{{\textbf G}}
\newcommand{\bS}{{\textbf S}}
\newcommand{\bY}{{\textbf Y}}
\newcommand{\bJ}{{\textbf J}}
\newcommand{\bK}{{\textbf K}}

\newcommand{\bL}{{\textbf L}}
\newcommand{\bM}{{\textbf M}}
\newcommand{\bN}{{\textbf N}}
\newcommand{\bI}{{\textbf I}}
\newcommand{\bH}{{\textbf H}}

\newcommand{\bO}{{\textbf O}}
\newcommand{\bQ}{{\textbf Q}}
\newcommand{\bW}{{\textbf W}}
\newcommand{\bR}{{\textbf R}}
\newcommand{\bT}{{\textbf T}}
\newcommand{\bZ}{{\textbf Z}}
\newcommand{\bX}{{\textbf X}}
\newcommand{\bx}{{\textbf x}}

\newcommand{\bu}{{\textbf u}}
\newcommand{\bV}{{\textbf V}}
\newcommand{\bU}{{\textbf U}}

\newcommand{\bs}{{\textbf s}}

\newcommand{\bp}{{\textbf p}}

\newcommand{\bq}{{\textbf q}}

\newcommand{\bv}{{\textbf v}}

\newcommand{\bfz}{{\mathbf 0}}

\newcommand{\cO}{ {\cal O} }

\newcommand{\bPhi}{ \boldsymbol{\Phi} }

\newcommand{\bPsi}{ \boldsymbol{\Psi} }

\newcommand{\blambda}{\boldsymbol{\lambda}}

\newcommand{\bmu}{\boldsymbol{\mu}}
\newcommand{\bsigma}{\boldsymbol{\sigma}}
\newcommand{\bSigma}{\boldsymbol{\Sigma}}

\newcommand{\bLambda}{\boldsymbol{\Lambda}}

\newcommand{\cR}{ {\cal R} }

\newcommand{\cL}{ {\cal L} }
\begin{document}
\title{\Large{Data-driven quadratic modeling in the Loewner framework from input-output time-domain measurements}}
\novelty{Inference of quadratic control models from input-output time-domain data that serve the objectives of reduction and identification in non-trivial equilibria.}
\author[$\ast$]{\large{D.~S. Karachalios}}
\affil[$\ast$]{Max Planck Institute for Dynamics of Complex Technical Systems, Magdeburg, Germany and current affiliation: Institute for Electrical Engineering in Medicine, University of Luebeck, Germany.\authorcr
  \email{karachalios@mpi-magdeburg.mpg.de}, \orcid{0000-0001-9566-0076}}
\author[$\dagger$]{I.~V. Gosea}
\affil[$\dagger$]{Max Planck Institute for Dynamics of Complex Technical Systems, Magdeburg, Germany.\authorcr
  \email{gosea@mpi-magdeburg.mpg.de}, \orcid{0000-0003-3580-4116}}
\author[$\dagger$]{L. Gkimisis}
\affil[$\dagger$]{Max Planck Institute for Dynamics of Complex Technical Systems, Magdeburg, Germany.\authorcr
  \email{gkimisis@mpi-magdeburg.mpg.de}}  
  \author[$\ddagger$]{A.~C. Antoulas}
\affil[$\ddagger$]{Max Planck Institute for Dynamics of Complex Technical Systems, Magdeburg, Germany and Rice University, Electrical and Computer Engineering Department, Houston, TX 77005, USA and Baylor College of Medicine, Houston, TX 77030, USA.\authorcr
  \email{aca@rice.edu}}  
\shorttitle{Data-driven quadratic modeling}
\shortauthor{Karachalios et al.}
\abstract{In this study, we present a purely data-driven method that uses the Loewner framework (LF) along with nonlinear optimization techniques to infer quadratic with affine control dynamical systems that admit Volterra series (VS) representations from input-output (i/o) time-domain measurements. The proposed method extensively employs optimization tools for interpolating the symmetric generalized frequency response functions (GFRFs) derived in the VS framework. The GFRF estimations are obtained from the Fourier spectrum (phase and amplitude) of the quasi-steady state system response under harmonic excitation. Appropriate treatment of these measurements under the developed framework allows the identification of low-order quadratic state-space models with non-trivial stable equilibria, such as in the Lorenz '63 forced system. We thus can achieve low-order global model identification for systems that can bifurcate to multiple equilibria after solely collecting measurements from a local stable operational regime. The developed framework is tested for several examples of increasing dimension and complexity, up to a test case of the viscous Burgers' equation with Robin boundary conditions. In the latter case, this study enforces new directions of employing data-driven reduced model inference that successfully can provide low-order accurate surrogate predictive models suitable for control. Future directions and open challenges conclude this work.}
\maketitle  
\section{Introduction}
\paragraph{General model order reduction approaches and data-driven emerging} Mathematical modeling should be consistent with physical laws when considering engineering applications. Developing surrogate models for robust simulation, design, and control is a worthwhile engineering task, and methods that can interpret or discover the governing equations \cite{brunton2022data} are essential. Numerical discretization of partial differential equations (PDEs) results in large-scale systems of ordinary differential equations (ODEs) \cite{ACA05}. Even with the recent development of high-performance computing (HPC) environments, the resulting dynamical models after the spatial discretization of complex PDEs still inherit high complexity to handle efficiently. Therefore, approximation of large-scale dynamical systems \cite{ACA05} is crucial for efficient simulation. The technique to reduce complexity is known as model order reduction (MOR) \cite{HandbookVol1,HandbookVol2,HandbookVol3}. There are many ways of reducing large-scale models, and each method is tailored to specific applications and aspects of complexity reduction. A good distinction among methods is the accessibility to a high-fidelity model that discriminates the methodological approaches into branches as intrusive or non-intrusive modeling, where various data-driven methods continuously emerge for the latter. For the intrusive case where a high-fidelity model is available, methods such as proper orthogonal decomposition (POD)\cite{PODWeissJulien}, balanced truncation (BT) (see the recent survey \cite{breiten20212}) moment matching (MM) (recent survey \cite{benner2021model}, and the references therein) were extensively used to build low-order surrogate models that approximate the original full-order model (FOM) without significant loss in accuracy, offering an error bound, e.g., BT, and a guarantee on stability (BT and some variants of MM). Additional MOR methods for nonlinear systems have also been developed (by extending the linear counterpart of BT, MM, or others) \cite{baur2014model,ABG20}. 

\paragraph{Data-driven snapshots based methods} On the other hand, the ever-increasing availability of data, i.e., measurements related to the original model, render non-intrusive techniques such as machine learning (ML) combined with model-based methods \cite{brunton2022data} as a promising modeling alternative. ML has demonstrated remarkable success in specific tasks, e.g., pattern recognition or feature extraction. However, issues such as interpretability and generalizability limit some existing ML methods. Recent endeavors that connect physics with machine learning methods \cite{RAISSI2019686} and learn the underlying operator as in \cite{li2021fourier}  are very promising. Data assimilation with model-based data from proprietary software through MOR techniques such as the dynamic mode decomposition (DMD)\cite{schmid_2010,ProBruKut2016,morHirHKetal20}, the operator inference (OpInf) \cite{morPehW16,morBenGKetal20,KHODABAKHSHI2022114296} has become popular. The aforementioned methods (OpInf, DMD) and others, such as the sparse identification of dynamical systems (SINDy)\cite{SINDY} assume state access, i.e., measurement availability for all system states, at all times. Besides rendering the model inference input-dependent, this assumption might be impossible to satisfy for many real-life applications.

\paragraph{Non-intrusive data-driven from i/o data} Towards the same aim of model discovery using exclusively input-output (i/o) data and for constructing models, the Loewner framework constitutes a non-intrusive method that deals directly with i/o data (real-world measurements used in most applied sciences, e.g., frequency, velocity, voltage, charge, or concentration). The Loewner framework aims to the inference (or realization as often termed in the automatic control community) of linear or nonlinear state-space systems and simultaneously offers the opportunity for data-driven complexity reduction. Out of many existing methods available, we mention those based on rational approximation; the Loewner framework \cite{ajm07}, vector fitting (VF) \cite{VF}, and the AAA algorithm \cite{NST18}. We refer the reader to the extensive analysis provided in \cite{ABG20} for more details on such methods. 

The realization aim in control theory refers to discovering the operators that define a dynamical system in state-space representation from i/o measurements. Such realization of state-space systems is not unique due to different coordinates as opposed to realizing a transfer function or GFRFs that obey a unique representation. The linear case is introduced in \cite{HOKALMAN} and has been extended further from \cite{JuangPappa}. For the nonlinear case, extensions to the realization algorithm through the subspace method and for the bilinear control system case can be found in \cite{Isidori1973, Favoreel, CHEN2000, Ramos2009} with the references within, and for the linear parametric-varying (LPV) systems in \cite{TothHossam} when the scheduling signals can be measured. Other methods for data-driven system identification or reduction based on nonlinear autoregressive moving average models with exogenous inputs (NARMAX) can be found in \cite{BillingsNARMAX,Chorin2015} and in connection with the Koopman operator and the Wiener projection in \cite{LIN2021109864}.

A viable alternative is to devise methods that directly learn the continuous in-time operators without introducing additional discretization errors over the temporal space. In particular, the Fourier transform, which maps the time domain to the frequency domain, given that the Nyquist criterion is respected, allows continuous signal reconstruction from a finite spectrum when the correct sampling frequency has been considered. Developing methods over the frequency domain where interpolation can be enforced has the advantage that the model's dimension scales mainly with the state dimension (i.e., spatial) as opposed to existing time-domain methods (e.g., subspace identification), where the model's dimension scales in addition with the time discretization dimension and leads inevitably to exponential complexity. Importantly, processing the time-domain output over the frequency domain has advantages in terms of sparse information and making it possible to utilize methods such as the Loewner framework.

\paragraph{Existing literature on the Lowner framework} The Loewner framework (LF) could be interpreted as an interpolatory MOR technique that identifies state-space systems for certain generalized nonlinear classes, particularly: bilinear systems \cite{AGI16}, linear switched systems \cite{GPAswitch}, linear parameter varying systems \cite{GPA21CDC}, quadratic-bilinear systems \cite{morGA15,morGosA18,morAntGH19}, and polynomial systems \cite{PSGoyBen}. The aforementioned nonlinear variants of the Loewner framework construct efficient surrogate models from i/o data, sampled from a high-fidelity model. The challenge of the existing nonlinear approaches is that they cannot be applied when the underlying model is inaccessible. In a natural measurement environment where i/o time-domain measurements can be collected from an actual plant, a unique, measurable set of kernels exists within the VS framework, i.e., the symmetric Volterra kernels. The symmetric GFRFs can be derived with the growing exponential approach tailored to the probing method \cite{Ru82}, under the assumption that harmonic excitation leads to harmonic output. The LF has already been extended to handle input-output time-domain data for the linear system case in \cite{PGW17}.

\paragraph{Scope of the proposed non-intrusive i/o method and applicability}
We introduce the scope and generality of our new method, which uses i/o time-domain measurements to infer quadratic state-space models with affine control after combining the Loewner and Volterra frameworks with nonlinear optimization techniques. The quadratic model class assumption results from first principles, where ODEs of a specific polynomial degree can be derived or efficiently approximated. In particular, dynamical control systems that belong to the quadratic class contain a wide variety of challenging problems with direct connection to engineering applications requiring efficient simulation, such as the incompressible Navier Stokes and variants such as the viscous Burgers' model for investigating fluid mechanics, the Lorenz system for weather prediction, the Lotka-Volterra for population dynamics, and SIRD models for epidemiology among others. Moreover, using algebraic techniques such as lifting \cite{Gu11,Gosea22mdpi}, the class of quadratic systems could incorporate wider classes of nonlinear systems with analytic nonlinearities, thus extending the applicability of the inference task hereby pursued. 

A wide class of nonlinear systems can be described through the Volterra Series (VS) (the Volterra Series framework is also known as the Volterra-Wiener) approach \cite{Ru82,Billings2013NonlinearSI}. The VS framework decomposes the nonlinear dynamics into an infinite series of convolutional operators involving generalized kernels with input multiplication. It can model weakly nonlinear dynamics due to the polynomial structure, potentially affecting its convergence properties. Thus, the VS representation is limited to a subset of nonlinear systems, i.e., weakly nonlinear for which self-sustained dynamics, phenomena such as sub-harmonic oscillations, continuous spectrum, or chaotic dynamics do not emerge \cite{Billings2003Subharmonic}. In line with the above limitation set by the VS representation, we further assume that a harmonic input will eventually lead to a harmonic system response (referred to as a quasi-steady state) associated with the so-called probing method. 

The proposed method can be applied when harmonic forcing (potentially with non-zero initial conditions) to a black box system reaches a stable equilibrium point with a quasi-steady state, and discrete spectrum measurements from the processed steady-state profile can be obtained. These measurements consist of the harmonics that the underlying system operates around the achievable stable equilibrium. The derived harmonics are described within the VS framework, where appropriate handling of those through kernel separation can lead to symmetric GFRF estimations. Further, the estimations of the first three symmetric GFRFs are connected with the interpolatory LF method, and together with nonlinear optimization methods, low-order quadratic models can be inferred with the potential to achieve global identification.

\paragraph{Compared to existing relevant literature}
In comparison with our prior studies in \cite{morKarGA19a,KARACHALIOS20227,gosea2021learning}, the advances in the newly proposed method are that we use generalized frequency measurements from high-order kernels (first three GFRFs) of the symmetric type that interpret the propagating harmonics in the time domain output, allowing for the inverse transformation (from time to frequency domain) via the Fourier transform. Further, we solve the resulting nonlinear optimization problems to achieve quadratic model inference that implicitly interpolates the whole VS after explicitly enforcing interpolation to a finite set of GFRFs (first three). Moreover, after measuring the local dynamical behavior and computing the translated non-zero equilibrium point, we can identify the global quadratic system that contains all the nontrivial equilibria.

\paragraph{Structure of the paper}
The rest of the paper is organized as follows: \Cref{sec:Analysis&Methods} presents the analysis and employed methods, progressing with handy examples from the linear case in \cref{sec:LF} up to the nonlinear in \cref{sec:Volterra} and quadratic with the data acquisition process in \cref{sec:Quadratic}. The method, along with the derived algorithms, is in \cref{sec:method}. The dynamics can be measured around nonzero equilibrium points; thus, analysis for such systems is included in \cref{sec:multiEqulibria}. In \Cref{sec:Results}, we test the method in more challenging benchmarks. Firstly, we illustrate different cases of identifying the forced Lorenz attractor with trivial and non-trivial equilibria, and the global identification feature is analyzed. Secondly, the method is also tested for the reduction/accuracy performance of the viscous Burgers' model with Robin boundary conditions. Finally, \Cref{sec:conclusions} summarizes the findings and discusses potential future research directions.

\paragraph{Notation} $\IR,~\IC$ fields of real and complex numbers and $\IR^{m\times n},~\IC^{m\times n}$ real and complex $m\times n$ matrices; $|\xi|$ absolute value of real or complex scalar; imaginary unit $j^2=-1$; for $s=a+jb$, $Re(s)=a,~Im(s)=b$ real and imaginary part of $s$, $\bar{s}=a-jb$ conjugate of $s$; $a_{ij}$ the $(i,j)$-th entry of $\bA$; $\bA(i:j,:)$ sub-matrix with rows $i,\ldots,j$ of $\bA$ and all columns; $\bA^\top$ transpose of $\bA$; $\bA^\star:=(\bar{\bA})^\top$, complex conjugate transpose; $\bA^{-1}$ inverse of nonsingular $\bA$; $\bI_n$ identity matrix of dimension $n$; $\sigma_{\text{max}}(\bA)$ the largest singular value; $\lVert u\rVert_p:=(\sum_{i=1}|u_i|^p)^{1/p}$ for $u\in\IC^n$ and $1\leq p<\infty$; $\bA\otimes\bB$ the Kronecker product of $\bA$ and $\bB$; $\texttt{vec}(\bA)$ vectorization operator; gradient $\nabla f:=(\partial_{x_1}f,\ldots,\partial_{x_n}f)^\top$.

\section{Analysis and methods}\label{sec:Analysis&Methods}
We start with linear system theory, providing explicit solutions and deriving the main mathematical notions of our method. Next, we introduce the LF as an interpolatory tool in the linear case and explain how one can identify minimal linear systems from time-domain measurements similar to \cite{PGW17}. The approximation of nonlinear systems with the VS is introduced, focusing on the quadratic state-space control system case. Analytical derivations are presented for extracting the symmetric GFRFs with the probing method. Further, the analysis of addressing the issue of local equilibrium measurements is detailed. Finally, the chapter concludes with the proposed method supported with concise algorithms summarizing the computational techniques.
\subsection{Linear time-invariant (LTI) systems}
We start our analysis with the linear time-invariant\footnote{The invariant concept for the state-space system representation has to do with the fact that the operator of the underlying system is temporal (time) independent.} (LTI) system in the single-input, single-output (SISO) form as
\begin{equation}
\label{eq:LTI}
    \left\{\begin{aligned}
    \bE\dot{\bx}(t)&=\bA\bx(t)+\bB u(t),\\
    y(t)&=\bC\bx(t)+\bD u(t),~\bx(0)=\bx_0=\textbf{0},~t\geq 0,
    \end{aligned}\right.
\end{equation}
where the state-dimension is $n$, and the matrices are: $\bE,~\bA\in\IR^{n\times n},~\bB,\bC^T\in\IR^{n\times 1},~\bD\in\IR$. Concentrating in the ordinary differential equation \cref{eq:LTI}, with a non-singular (i.e., invertible) $\bE$ matrix, we have the following explicit solution

\begin{equation}
    \begin{aligned}
    \dot{\bx}(t)=\bE^{-1}\bA\bx(t)+\bE^{-1}\bB u(t)&\Rightarrow
    \frac{d}{dt}\left[e^{-\bE^{-1}\bA t}\bx(t)\right]=e^{-\bE^{-1}\bA t}\bE^{-1}\bB u(t)\Rightarrow\\
    \int_{0}^t\frac{d}{d\tau}\left[e^{-\bE^{-1}\bA \tau}\bx(\tau)\right]d\tau&=\int_{0}^{t}e^{-\bE^{-1}\bA \tau}\bE^{-1}\bB u(\tau)d\tau\Rightarrow\\
    e^{-\bE^{-1}\bA t}\bx(t)-e^{-\bE^{-1}\bA\cdot 0}\bx(0)&=\int_{0}^{t}e^{-\bE^{-1}\bA \tau}\bE^{-1}\bB u(\tau)d\tau\Rightarrow\\
    \bx(t)=e^{\bE^{-1}\bA t}\bx(0)+&\int_{0}^{t}e^{\bE^{-1}\bA (t-\tau)}\bE^{-1}\bB u(\tau)d\tau\Rightarrow(\text{with}~\sigma=t-\tau)\\
     y(t)=\bC e^{\bE^{-1}\bA t}\bx_0+&\int_{0}^{t}\underbrace{\bC e^{\bE^{-1}\bA \sigma}\bE^{-1}\bB}_{h_1(\sigma)} u(t-\sigma)d\sigma+\bD u(t).
    \end{aligned}
\end{equation}

The input-output solution of an LTI system with zero-initial conditions $\bx_0=\textbf{0}$ and a zero feed-forward term $\bD=0$, results to the convolution integral $y(t)=\int_{0}^{t}h_1(\sigma)u(t-\sigma)d\sigma$.
On the other hand, direct application of the Laplace transform $\cL(\cdot)$ in\cref{eq:LTI} with the complex frequency $s\in\IC$, will give
\begin{equation}\label{eq:LTIlap}
    \left\{\begin{aligned}
    \cL\left[\bE\dot{\bx}(t)\right]&=\cL\left[\bA\bx(t)\right]+\cL\left[\bB u(t)\right],\\
    \cL\left[y(t)\right]&=\cL\left[\bC\bx(t)\right]+\cL\left[\bD u(t)\right],
    \end{aligned}\right.\Rightarrow\left\{\begin{aligned}
    s\bE\bX(s)-\bE\bx_0&=\bA\bX(s)+\bB U(s),\\
    \bY(s)&=\bC\bX(s)+\bD U(s).
    \end{aligned}\right.
\end{equation}
Solving the algebraic equation in \cref{eq:LTIlap} w.r.t $\bX(s)$, and substituting to the output equation of the system in \cref{eq:LTIlap} as $Y(s)=\left[\bC(s\bE-\bA)^{-1}\bB+\bD\right]U(s)$, we conclude to the definition of the transfer function (1st GFRF in VS; thus, denoted $H_1$) as
\begin{equation}\label{eq:H1}
        H_1(s):=\frac{Y(s)}{U(s)}=\left[\bC(s\bE-\bA)^{-1}\bB+\bD\right],~~H_1:\IC\rightarrow\IC.
\end{equation}
The important first observation is that the transfer function \cref{eq:H1} involves all the matrices that define a linear system in state-space form over the frequency-domain, e.g., the imaginary-axis of the complex plane is with $s=j\omega,~\omega\in\IR$. The transfer function can also be considered an invariant measure and thus unique under different state-space realizations obtained under different coordinates. For the case where $\bD=0$, the Laplace transform of the impulse response $h_1(\sigma),~\sigma\in\IR_+$, is the transfer function $H_1(s),~s\in\IC$. For a concise representation of the following substantial quantities, we introduce the resolvent as
\begin{equation}\label{eq:resolvent}
    \bPhi(s)=(s\bE-\bA)^{-1}\in\IC^{n\times n}.
\end{equation}
\subsection{The Loewner framework for LTI systems}\label{sec:LF}
We start with the Loewner framework (LF) in the linear case utilizing similar notions from \cite{aa90,birkjour,morKarGA19a}. The LF is an interpolatory method that infers potentially reduced models whose transfer function matches/interpolates that of the underlying original system's (FOM) transfer function at selected interpolation points. To enforce this interpolation, and in the case of SISO systems, the following rational scalar interpolation problem has to be formulated. 

Consider the set of complex-valued data with $f$ a complex rational function as $\{\left(s_k,f(s_k)\right)\in\IC\times\IC:k=1,\ldots,2n)\}$. We partition these data into two disjoint subsets: $$\bS=[\underbrace{s_1,\ldots,s_n}_{\mu},\underbrace{s_{n+1},\ldots,s_{2n}}_{\lambda}],~\bF=[\underbrace{f(s_1),\ldots,f(s_n)}_{\IV},\underbrace{f(s_{n+1}),\ldots,f(s_{2n})}_{\IW}],$$
where $\mu_i=s_i$, $\lambda_i=s_{n+i}$, $v_{i}=f(s_i)$, $w_{i}=f(s_{n+i})$ for $i=1,\ldots,n$.\\[1mm]
The objective is to find $H(s)\in\IC$ such that:
\begin{equation}\label{eq:interpol}
    H(\mu_i)=v_{i},~i=1,\ldots,n,~\text{and}~H(\lambda_j)=w_{j},~j=1,\ldots,n.
\end{equation}
The {left data set} is denoted as: $\mathbb{M}=\left[\mu_1,\cdots,\mu_n\right]\in\IC^{1\times n}$, and $\IV=\left[v_1,\cdots,v_n\right]^{T}\in\IC^{n\times 1}$,
while the {right data set} as: $\mathbb{\bLambda}=\left[\lambda_1,\cdots,\lambda_n\right]^{T}\in\IC^{n\times 1}$, and $\IW=[w_1,\cdots, w_n]\in\IC^{1\times n}$.
The set of interpolation points can be determined by the problem or selected to achieve the given model reduction goals, e.g., approximation within a specific frequency bandwidth. 
\subsubsection{The Loewner matrices}\label{sec:LoewnerMatrices}
Given the row array of complex numbers $(\mu_j,v_j)$, $j=1,\ldots,{n}$, and the column array, $(\lambda_i,w_i)$, $i=1,\ldots,{n},$ (with $\lambda_i$ and the $\mu_j$ mutually distinct)
the associated {Loewner matrix} $\IL$ and the {shifted Loewner matrix} $\sIL$ are defined as:
\begin{equation}\label{eq:ComplexLoewnerMats}
\IL\!=\!\left[\!\begin{array}{ccc}
\frac{v_1-w_1}{\mu_1-\lambda_1} & 
\cdots &
\frac{v_1-w_{n}}{\mu_1-\lambda_{n}} \\
\vdots & \ddots & \vdots \\
\frac{v_n-w_1}{\mu_n-\lambda_1} & 
\cdots &
\frac{v_n-w_{n}}{\mu_n-\lambda_{n}} 
\end{array}\!\right]\!\in\!\IC^{n\times n},~
\sIL\!=\!\left[\!\begin{array}{ccc}
\frac{\mu_1v_1-\lambda_1w_1}{\mu_1-\lambda_1} & 
\cdots &
\frac{\mu_1v_1-\lambda_{n}w_{n}}{\mu_1-\lambda_{n}} \\
\vdots & \ddots & \vdots \\
\frac{\mu_n v_n-\lambda_1 w_1}{\mu_n-\lambda_1} & 
\cdots &
\frac{\mu_n v_n-\lambda_{n}w_{n}}{\mu_n-\lambda_{n}} 
\end{array}\!\right]\!\in\!\IC^{n\times n}.    
\end{equation}
\begin{definition}\label{def:McMillan}
If $g$ is a rational function, i.e., $g(s)=\frac{p(s)}{q(s)}$, for appropriate polynomials $p$, $q$, the McMillan degree or the complexity of $g$ is~$\mbox{deg}\,(g)=\max\{\mbox{deg}(p),\mbox{deg}(q)\}$.
\end{definition}

Now, if $w_i=g(\lambda_i)$, and $v_j=g(\mu_j)$, are 
\textit{samples} of a rational function $g$, the \textit{main property} 
of the Loewner matrices asserts in the following \cref{the:aa90} from \cite{aa90}.
\begin{theorem}\label{the:aa90}
Let $\IL$ be as above. If $n\geq {deg}\,g$, then $\,{rank}\, \IL = {\deg}\,g$.
\end{theorem}
 The result in \cref{the:aa90} indicates that the rank of the matrix $\IL$, which is constructed directly from measurements, encodes the complexity of the underlying rational function $g$ in terms of the McMillan degree \cref{def:McMillan}. Furthermore, the same result holds for
matrix-valued functions $g$ that emerge in the multi-input multi-output (MIMO) case. Construction of the rational interpolants associated with the revealed complexity degree follows.
\subsubsection{Construction of interpolants}
If the pencil\footnote{A linear matrix pencil, denoted by $(A, B)$, plays an important role in linear algebra, i.e., in finding the eigenvalues of $(A, B)$ numerically.} $(\IL_s,~\IL)$ is regular, then the quadruple $(\bE=-\IL,~\bA=-\IL_s,~\bB=\IV,~\bC=\IW)$, is a minimal realization of an interpolant for the data, i.e.,
\begin{equation}\label{eq:illinterpolant}
    H(s)=\IW(\IL_s-s\IL)^{-1}\IV.
\end{equation}
Otherwise, as shown in \cite{aa90}, the problem in \cref{eq:interpol} has a solution provided that $$\texttt{rank}\,\left[s\IL-\IL_s\right]=\texttt{rank}\,\left[\IL,~\IL_s\right]=\texttt{rank}\,\left[\begin{array}{c}\IL~~\IL_s\end{array}\right]^T={r},$$ for all $s\in\{\mu_i\}\cup\{\lambda_j\}$. 
Consider then the compact/thin singular value decomposition (SVD): $\left[\IL,~\IL_s\right]=\bY\widehat{\Sigma}_{ {r}}\tilde{\bX}^*,~\left[\begin{array}{c}\IL~~\IL_s\end{array}\right]^T = {\tilde\bY}\Sigma_{ {r}} \bX^*$, where $\widehat{\Sigma}_{ {r}}$, $\Sigma_{ {r}}$ $\in\IR^{{{r}}\times{r}}$,~
$\bY \in\IC^{n\times{r}}$,$\bX\in\IC^{n\times{r}}$,~$\tilde{\bY}\in\IC^{2n\times{r}}$,~$\tilde{\bX}\in\IC^{r\times{2n}}$. The order $r$ can be chosen as the numerical rank (as opposed to the exact rank) of the Loewner pencil. Performing the SVD, the minimal degree as the (numerical) rank of the Loewner pencil can be revealed. Therefore, a selection of the truncation degree $r$ comprises a trade-off between the accuracy and dimensionality of the derived reduced model. Suppose the neglected singular values are close to machine precision. In that case, the interpolation error is close to machine precision, and identification has been achieved (given sufficient data). Still, the degree of the "reduced" model, in that case, could be high after including all the non-zero numerical spectrum that the singular values encode, which can lead to an identified high-dimensional complex model. On the other hand, if a truncation degree $r$ is considered with the rest of the singular values as $\sigma_{r+k}\neq 0,~k=1,...$, some interpolation error is introduced, but the model can benefit from the low-dimensionality without affecting its accuracy significantly. The estimated ROM matrices of the inferred model will be denoted with hats, e.g., $\hat{\bA}$. The LF realizes the operator of a linear system at a specific random coordinate system. The transfer function stays invariant to any coordinate transformation, similar to the Markov parameters \cite{Isidori1973}. We introduce the following \cref{the:LoewnerRealization}.
\begin{theorem}\label{the:LoewnerRealization}
The quadruple $(\hat{\bA},\hat{\bB},\hat{\bC},\hat{\bE})$ of size 
~$ {r}\times  {r}$, ~$ {r}\times  {1}$, ~$ {1}\times r$, ~$r\times  {r}$,~ given by:
\begin{equation*} \label{redundant}
\hat{\bE}  = -\bY^T\IL \bX ,~~\hat{\bA}  = -\bY^T\IL_s \bX, ~~\hat{\bB}  = \bY^T\IV,~~
\hat{\bC}  = \IW \bX ,
\end{equation*}
is a descriptor realization\footnote{A descriptor realization is the quadruple of the matrices that define a linear system $(\bA,\bB,\bC,\bE)$. Usually, we refer to such systems as descriptor systems when $\bE$ is not the identity matrix or is non-invertible (differential algebraic equations). The Loewner framework provides a descriptor realization. Alike the pencil $(\bA,\bE)$, the Loewner pencil is $(\IL_s,\IL)$, and both contain the same eigenvalues.} of an (approximate) interpolant of the data with McMillan degree $r=rank(\IL)$, where
\begin{equation}
    \hat{H}(s)=\hat{\bC}(s\hat{\bE}-\hat{\bA})^{-1}\hat{\bB}.
\end{equation}
\end{theorem}

For more details on the construction/identification of linear time-invariant systems with the LF, we refer the reader to \cite{ABG20,birkjour,morKarGA19a} where both the SISO and MIMO cases are addressed together with analysis on technical aspects such as distribution of interpolation points, splitting of measurements, construction of real-valued models, etc.) together with concise algorithms.
\subsubsection{An introductory linear example}\label{sec:linearintro}
Consider the following SISO minimal linear system of dimension $n=2$ with non-zero initial condition
\begin{equation}\label{eq:linsys}
    \left\{\begin{aligned}
        \dot{\bx}(t)&=\left[\begin{array}{cc}
           -1  & 0\\
           0 & -2
        \end{array}\right]\bx(t)+\left[\begin{array}{c}
             1  \\
             1 
        \end{array}\right]u(t),\\
        y(t)&=\left[\begin{array}{cc}
            1 & 1 \end{array}\right]\bx(t),~y_0=y(0)=\bC\bx(0)=\left[\begin{array}{cc}
                1 & 1 \\
            \end{array}\right]\left[\begin{array}{c}
                 0.5  \\
                 0 
            \end{array}\right]=0.5,~t\geq 0.
    \end{aligned}\right.
\end{equation}
The task is to infer the system in \cref{eq:linsys} from time-domain input-output measurements obtained under harmonic excitation.

\begin{figure}[h!]
\centering
    \includegraphics[scale=0.14]{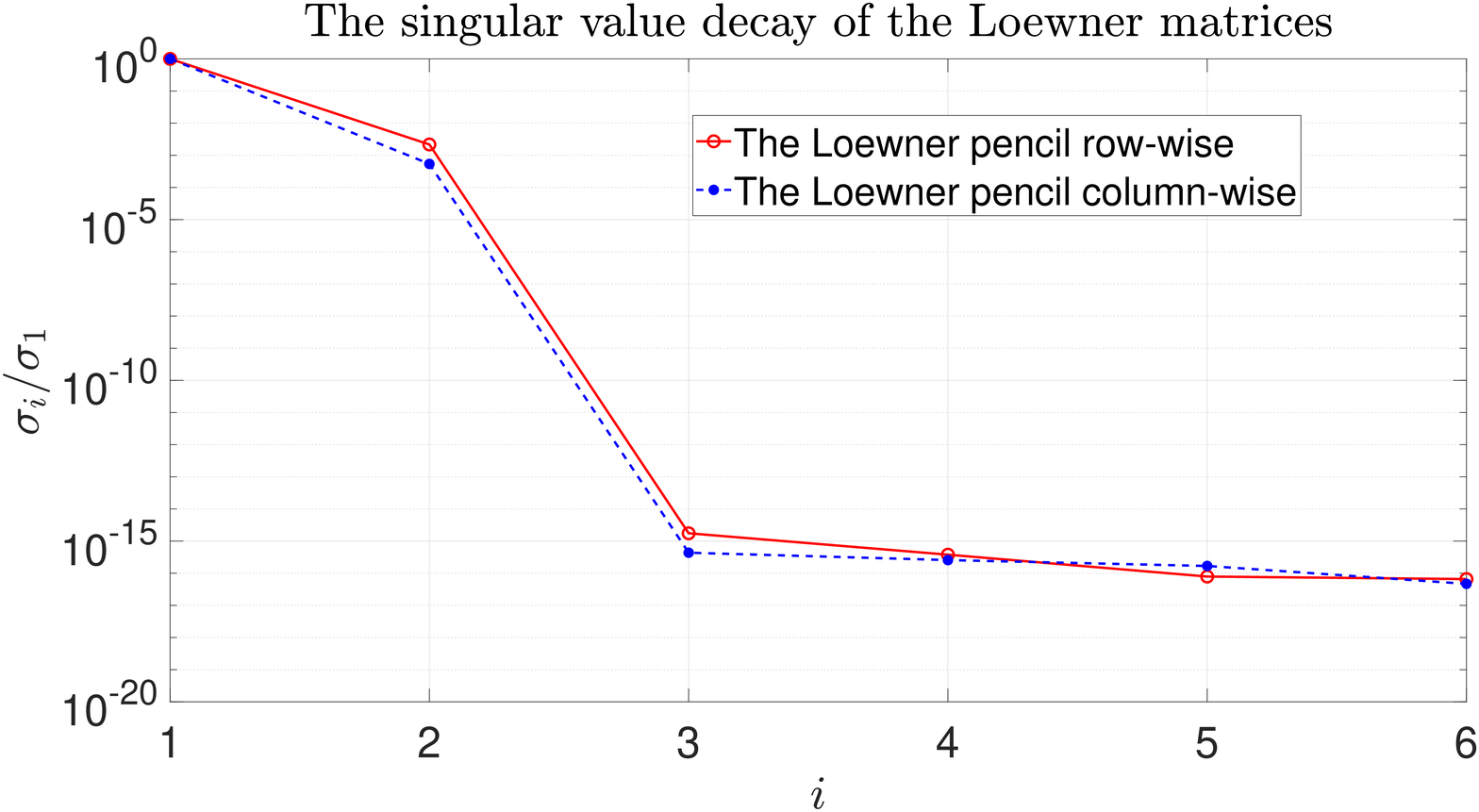}
     \includegraphics[scale=0.14]{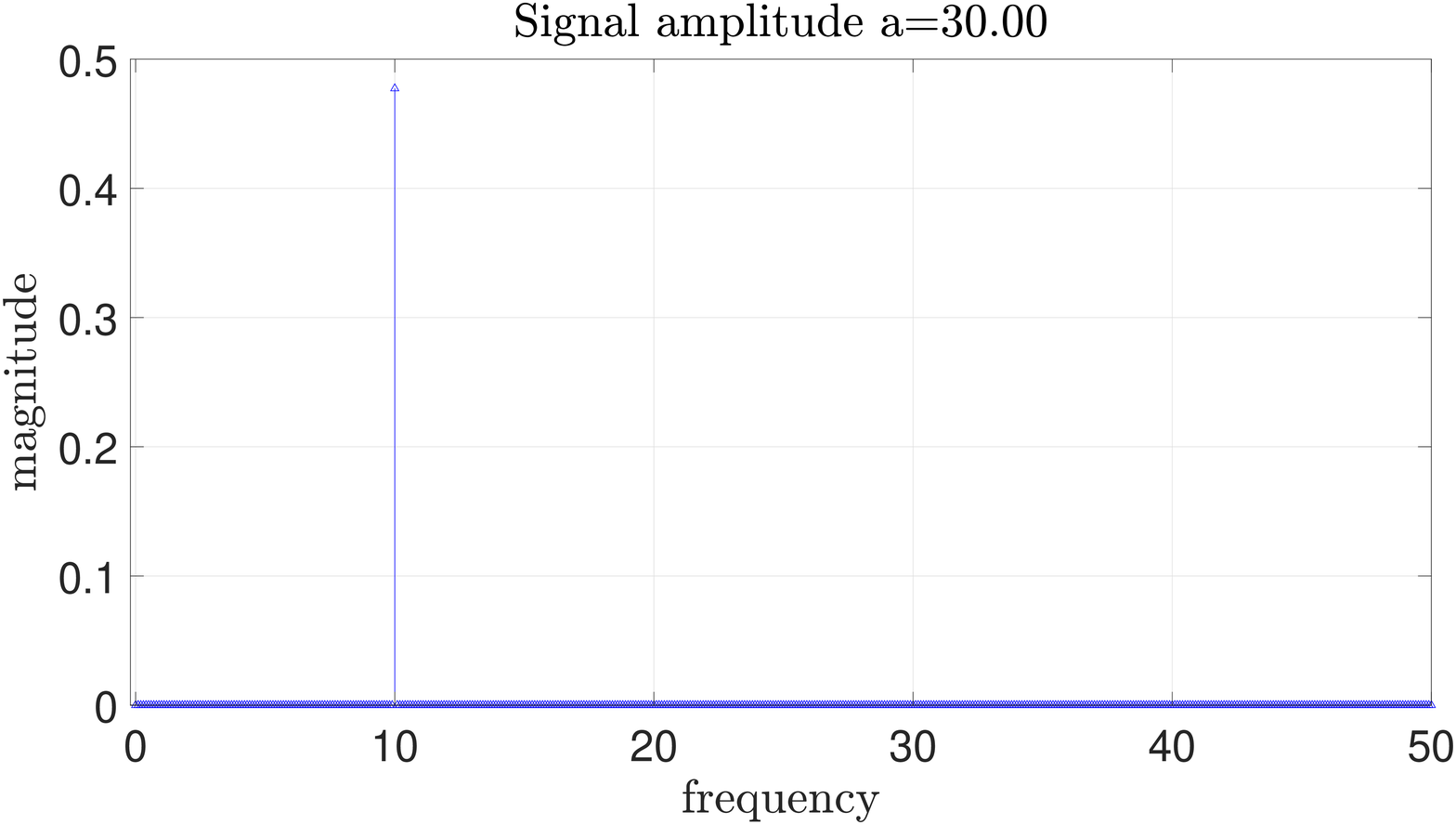}
    \caption{Left: The decay of singular values of Loewner matrices that exploit the minimal degree $r=2$. Right: Output evolution $y(t)$ with input as $u(t)=30\cos(2\pi\cdot 10 \cdot t)$ where $\alpha=30,~\omega_1=10~\text{(Hz)}$.}
    \label{fig:scalarlinear1}
\end{figure}

\paragraph{Data acquisition in the linear case} We excite the system in \cref{eq:linsys} with $u(t)=\alpha\cos(2\pi\omega_1 t)$, and we get the following power spectrum in \cref{fig:scalarlinear1}(right). After applying fast Fourier transform (FFT) to the quasi-steady state profile, the single-sided power spectrum provides the complex phase from which we estimate the transfer function (1st Volterra kernel) as $H_1(j2\pi\omega_1)=2P(2\pi j\omega_1)/\alpha\in\IC$. We repeat the procedure until we obtain enough estimations for the $H_1$. As the underlying system has a minimal dimension $n=2$, four $(4)$ measurements are sufficient for the Loewner framework to identify the system. Here, to illustrate how the LF method performs data-driven reduction, we take more measurements and let the SVD algorithm decide on the minimal truncation degree. Suppose we collect more measurements, i.e., six $(6)$ pairs of complex points, thus, $\mathcal{D}:=\left\{\left(2\pi j\omega_1,H_1(2\pi j\omega_1)\right),\left(2\pi j\omega_2,H_1(2\pi j\omega_2)\right),\ldots,\left(2\pi j\omega_6,H_1(2\pi j\omega_6)\right)\right\}.$

We split the data set $\mathcal{D}$ into two disjoint sets according to \cref{tab:datasplit}, and we form the complex Loewner matrices of dimension three $(3)$ as in \cref{eq:ComplexLoewnerMats}.
\begin{table}[h!]
    \centering
    \begin{tabular}{cc|cc}
       Left data &  & Right data & \\\hline
        $\mu_1=2\pi j\omega_1$ & $\upsilon_1=H_1(2\pi j\omega_1)$ & $\lambda_1=2\pi j\omega_4$ & $w_1=H_1(2\pi j\omega_4)$ \\
         $\mu_2=2\pi j\omega_2$ & $\upsilon_2=H_1(2\pi j\omega_2)$ & $\lambda_2=2\pi j\omega_5$ & $w_2=H_1(2\pi j\omega_5)$ \\
          $\mu_3=2\pi j\omega_3$ & $\upsilon_3=H_1(2\pi j\omega_3)$ & $\lambda_3=2\pi j\omega_6$ & $w_3=H_1(2\pi j\omega_6)$ \\
    \end{tabular}
    \caption{Complex data and a splitting scheme for the Loewner framework.}
    \label{tab:datasplit}
\end{table}

The underlying system is real; thus, for any measurement $s\in\IC$, it holds $H(s)=-\bar{H}(\bar{s})$ where $\bar{s}$ denotes the conjugate of $s$. With this real symmetric property, we can transform complex data to real as in \cite{morKarGA19a}.  

The conversion from complex to real data holds for both the Loewner matrices. More details on deriving the equivalent real Loewner model can be found in \cite{morKarGA19a}. 
In our example, we proceed by including the six conjugate pairs of measurements to produce a real Loewner model of dimension six $(6)$. To obtain the interpolant in \cref{eq:illinterpolant}, we first perform the SVD of the Loewner matrices. Suppose the measurements are $\omega_i=5i,~i=1,\ldots,6$.

In \cref{fig:scalarlinear1}(left), the singular value decay gives the minimal degree of the underlying linear model as $r=2$. Thus, we use the left and right singular matrices to project the singular Loewner model of dimension $6$ to $r=2$. By using the \cref{the:LoewnerRealization}, the reduced linear system has the following matrices with $5$-digits precision
\begin{equation}\label{eq:scalarlinearreduced}
\hat{\bA}_r=\left[\begin{array}{cc} 5.2775 & -0.30554\\ 149.52 & -8.2775 \end{array}\right],~\hat{\bB}_r=\left[\begin{array}{c} -63.79\\ 4.4346 \end{array}\right],~\hat{\bC}_r=\left[\begin{array}{cc} -0.031254 & 0.0014167 \end{array}\right].
\end{equation}
The reduced system is
\begin{equation}\label{eq:identified}
    \left\{\begin{aligned}
        \dot{\hat{\bx}}(t)&=\hat{\bA}\hat{\bx}(t)+\hat{\bB}u(t),\\
        \hat{y}(t)&=\hat{\bC}\hat{\bx}(t),~\hat{\bx}(0)=\hat{\bx}_0=?
    \end{aligned}\right.
\end{equation}
The linear matrix in \cref{eq:scalarlinearreduced} has the following eigenvalues $\texttt{eig}(\hat{\bA})=(-1,-2)$ that match the original in \eqref{eq:linsys} with machine precision. To compare the two systems (original and identified), we need to align the two coordinate systems to make any use of the original observable initial conditions $y_0$. To infer $\hat{\bx}_0$, transient response data is necessary because they encode information on the system's initial conditions, in contrast to steady-state system response. Thus, we simulate the corresponding autonomous system with $u(t)=0$. The transient response of the original system until the dynamics converge to the zero equilibrium is $y(t)=\bC e^{At}\bx_0$. We must enforce $\hat{y}(t)=y(t),\forall t$. To infer the initial conditions $\hat{\bx}_0$, we need to specify unknowns equal to the discovered dimension $r$ and can be easily generalized for any $r\in\IZ_+$. In our case, where $r=2$, we write $\hat{\bx}_0=\left[\begin{array}{cc}
    \hat{x}_{01} & \hat{x}_{02} \\ 
\end{array}\right]^T$. To specify $2$ degrees of freedom, we need at least $2$ equations. We enforce identical transient phases in the absence of the controller $u(t)=0$ as
\begin{equation}
    \begin{aligned}
        \hat{y}(t)&=y(t)\Leftrightarrow\hat{\bC}e^{\hat{\bA}t}\hat{\bx}_0=y(t)\Leftrightarrow\hat{\bC}e^{\hat{\bA}t}\left[\begin{array}{cc}
             \hat{x}_{01}  \\
             \hat{x}_{02}             
        \end{array}\right]=y(t),~\text{by using}~t_1,t_2\Rightarrow\\
        &\underbrace{\left[\begin{array}{c}
             \hat{\bC}e^{\hat{\bA}t_1}  \\
             \hat{\bC}e^{\hat{\bA}t_2} 
        \end{array}\right]}_{\in\IR^{2\times 2}}\left[\begin{array}{cc}
             \hat{x}_{01}  \\
             \hat{x}_{02}             
        \end{array}\right]=\underbrace{\left[\begin{array}{cc}
             y(t_1)  \\
             y(t_2)             
    \end{array}\right]}_{\text{measurements}}\Rightarrow\left[\begin{array}{cc}
             \hat{x}_{01}  \\
             \hat{x}_{02}             
        \end{array}\right]=\left[\begin{array}{c}
             \hat{\bC}e^{\hat{\bA}t_1}  \\
             \hat{\bC}e^{\hat{\bA}t_2} 
        \end{array}\right]^{-1}\left[\begin{array}{cc}
             y(t_1)  \\
             y(t_2)             
        \end{array}\right].
    \end{aligned}
\end{equation}
For our example, we use $t_1=0,~t_2=1$. Thus, the above inversion is well-conditioned, and we compute the initial condition to the hatted system with five-digit precision as
\begin{equation}\label{eq:initialconditionhat}
    \hat{\bx}_0=\left[\begin{array}{c} -232.79\\ -4782.8 \end{array}\right].
\end{equation}
The original system in \cref{eq:linsys} and the identified \cref{eq:identified} together with the inferred initial conditions \cref{eq:initialconditionhat} are equivalent in all dynamical regimes (transient and steady state). They differ only by one similarity transform (\cref{app:align}). To illustrate this argument, we construct such a transformation from \cref{app:align}, that the system's coordinates are aligned (and thus identical to) the original system's operator, it remains
\begin{equation}
    \Psi^{-1}\hat{\bx}_0= \Psi^{-1}\left[\begin{array}{c} -232.79\\ -4782.8 \end{array}\right]=\left[\begin{array}{c} 0.5\\ 0 \end{array}\right].
\end{equation}

Note that in the linear case, with only one multi-harmonic excitation, an injenctive input-output harmonic map with enough measurements to realize the transfer function from a single input within the LF method similar to \cite{PEHERSTORFER2016196} can be obtained. In contrast, for the nonlinear case, due to harmonic intermodulation, the harmonic indexing and kernel separation make the task of estimating GFRFs with a multi-harmonic excitation challenging. Finally, as we have analyzed the mechanism in dealing with linear system identification via the LF, which will play a significant role in identifying the minimal linear subsystem from the nonlinear class under consideration, nonlinear analysis is ready to emerge.

\subsection{Nonlinear systems theory with the Volterra series (VS) representation}\label{sec:Volterra}
In this section, we introduced the VS representation of nonlinear systems that admit such polynomial series expansion, i.e., weakly nonlinear, where the energy comes from the controller and the dynamics evolve after possible bifurcation around achievable local stable equilibrium points. For the SISO case, i.e., $u\rightarrow\boxed{\Sigma_{\text{NL}}}\rightarrow y$, the mathematical notion of the VS contains an infinite series expansion based on multi-convolution operations as
\begin{equation}\label{eq:Volterra}
y(t)=\sum_{n=1}^{\infty}y_n(t),~y_n(t)=\int_{0}^{\infty}\cdots\int_{0}^{\infty} h_n(\tau_1,\ldots,\tau_n)\prod_{i=1}^{n}u(t-\tau_i)d\tau_i,
\end{equation}
where $h_n(\tau_1,\ldots,\tau_n)$ is a real-valued function of $\tau_1,\ldots,\tau_n$ known as the $n$th-order time-domain Volterra kernel. After a multivariate Laplace transform to the time-domain kernels $h_n(\tau_1,\ldots,\tau_n)$, the $n$th-order generalized frequency response function (GFRF) is defined as
\begin{equation}\label{eq:kernels1}
H_n(s_1,\ldots,s_n)=\int_{0}^{\infty}\cdots\int_{0}^{\infty}h_n(\tau_1,\ldots,\tau_n)e^{-\sum_{i=1}^{n}s_i\tau_i}d\tau_1\cdots d\tau_n.
\end{equation}

The above mathematical formulation is general and covers many nonlinear systems. One way to derive closed forms of GFRFs is to assume a physics-based structure of the underlying system. In this direction, quadratic systems are of particular interest. This is due to the variety of applications where the analytical state nonlinearity is quadratic (e.g., Navier Stokes, Burgers', etc.), as well as the possibility of embedding any system with analytical nonlinearity to a quadratic structure through lifting techniques \cite{Gu11,Gosea22mdpi}. As detailed in \cite{Gu11}, lifting strategies can result in polynomial ODEs systems of quadratic order or in systems of differential-algebraic equations (DAEs) where the non-invertible $\bE$ in the latter case makes the problem quite challenging even when state access measurements (snapshots) can be accessed \cite{KHODABAKHSHI2022114296}. We proceed with general nonlinear systems that result after possible lifting embedding to quadratic ODEs.

\subsection{The quadratic dynamical system}\label{sec:Quadratic}
We continue our analysis with the general state-space representation of the quadratic dynamical system and for the SISO case as
\begin{equation}\label{eq:qsys}
    \left\{\begin{aligned}
    \dot{\bx}(t)&=\bA\bx(t)+\bQ(\bx(t)\otimes\bx(t))+\bB u(t),\\
    y(t)&=\bC\bx(t),~\bx(0)=\bx_0=\textbf{0},~t\geq 0,
    \end{aligned}\right.
\end{equation}
where the state dimension is $n$, and the matrices are: $\bA\in\IR^{n\times n}$, $\bQ\in\IR^{n\times n^2},~\bB,\bC^T\in\IR^{n\times 1}$. The term $\bQ(\bx(t)\otimes\bx(t))$ denotes the product of the quadratic matrix with the Kronecker product. The Kronecker product $\otimes$ is defined as in the following simple case $\left[\begin{array}{cc}
    x_1 & x_2 \\
\end{array}\right]\otimes\left[\begin{array}{cc}
    x_1 & x_2 \\
\end{array}\right]=\left[\begin{array}{cccc}
    x_1^2 & x_1x_2 & x_2x_1 & x_2^2 \\
\end{array}\right]$. Due to the commutative property, the matrix $\bQ$ denotes the Hessian of the right-hand side in \cref{eq:qsys} and exhibits a symmetric structure. In particular, for any two arbitrary vectors $\bu,\bv\in\IR^n$, we can always ensure that it holds
\begin{equation}\label{eq:Qsym}
    \bQ(\bu\otimes\bv)=\bQ(\bv\otimes\bu).
\end{equation}

We now focus on systems with quadratic nonlinearities \cref{eq:qsys}. Different ways of deriving generalized frequency response functions (GFRFs) exist. We want to obtain closed formulations of GFRFs that can be measured under harmonic excitation, interpreted in the frequency domain, and then associated with the quadratic operator matrices. To be more mathematically precise, suppose a two-tone excitation signal like $u(t)=\cos(\omega_1 t)+\cos(\omega_2 t)$; thus, the measurement should correlate with $H_2(\omega_1,\omega_2)$. If we interchange the frequencies, the input stays the same, and the measurement correlates with $H_2(\omega_2,\omega_1)$. Consequently, these two measurements cannot differ, which explains precisely the symmetric property. Therefore, our first aim is to derive the symmetric GFRFs for the quadratic case that can be processed from the time domain to the frequency domain. The second aim is to use these measurements to infer the matrices $(\bA,~\bQ,~\bB,~\bC)$ that define the assumed quadratic system class.  

\subsubsection{Deriving higher-order transfer functions for the quadratic system}\label{sec:2.5}
As detailed in \cite{Ru82}, different ways exist to extract different types of GFRFs. One way is the variational approach, where the structure of the triangular kernels (or regular kernels) can be revealed through Picard iterations. In particular, regular GFRFs can be derived after shifting the frequency domain of triangular GFRFs. Despite the ease of using intrusive models, the regular Volterra kernels cannot be measured directly from the time domain. Therefore, we continue deriving GFRFs, namely with the growing exponential approach (i.e., the probing method) for treating the issue of kernel estimation from time domain harmonic data. By probing the system and after processing the quasi-steady state time evolution in the frequency domain via the Fourier (a special case of Laplace over the imaginary axis) transform, the time domain signal is decomposed to harmonics that scale and shift w.r.t. the symmetric GFRFs. Methods for estimating these symmetric GFRFs (e.g., kernel separation) were introduced in \cite{Boyd1983MeasuringVK,morKarGA19,Xpar2006}.
\paragraph{\textbf{The Probing method}} It was shown by Rugh \cite{Ru82} and Billings \cite{Billings2013NonlinearSI} that for nonlinear systems which are described by the VS \cref{eq:Volterra} and excited by a combination of complex exponentials $u(t)=\sum_{i=1}^{R}e^{s_i t},\quad1\leq R$, the output response can be written as

\begin{equation}
\label{eq:output}
\begin{aligned}
y(t)&=\sum_{n=1}^{\infty}\sum_{i_1=1}^{R}\cdots\sum_{i_n=1}^{R}{H}_{n}(s_{i_1},\ldots,s_{i_n})e^{(s_{i_1}+\cdots+s_{i_n})t}\\
&=\sum_{n=1}^{\infty}\sum_{m(n)}\tilde{H}_{m_{1}(n)\cdots m_{R}(n)}(s_{1},\ldots,s_{R})e^{(m_{1}(n)s_1+\cdots+m_{R}(n)s_{R})t},
\end{aligned}
\end{equation}
where $\sum_{m(n)}$ indicates an $R$-fold sum over all integer indices $m_{1}(n),\ldots,m_{R}(n)$ such that $0\leq m_{i}(n)\leq n,~m_{1}(n)+\cdots+m_{R}(n)=n$, and
\begin{equation}\label{eq:kernels2}
\tilde{H}_{m_{1}(n)\cdots m_{R}(n)}(s_1,\ldots,s_{R})=\frac{n!}{m_{1}(n)!\cdots m_{R}(n!)}H_{n}(\underbrace{s_{1},\ldots,s_{1}}_{m_{1}(n)},\ldots,\underbrace{s_{R},\ldots,s_{R}}_{m_{R}(n)}).
\end{equation}

 Note that $\tilde{H}_n$ is the weighted GFRF, corresponding to $H_n$ as in \cref{eq:kernels2}; The former scales with the factor $\frac{n!}{m_{1}(n)!\cdots m_{R}(n!)}$. Note also that different input amplitudes can be considered, as in \cite{VoltMacrXpar,Billings2013NonlinearSI} where the amplitude shift allows the kernel separation. 

 Following the derivations in \cref{app:GFRFs}, we obtain the set of the first three symmetric GFRFs 
 \begin{equation}\label{eq:symGFRFs}
     \begin{aligned}
         H_1(s_1)&=\bC\underbrace{\bPhi(s_1)\overbrace{\bB}^{\bR_1}}_{\bG_1(s_1)},\\
         H_2(s_1,s_2)&=\bC\underbrace{\frac{1}{2}\bPhi(s_1+s_2)\bQ\overbrace{\left[\bG_1(s_1)\otimes\bG_1(s_2)+\bG_1(s_2)\otimes\bG_1(s_1)\right]}^{\bR_2(s_1,s_2)}}_{\bG_2(s_1,s_2)},\\
         H_3(s_1,s_2,s_3)&=\bC\frac{1}{6}\bPhi(s_1+s_2+s_3)\bQ\overbrace{\left[+~\bG_1(s_1)\otimes\bG_2(s_2,s_3)+\bG_2(s_2,s_3)\otimes\bG_1(s_1)+\right.}^{\bR_3(s_1,s_2,s_3)}\\
         &~~~~~~~~~~~~~~~~~~~~~~~~~~~~~~~~~+\left.\bG_1(s_2)\otimes\bG_2(s_1,s_3)+\bG_2(s_1,s_3)\otimes\bG_1(s_2)+\right.\\
         &~~~~~~~~~~~~~~~~~~~~~~~~~~~~~~~~~+\left.\bG_1(s_3)\otimes\bG_2(s_1,s_2)+\bG_2(s_1,s_2)\otimes\bG_1(s_3)\right].
     \end{aligned}
 \end{equation}
 
 At this point, we illustrate some properties the symmetric GFRFs in \cref{eq:symGFRFs} inherit. 
\begin{itemize}
    \item \textbf{Symmetry}: As it is evident, any permutation of the set $(s_1,s_2,\ldots,s_n)$ will result to the same evaluation of the $H_n(s_1,s_2,\ldots,s_n)$ and $\bG_n(s_1,s_2,\ldots,s_n)$.
    \item \textbf{Decompositions}: In control theory, the reachability and observability conditions are connected theoretically with the minimality degree \cite{ACA05,ABG20}. In simple terms, reachability signifies the ability to reach a state by applying suitable input to the system.  On the other hand, observability signifies the ability to observe some system state. Introducing the general reachability $\cR$ and observability $\cO$ matrices similar to \cite{ABG20}, and for the symmetric GFRFs case, a more concise representation of the GFRFs \cref{eq:symGFRFs} can be derived. As we have introduced the generalized reachability matrices in \cref{app:GFRFs} and in \eqref{eq:symGFRFs} as $\bR_1,~\bR_2,~\bR_3$, the corresponding general observability matrices remain
    \begin{equation}
        \begin{aligned}
        \bO_1(s_1)&=\bC\bPhi(s_1),\\
        \bO_2(s_1,s_2)&=\frac{1}{2}\bC\bPhi(s_1+s_2),\\
        \bO_3(s_1,s_2,s_3)&=\frac{1}{6}\bC\bPhi(s_1+s_2+s_3).\\
        \end{aligned}
    \end{equation}

Based on the above observations and notions, we can better represent the input to state symmetric GFRFs by exploiting their structure. To distinguish the position of the quadratic matrix in the 3rd GFRFs, the superscripts $(\cdot)^\ell$-left and $(\cdot)^r$-right are introduced.
\begin{equation}
    \begin{aligned}
    \bG_1(s_1)&=\bPhi(s_1)\bR_1,\\
    \bG_2(s_1,s_2,\bQ)&=\frac{1}{2}\bPhi(s_1,s_2)\bQ\bR_2(s_1,s_2),\\
    \bG_3(s_1,s_2,s_3,\bQ^{\ell},\bQ^{r})&=\frac{1}{6}\bPhi(s_1,s_2,s_3)\bQ^{\ell}\bR_3(s_1,s_2,s_3,\bQ^{r})
    \end{aligned}
\end{equation}
and for the input to output GFRFs as
\begin{equation}\label{eq:ioTFscript}
    \begin{aligned}
    H_1(s_1)&=\bO_1(s_1)\bR_1,\\
    H_2(s_1,s_2,\bQ)&=\bO_2(s_1,s_2)\bQ\bR_2(s_1,s_2),\\
    H_3^{(ij)}(s_1,s_2,s_3,\bQ^{(i)},\bQ^{(j)})&=\bO_3(s_1,s_2,s_3)\bQ^{(i)}\bR_3(s_1,s_2,s_3,\bQ^{(j)}).
    \end{aligned}
\end{equation}
\item \textbf{The generalized reachability matrix $\bR_3(\bQ)$ is linear w.r.t the quadratic operator $\bQ$}. Assume $\lambda_1,~\lambda_2\in\IR$ and $\bQ_1,~\bQ_2\in\IR^{n\times n^2}$. Then, it holds 
\begin{itemize}
    \item \textbf{Linear property:} $\bR_{3}(\lambda_1\bQ_1+\lambda_2\bQ_2)=\lambda_1\bR_{3}(\bQ_1)+\lambda_2\bR_{3}(\bQ_2)$.
    \begin{proof} By neglecting the similar-structured terms (s.s.t), we can prove the following:
    \begin{equation*}
        \begin{aligned}
        &\bR_{3}(s_1,s_2,s_3,\lambda_1\bQ_1+\lambda_2\bQ_2)=\bG_1(s_1)\otimes\bG_2(s_2,s_3,\lambda_1\bQ_1+\lambda_2\bQ_2)+s.s.t.\\
        &=\bG_1(s_1)\otimes\frac{1}{2}\bPhi(s_1,s_2)(\lambda_1\bQ_1+\lambda_2\bQ_2)\bR_2(s_1,s_2)+s.s.t.\\
        &=\bG_1(s_1)\otimes\frac{1}{2}\bPhi(s_1,s_2)\lambda_1\bQ_1\bR_2(s_1,s_2)+\bG_1(s_1)\otimes\frac{1}{2}\bPhi(s_1,s_2)\lambda_2\bQ_2\bR_2(s_1,s_2)+s.s.t.\\
        &=\lambda_1\bG_1(s_1)\otimes\bG_2(s_2,s_3,\bQ_1)+\lambda_2\bG_1(s_1)\otimes\bG_2(s_2,s_3,\bQ_2)+s.s.t.\\
        &=\lambda_1\bR_{3}(s_1,s_2,s_3,\bQ_1)+\lambda_2\bR_{3}(s_1,s_2,s_3,\bQ_2).
        \end{aligned}
    \end{equation*}
    \end{proof}
\end{itemize}
\end{itemize}

Starting from the quadratic dynamical system in \cref{eq:qsys}, we have derived all the quantities of interest (with their corresponding properties) for setting up our method in what follows. Equivalent descriptions between the time and frequency domain representations have been addressed for this problem using the VS. Closing this section, as we deal with a data-driven method, we proceed in explaining the data acquisition process starting from a solid, consistent theoretical background that was exploited in the previous derivations and in the Volterra framework for representing nonlinear systems.

\subsubsection{Data acquisition with harmonic probing}\label{sec:dataaquisition}
For a given multi-harmonic input in the time domain with driving frequencies $(s_1,s_2,\ldots,s_n)\in\IC^n$, the idea is to estimate the GFRFs, i.e., $H_n(s_1,s_2,\ldots,s_n)$ from the processed frequency power spectrum. These GFRFs encode system-invariant information that can lead to identification or, in general, inference since they are directly associated with the operator matrices from the underlying quadratic system. When the excitation signal is $u(t)=\sum_{i=1}^Re^{s_i t}$, the output can be decomposed to a VS \cref{eq:output}. We showcase the following representative cases:
\begin{itemize}
    \item Case 1: A single complex harmonic $$u(t)=e^{2\pi j\omega t}\rightarrow\boxed{\Sigma_{\text{NL}}}\rightarrow \sum_{n=0}^{\infty}H_n(\underbrace{\omega,\ldots,\omega}_{n~\text{times}})e^{2\pi nj\omega t},$$
where the harmonic indexing and GFRFs (kernels) separation are direct.
    \item Case 2: A double complex harmonic with $\omega_1,\omega_2>0$ 
    $$
    u(t)=e^{2\pi j\omega_1 t}+e^{2\pi j\omega_2 t}\rightarrow\boxed{\Sigma_{\text{NL}}}\rightarrow \sum_{n=0}^{\infty}H_{n}(\underbrace{\omega_1,\ldots,\omega_1}_{m_1(n)~\text{times}},\underbrace{\omega_2,\ldots,\omega_2}_{m_2(n)~\text{times}})e^{2\pi j(m_1(n)\omega_1+m_2(n)\omega_2)t},
    $$
    where the propagating harmonics are infinite and precisely the integer multiples of $m_1(n)\omega_1+m_2(n)\omega_2,~m_1(n),m_2(n)\in\IZ_+$ with $m_1(n)+m_2(n)=n$. For instance, if $n=2$, the combinations are $(2,0),(1,1),(0,2)$. These correspond to the separable kernels $H_2(2\omega_1)$, $H_2(\omega_1,\omega_2)$, $H_2(2\omega_2)$ that can be inferred directly from the power spectrum. Of course, the wise choice of $\omega_1,~\omega_2$ can make the kernel separation exact.
     \item Case 3: A single real harmonic $2a\cos(2\pi\omega t)$ where  from Euler's formula with $\omega>0$, we have
    $$
    u(t)=ae^{2\pi j\omega t}+ae^{-2\pi j\omega t}\rightarrow\boxed{\Sigma_{\text{NL}}}\rightarrow \sum_{n=0}^{\infty}a^nH_{n}(\underbrace{\omega,\ldots,\omega}_{m_1(n)~\text{times}},\underbrace{-\omega,\ldots,-\omega}_{m_2(n)~\text{times}})e^{2\pi j(m_1(n)-m_2(n))\omega)t},
    $$
    where the propagating harmonics are infinite and precisely with values of integer multiples of $(m_1(n)-m_2(n))\omega,~m_1(n),m_2(n)\in\IZ_+$ with $m_1(n)+m_2(n)=n$. Immediately, each harmonic is a superposition of infinite kernels as for $\omega$ holds 
    \begin{equation}\label{eq:harmonic}
    Y_{1}(\omega)=aH_{1}(\omega)+a^3w_3H_3(\underbrace{-\omega,\omega,\omega}_{2\omega-\omega=\omega})+a^5w_5H_5(\underbrace{-\omega,+\omega,-\omega,\omega,\omega}_{3\omega-2\omega=\omega})+...
\end{equation}
\end{itemize}

In practice, one could excite the system with a multi-tone input of low amplitude such that higher harmonic intermodulation is negligible. When the excitation activates more harmonics, kernel separation with the known techniques from \cite{BoydChua85,Boyd1983MeasuringVK,morKarGA21a} are available. Equation \eqref{eq:harmonic} suggests kernel separation by sweeping the amplitude, which can be done up to a degree with multi-tone input. The polynomial structure in the amplitude $a$ easily can lead to ill-conditioning due to the Vandermonde construction, but such an issue can be addressed from \cite{ArnoldiNaka}. Moreover, to this direction of harmonic detection and kernel separation, machinery such \cite{Xpar2006} can make this process applicable even for more complex cases based on small signal multi-harmonic perturbations that lead to thorough linearization in handling efficiently the harmonic distortion.
\subsection{Method for quadratic modeling from i/o time domain data}\label{sec:method}
Next, we introduce the proposed method for computing quadratic state-space models from the first three symmetric GFRFs \cref{eq:symGFRFs}, which can be measured from time-domain harmonic excitation in two steps. The first step infers the minimal linear subsystem with the LF, and the second step employs linear and nonlinear optimization schemes to infer the quadratic matrix. Combining these, a realization of the operator that defines the quadratic system can be achieved.
\subsubsection{Identification of the linear subsystem with the Loewner framework}\label{sec:linerminimal}
Using measurements of the 1st harmonic, we can identify the minimal linear subsystem of order $r\leq n$, as $(\hat{\bA},~\hat{\bB},~\hat{\bC})$ ($\hat{\bE}$ invertible) with the Loewner framework as illustrated in \cref{sec:linearintro}. Further, using the identified linear subsystem, we can formulate optimization problems where estimations of the quadratic operator can be achieved after incorporating higher-order harmonics associated with the higher-order GFRFs, i.e., $H_2,~H_3$. We acquire and solve these optimization problems in two steps: solving an under-determined linear optimization problem that obeys infinite parameterized solutions due to non-empty null space within a least-squares setting and solving a non-linear optimization problem with the Newton method for fixing the remaining parameters.

\subsubsection{Estimation of the quadratic operator from the 2nd kernel}\label{sec:estimationQfromH2}
Identification of the minimal linear subsystem $(\hat{\bA},\hat{\bB},\hat{\bC})$ of order $r$ as described in \cref{sec:linerminimal} allows the construction of the reduced resolvent $\hat{\bPhi}(s)=(s\hat{\bI}-\hat{\bA})^{-1}\in\IC^{r\times r}$, and the 2nd GFRFs with the unknown operator $\hat{\bQ}$ can be written as:
\begin{equation}
\begin{aligned}
\hat{H}_2(s_1,s_2)&=\underbrace{\frac{1}{2}\hat{\bC}\hat{\bPhi}(s_1+s_2)}_{\hat{\bO}_2(s_1,s_2)}\hat{\bQ}\underbrace{\left[(\hat{\bPhi}(s_1)\hat{\bB})\otimes(\hat{\bPhi}(s_2)\hat{\bB})+(\hat{\bPhi}(s_2)\hat{\bB})\otimes(\hat{\bPhi}(s_1)\hat{\bB})\right]}_{\hat{\bR}_2(s_1,s_2)}=\\
&=\hat{\bO}_2(s_1,s_2)\hat{\bQ}\hat{\bR}_2(s_1,s_2)=\left(\hat{\bO}_2(s_1,s_2)\otimes\hat{\bR}_2^{T}(s_1,s_2)\right)\texttt{vec}(\hat{\bQ}).
\end{aligned}
\end{equation}
The way of estimating the quadratic operator $\hat{\bQ}$ comes after enforcing interpolation with the 2nd harmonic (2nd kernel) over a 2D grid of selected measurements $\left(s_1^{(k)},s_2^{(k)}\right)$. Thus, we enforce
\begin{equation}
    \underbrace{H_{2}\left(s_1^{(k)},s_2^{(k)}\right)}_{k^{\text{th}}~\text{measurement}}=\hat{H}_{2}\left(s_1^{(k)},s_2^{(k)}\right),
\end{equation}
{and we construct the following linear optimization problem that is solvable by minimizing the $2$-norm (least-squares) similarly to the quadratic-bilinear case in \cite{KARACHALIOS20227}}. Collecting $k$ pairs of measurements $(s_{1}^{(k)},s_{2}^{(k)})$, we conclude to:
\begin{equation}\label{eq:H2LS}
\underbrace{\left[\begin{array}{c}
H_{2}(s_{1}^{(1)},s_{2}^{(1)})\\[1mm]
H_{2}(s_{1}^{(2)},s_{2}^{(2)})
\\[1mm]
\vdots\\[1mm]
H_{2}(s_{1}^{(k)},s_{2}^{(k)})
\end{array}\right]}_{(k\times 1)}=\underbrace{\left[\begin{array}{c}
\hat{\bO}_2^{(1)}\otimes\hat{\bR}_2^{T(1)}\\[1mm]
\hat{\bO}_2^{(2)}\otimes\hat{\bR}_2^{T(2)}
\\[1mm]
\vdots\\[1mm]
\hat{\bO}_2^{(k)}\otimes\hat{\bR}_2^{T(k)}
\end{array}\right]}_{\bM:~ (k\times r^3)}\underbrace{\texttt{vec}(\hat{\bQ})}_{(r^3\times 1)}.
\end{equation}
The quadratic operator inherits symmetries, e.g., the terms $x_i x_j$ and $x_j x_i$ appear twice in the product $\bx \otimes \bx$. These symmetries are known by construction \cref{eq:Qsym} and can be appropriately handled. Despite considering these symmetries, a non-unique quadratic representation of the original system remains. The information from the 2nd $H_2$ is insufficient to infer the entirety of the quadratic operator entries. Algebraically, this can be explained by the rank deficiency of the least square matrix $\bM\in\IR^{k\times r^3}$. Real symmetry can be enforced in \cref{eq:H2LS} by including the conjugate counterparts. The above problem motivates the usage of higher harmonics (GFRFs) where the remaining parameters of the above under-determined problem can be estimated. In particular, evaluating the quadratic operator $\hat{\bQ}$ can be further parameterized with the non-empty null space we computed from the \cref{eq:H2LS} least-squares problem. 

\begin{remark}[Null space dimension] The dimension of the null space in \cref{eq:H2LS} for enough measurements $k\geq r$ and after enforcing the quadratic symmetries of the minimal quadratic system $(\bA_r,\bQ_r,\bB_r,\bC_r)$ scales with the state dimension $r$ as $\texttt{dim}(\texttt{ker}(\bM))=(r^3-2r^2+r)/2$. These are the remaining free parameters in the quadratic operator $\hat{\bQ}$ that cannot be estimated from \cref{eq:H2LS}. Therefore, the main idea is potentially obtaining the missing parameters from the higher-order GFRFs. The asymptotic behavior of the $\texttt{ker}(\bM)$ indicates that after enforcing the quadratic symmetries, only half of the entries in $\hat{\bQ}$ can be detected. As detailed next, this endeavor will lead to solving nonlinear quadratic vector equations. A unique solution for inferring the quadratic operator in this framework is not guaranteed in the general case. Indeed, for an arbitrary nonlinear system, we have no reason to expect that the 3rd GFRF kernel will be sufficient to infer the quadratic term. However, the aforementioned limitation to weakly nonlinear systems set by the VS representation significantly narrows the breadth of systems where inference through the 3rd GFRF is investigated. In addition, the problem formulation contains the optimal solution that could lead to identification or optimal model inference in the reduced dimension. However, the exact identification problem regarding how many GFRFs are needed remains open. Thus, we proceed with our method and seek to construct at least (sub)optimal models capable of explaining higher GFRFs with the potential of interpolating the whole Volterra series.
\end{remark}

The inferred quadratic operator has $r^3$ unknowns (less due to symmetries $(r^2(r+1)/2)$). If the rank of the matrix $\bM$ is $\texttt{rank}(\bM)=p<r^3$, the parametric solution of $\hat{\bQ}$ that we obtain from $H_{2}$ measurements with the dimension of the null space $m=r^3-p$ can be written as:
\begin{equation}\label{eq:Qsk}
\hat{\bQ}=\underbrace{\hat{\bQ}_{0}}_{\text{rank solution}}+\underbrace{\sum_{i=1}^{m}\lambda_{i}\hat{\bQ}_{i}}_{\text{parameterization}}
\end{equation}
The above splitting \cref{eq:Qsk} can be considered the same when the operators $\hat{\bQ}_0,~\hat{\bQ}_i,~i=1,\ldots,m$ are represented as vectors after vectorization due to the linear property of $\texttt{vec}(\cdot)$\footnote{The vectorization is row-wise, $vec(\bQ)=\left[\begin{array}{ccc}
\bQ(1,1:r^2) & \cdots & \bQ(r,1:r^2)
\end{array}\right]^T\in\IR^{r^3\times 1}$.}.
\subsubsection{Inferring the quadratic operator from the 3rd GFRF}
From the parameters $\lambda_i$ in \cref{eq:Qsk}, we also search those that explain the interpolation of the 3rd kernel. Therefore, we can write:
\begin{equation}\label{eq:hatH3}
    \hat{H}_{3}(s_1,s_2,s_3)=\hat{\bO}_{3}(s_1,s_2,s_3)\hat{\bQ}\hat{\bR}_{3}(s_1,s_2,s_3,\hat{\bQ}),
\end{equation}
and substituting \cref{eq:Qsk} in \cref{eq:hatH3}, due to the linear property of the operator $\bR_3$ as explained in \cref{sec:2.5}, we can derive
\begin{equation}
    \begin{aligned}
    &\hat{H}_{3}(s_1,s_2,s_3,\hat{\bQ},\hat{\bQ})=\hat{\bO}_{3}(s_1,s_2,s_3)\left(\hat{\bQ}_{0}+\sum_{i=1}^{m}\lambda_{i}\hat{\bQ}_{i}\right)\hat{\bR}_{3}\left(s_1,s_2,s_3,\hat{\bQ}_0+\sum_{i=1}^{m}\lambda_{i}\hat{\bQ}_{i}\right)\\
    &=\hat{\bO}_{3}(s_1,s_2,s_3)\hat{\bQ}_{0}\hat{\bR}_{3}\left(s_1,s_2,s_3,\hat{\bQ}_0\right)+\hat{\bO}_{3}(s_1,s_2,s_3)\hat{\bQ}_{0}\hat{\bR}_{3}\left(s_1,s_2,s_3,\sum_{i=1}^{m}\lambda_{i}\hat{\bQ}_{i}\right)+\\
    &\hat{\bO}_{3}(s_1,s_2,s_3)\left(\sum_{i=1}^{m}\lambda_{i}\hat{\bQ}_{i}\right)\hat{\bR}_{3}\left(s_1,s_2,s_3,\hat{\bQ}_0\right)\\
    &+\hat{\bO}_{3}(s_1,s_2,s_3)\left(\sum_{i=1}^{m}\lambda_{i}\hat{\bQ}_{i}\right)\hat{\bR}_{3}\left(s_1,s_2,s_3,\sum_{i=1}^{m}\lambda_{i}\hat{\bQ}_{i}\right)
    =\hat{H}_{3}^{(00)}(s_1,s_2,s_3,\hat{\bQ}_0,\hat{\bQ}_0)+\\
    &\sum_{i=1}^m\lambda_i\left(\hat{H}_{3}^{(i0)}(s_1,s_2,s_3,\hat{\bQ}_{i},\hat{\bQ}_0)+\hat{H}_{3}^{(0i)}(s_1,s_2,s_3,\hat{\bQ}_0,\hat{\bQ}_i)\right)+\\
    &+\sum_{i=1}^{m}\sum_{j=1}^{m}\lambda_i\lambda_j\hat{H}_{3}^{(ij)}(s_1,s_2,s_3,\hat{\bQ}_i,\hat{\bQ}_j),
    \end{aligned}
\end{equation}
where the superscript notation is similar to \cref{eq:ioTFscript}. The above problem can be written as a quadratic problem for root finding. We introduce the following notation: 
\begin{equation*}
    \begin{aligned}
        (\mathcal{A})_{i,j}&=\hat{H}_{3}^{(ij)}(s_{1},s_{2},s_{3},\hat{\bQ}_i,\hat{\bQ}_j),\\
        (\mathcal{B})_i&=\hat{H}_{3}^{(is)}(s_{1},s_{2},s_{3},\hat{\bQ}_i,\hat{\bQ}_0)+\hat{H}_{3}^{(si)}(s_{1},s_{2},s_{3},\hat{\bQ}_0,\hat{\bQ}_i),\\
        \mathcal{C}&=\hat{H}_{3}^{(ss)}(s_{1},s_{2},s_{3},\hat{\bQ}_0,\hat{\bQ}_0)-\hat{H}_{3}(s_{1},s_{2},s_{3},\hat{\bQ},\hat{\bQ}).
    \end{aligned}
\end{equation*}
 
We reformulate the problem by denoting $\blambda=\left[\begin{array}{cccc}
\lambda_{1} & \lambda_{2} & \cdots & \lambda_{m}\end{array}\right]^T$. 

The dimensions for a single measurement triplet $(s_1,s_2,s_3)$ are as follows $\mathcal{A}\in\IC^{n\times n},~\mathcal{B}\in\IC^{1\times n},~\mathcal{C}\in\IC$.
\begin{equation}\label{eq:quadraticvector}
\blambda^T\mathcal{A}\blambda+\mathcal{B}\blambda+\mathcal{C}=0.
\end{equation}

We can rewrite equation \eqref{eq:quadraticvector} in a more convenient format after vectorizing $\mathcal{A}$ that will allow arbitrary measurement augmentation as follows:
\begin{equation}
\texttt{vec}(\mathcal{A})(\blambda\otimes\blambda)+\mathcal{B}\blambda+\mathcal{C}=0.
\end{equation}

To enforce interpolation from the 3rd kernel, we equate 
\begin{equation}
\underbrace{H_3\left(s_1^{(k)},s_2^{(k)},s_3^{(k)}\right)}_{\text{k:~measurements}}=\hat{H}_3\left(s_1^{(k)},s_2^{(k)},s_3^{(k)},\hat{\bQ},\hat{\bQ}\right).   
\end{equation}
Further, by adding $k$ measurements along with the conjugate counterparts, we can transform the data matrices with real entries, and we result to 
\begin{equation}\label{eq:QuadVec}
\underbrace{\left[\begin{array}{c}
\texttt{vec}(\mathcal{A}_{1})\\
\texttt{vec}(\mathcal{A}_{2})\\
\vdots\\
\texttt{vec}(\mathcal{A}_{k})
\end{array}\right]}_{\bW\in\IR^{k\times m^2}}(\blambda\otimes\blambda)+\underbrace{\left[\begin{array}{c}
\mathcal{B}_{1}\\
\mathcal{B}_{2}\\
\vdots\\
\mathcal{B}_{k}
\end{array}\right]}_{\bZ\in\IR^{k\times m}}\blambda+\underbrace{\left[\begin{array}{c}
\mathcal{C}_{1}\\
\mathcal{C}_{2}\\
\vdots\\
\mathcal{C}_{k}
\end{array}\right]}_{\bS\in\IR^{k\times 1}}=\textbf{0}.
\end{equation}
The above equation can be written as $\bF(\blambda)=\textbf{0}$ by denoting $\bF(\cdot):\IR^{m}\rightarrow\IR^{k}$ as
\begin{equation}\label{eq:qvo}
\bF(\blambda)=\bW(\blambda\otimes\blambda)+\bZ\blambda+\bS,~\blambda\in\IR^{m}.
\end{equation}
The derivative (Jacobian) w.r.t the real vector $\blambda$ is:
\begin{equation}\label{eq:jqvo}
\bJ(\blambda)=\nabla_{\blambda}\bF(\blambda)=\bW(\blambda\otimes\bI+\bI\otimes\blambda)+\bZ.
\end{equation}
We seek the solution of \cref{eq:QuadVec}; thus, by introducing the Newton iterative procedure (fixed point iterations), starting with an initial seed $\blambda_0$, we can result in $\mathcal{F}(\blambda_{n+1})\rightarrow 0$ as $n\rightarrow\infty$. The iterations are described next:
\begin{equation}\label{eq:Newton_qve}
\blambda_{n+1}=\blambda_{n}-\bJ^{-1}(\blambda_{n})\bF(\blambda_{n}).
\end{equation}
Finally, upon Newton's method convergence, we obtain the vector $\blambda^*,~(\bF(\blambda^*)\approx\mathbf{0})$ from \cref{alg:qve}, which will lead to a better estimation of $\bQ$ that incorporates the measurements from the 3rd kernel. We notice in many situations that the error between the reduced and original systems improves significantly when the residual $\gamma$ of Newton's method remains small. Moreover, in many cases, identifying the original operator $\bQ$ is possible, as we will illustrate in our low-order example (i.e., the Lorenz attractor model).
\begin{algorithm}[!ht]
\caption{Solution $\blambda$ of the quadratic vector equation $\bF(\blambda)=\mathbf{0}$.}
\label{alg:qve}
\begin{algorithmic}
\STATE{Define: $\bW\in\IR^{k\times m^2},~\bZ\in\IR^{k\times m},~\bS\in\IR^{k\times 1}$ and the hyperparameters $\eta,~\gamma_0$.}
\STATE{Choose an initial random seed: $\blambda\in\IR^{m}$.}
\WHILE{$\gamma>\gamma_0$}
\STATE{Compute $\bF(\blambda)$ from \cref{eq:qvo} and $\bJ(\blambda)$ from \cref{eq:jqvo}}.
\STATE{Update $\blambda\leftarrow(\blambda-\bJ^{\#}\bF(\blambda))$,~$\#$ is: $"-1"$ or the Moore-Penrose pseudo-inverse (threshold $\eta$)}.
\STATE{Compute the residue $\lVert\bF(\blambda)\rVert=\gamma$}.
\ENDWHILE
\RETURN $\blambda$
\end{algorithmic}
\end{algorithm}

\begin{remark}[Reduction scheme for the quadratic vector equation]\label{rem:ProjectQVE}
    To specify $m$ unknowns from $N$ quadratic equations, instead of solving directly with \cref{alg:qve}, these equations can be projected in the appropriate dimension (dimension of the null space) $m$ with the SVD. To balance the interpolation goal approximately as the original system's dynamics evolve into a nonlinear manifold of a high dimension instead of only a reduced $r$. Thus, enforcing the exact interpolation for the third kernel might be impossible and result in nonsmooth fitting. Therefore, we consider the following augmented matrix $\bK=[\bW~~\bZ~~\bS]$, where SVD decomposes as $\bK=\bU\bSigma\bV^T$. We project with the left $m$ singular-vectors as $\hat{\bK}=\bU(:,1:m)^T\bK=\left[\begin{array}{ccc}
    \hat{\bW} & \hat{\bZ} & \hat{\bS}\end{array}\right]$ and we solve the reduced quadratic vector equation with \cref{alg:qve}.
\end{remark}

\begin{remark}[Identification of the minimal quadratic system]\label{rem:identification}
The proposed method works with measurements from the first three symmetric GFRFs that can be inferred from the time domain. All the parameters to achieve identification have been exploited. Regarding the uniqueness of the solution vector $\blambda$ in \cref{alg:qve}, it could be the case where the important degrees of freedom can be fixed from the first three GFRFs and achieve quadratic identification, but in the general case more information should be included from even higher kernels $(H_4,~H_5,\ldots)$. Consequently, \cref{alg:qve} should be upgraded in solving the generalized vector equation coupled in an arbitrary finite polynomial degree not specified in this study. In such cases, the complexity increases in both data acquisition and algorithmic complexity, where specifying the possible finite kernel degree to which quadratic identification can be achieved is left as an open problem. Although such a theoretical statement seems complex and has a high computational burden, the proposed method works toward the identification goal using information from the first most significant kernels. In the VS framework, higher kernels/harmonics carry less important information under the assumption of weakly nonlinear dynamics. Thus, an efficient approximation can be provided, as detailed in our following examples.   
\end{remark}
\subsubsection{The algorithm for quadratic modeling from i/o time-domain data}\label{sec:algo2}
Here, we present a concise algorithm summarizing the procedure for constructing quadratic state space models from harmonic data (samples of the symmetric kernels $H_1,~H_2,~H_3$). Measuring (symmetric) Volterra kernels is, by no means, a new topic. However, although previously addressed in \cite{Boyd1983MeasuringVK,morKarGA21a,VoltMacrXpar}, it remains a non-trivial task. The main difficulty has to do with the separation of commensurate frequencies. In other words, each of the propagating harmonics consists of a series of kernels and, therefore, evaluating the symmetric GFRFs requires kernel separation with an amplitude shifting \cite{morKarGA21a,VoltMacrXpar}. Towards this aim, X-parameters in \cite{Xpar2006}, and the references within, represent a direct generalization of the classical S-parameters (for linear dynamics) to the nonlinear case. With this agile machinery, estimations of the higher Volterra kernel can be made in a true engineering setup as in \cite{VoltMacrXpar}, and the quadratic state-space surrogate model can be inferred from the proposed method. The following algorithm can use such information (from the X-parameters) to construct interpretable quadratic models.
\begin{algorithm}[!ht]
\caption{Quadratic modeling from time-domain data}
\label{alg:Qmodel}
\begin{algorithmic}
\STATE{Input: \# Measurements of the symmetric GFRFs $H_1(s_1),~H_2(s_1,s_2),~H_3(s_1,s_2,s_3)$.}
\STATE{Define a truncation order $r$ with SVD from the Loewner matrix $\IL$ (minimal linear).}
\STATE{Realize the minimal linear subsystem $(\hat{\bA},~\hat{\bB},~\hat{\bC})$ of order $r$}.
\STATE{Estimate the $\hat{\bQ}_0\in\IR^{r\times r^2}$ from \cref{eq:H2LS}} by minimizing the $2-$norm error (least-squares).
\STATE{Update the $\hat{\bQ}\in\IR^{r\times r^2}$ from \cref{eq:Qsk} after solving \cref{eq:Newton_qve} with \cref{alg:qve}}.
\RETURN the state-space quadratic model $(\hat{\bA},~\hat{\bQ},~\hat{\bB},~\hat{\bC})$.
\end{algorithmic}
\end{algorithm}
\subsection{Quadratic state-space systems with multiple equilibrium points}\label{sec:multiEqulibria}
Quadratic systems can bifurcate to different equilibrium points. Thus, when measuring, multi-operational points can be revealed. To illustrate this phenomenon mathematically, we write the quadratic system \cref{eq:qsys} after shifting it with the non-zero equilibrium state $\bx_e$. We denote the new state variable $\tilde{\bx}(t)=\bx(t)-\bx_e$, where the tilde notation $\tilde{(\cdot)}$ in that case should enhance the memory that the corresponding matrix/state operates close to the equilibrium. The system remains
\begin{equation}\label{eq:bifurcatedmodel}
     \begin{aligned}
    \dot{\bx}(t)&=\bA\bx(t)+\bQ(\bx(t)\otimes\bx(t))+\bB u(t)\Rightarrow\\
    \dot{\tilde{\bx}}(t)&=\bA(\tilde{\bx}(t)+\bx_e)+\bQ\left((\tilde{\bx}(t)+\bx_e)\otimes(\tilde{\bx}(t)+\bx_e)\right)+\bB u(t)\Rightarrow\\
    \dot{\tilde{\bx}}(t)&=\bA\tilde{\bx}(t)+2\bQ(\bx_e\otimes\tilde{\bx}(t))+\bQ(\tilde{\bx}(t)\otimes\tilde{\bx}(t))+
    \bA\bx_e+\bQ(\bx_e\otimes\bx_e)+\bB u(t)\Rightarrow\\
    \dot{\tilde{\bx}}(t)&=\underbrace{\left(\bA+2\bQ(\bx_e\otimes\bI)\right)}_{\tilde{\bA}}\tilde{\bx}(t)+\bQ(\tilde{\bx}(t)\otimes\tilde{\bx}(t))+\underbrace{\bA\bx_e+\bQ(\bx_e\otimes\bx_e)}_{{\bL}}+\bB u(t).
 \end{aligned}
\end{equation}

Note that at the equilibrium equation holds $\bL:=\bA\bx_e+\bQ(\bx_e\otimes\bx_e)=\mathbf{0}$ as in the absence of the controller $u(t)$, and with zero initial conditions, e.g., $\bx_0=\mathbf{0}$, dynamics remain zero. We do not address situations with a limit cycle, e.g., systems with purely imaginary eigenvalues, which describe self-sustained dynamics. As a result, the quadratic system that we measure after reaching the equilibrium state $\bx_e$ is the following:
\begin{equation}\label{eq:QsysEquilibrium}
    \left\{\begin{aligned}
        \dot{\tilde{\bx}}(t)&=\tilde{\bA}\tilde{\bx}(t)+\bQ(\tilde{\bx}(t)\otimes\tilde{\bx}(t))+\bB u(t),\\
        y(t)&=\bC\tilde{\bx}(t)+\bC\bx_e,~\tilde{\bx}_0=-\bx_e,~t\geq 0.
    \end{aligned}\right.
\end{equation}
\begin{remark}[System invariant information under bifurcations]
    The system in \cref{eq:QsysEquilibrium} suggests that around the new equilibrium state point $\bx_e$, the matrices $(\bQ,~\bB,~\bC)$ stay the same compared to the global system's matrices up to a coordinate transformation, and only the linear matrix changes to $\tilde{\bA}=\bA+2\bQ(\bx_e\otimes\bI)$, along with a linear trend known as DC\footnote{DC: direct current, which describes the non-periodic term (zero frequency) in the power spectrum.}. The generalized Markov parameters of the system that contain only the matrices $(\bQ,~\bB,~\bC)$ are the same around any arbitrary equilibrium $\bx_e$ and any coordinate system to which the original system bifurcates.  
\end{remark}
\textbf{Two equilibrium points case:} To emphasize the operation of the dynamics around the equilibrium points, we keep the notation to all matrices and states with a tilde. Let's assume that the original global quadratic model can bifurcate to the two different equilibrium points $\tilde{\bx}_e^{(q)},~q=1,2$. The dynamical behavior can be explained locally and in every equilibrium point by the following systems
\begin{equation}\label{eq:syss}
    \left\{\begin{aligned}
\dot{\tilde{\bx}}_q(t)&=\tilde{\bA}_q\tilde{\bx}_q(t)+\tilde{\bQ}_q(\tilde{\bx}_q(t)\otimes\tilde{\bx}_q(t))+\tilde{\bB}_q u(t),\\
        y_q(t)&=\tilde{\bC}_q\tilde{\bx}_q(t),~\tilde{\bx}_q(0)=-\tilde{\bx}_e^{(q)},
        \end{aligned}\right.,~q=1,2.
\end{equation}

\textbf{Property}: For the two systems in \cref{eq:syss} holds
    \begin{equation}\label{eq:syssq}
\underbrace{\tilde{\bA}_q}_{\text{local}}=\underbrace{\bA_q}_{\text{global}}+2\tilde{\bQ}_q(\tilde{\bx}_e^{(q)}\otimes\bI),~q=1,2.
    \end{equation}

\begin{figure}[!ht]
    \centering
    \includegraphics[scale=0.26]{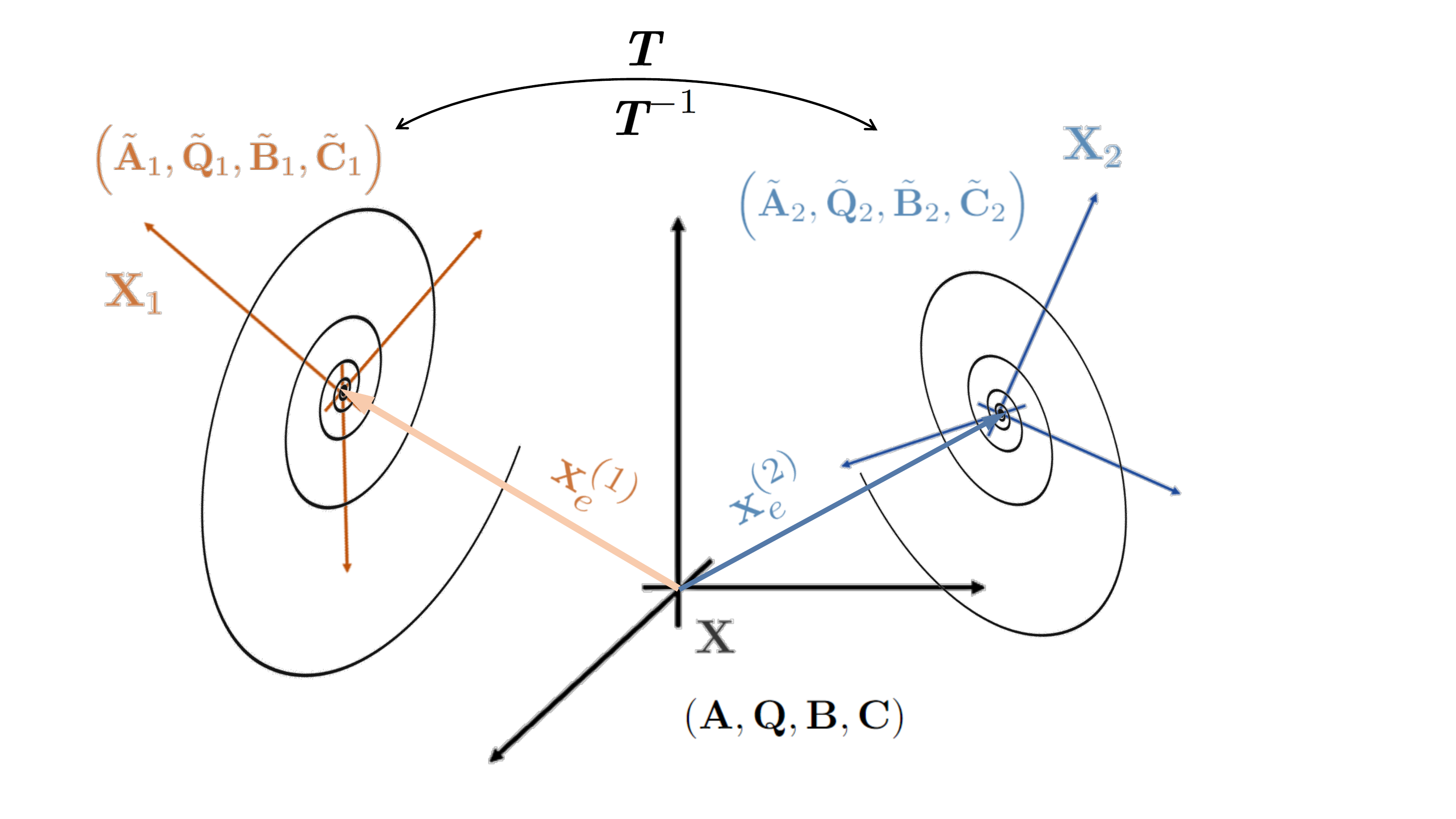}
    \vspace{-2mm}
    \caption{\cref{lem:2.6} through a schematic.}
    \label{fig:align}
\end{figure}

\begin{remark}\label{rem:MarkovQBC}
The Markov parameters involving the quadratic and input-output operators are the same. $\tilde{\bC}_1\tilde{\bB}_1=\tilde{\bC}_2\tilde{\bB}_2,~\text{and}~\tilde{\bC}_1\tilde{\bQ}_1(\tilde{\bB}_1\otimes\tilde{\bB}_1)=\tilde{\bC}_2\tilde{\bQ}_2(\tilde{\bB}_2\otimes\tilde{\bB}_2)$.
\end{remark}
According to \cref{rem:MarkovQBC}, similar information from the original system is encoded in both systems \cref{eq:syss}. Therefore, a similarity transformation $\bT$ exists, which aligns the two systems w.r.t the original coordinates.

\begin{lemma}[Coordinate alignment with unknown global linear matrices]\label{lem:2.6}
There exists a transformation matrix $\bT$ such that the two triplets given by $(\tilde{\bQ}_1,\tilde{\bB}_1,\tilde{\bC}_1)$ and by $(\tilde{\bQ}_2,\tilde{\bB}_2,\tilde{\bC}_2)$ (resulting after a global model bifurcates to different equilibrium points) can be aligned simultaneously with the original operators but to different coordinates, geometrically in \cref{fig:align}, and algebraically as
\begin{equation}\label{eq:Tsys}
\begin{aligned}
\tilde{\bQ}_1&=\bT\tilde{\bQ}_2(\bT^{-1}\otimes\bT^{-1}),~\tilde{\bB}_1=\bT^{-1}\tilde{\bB}_2\Leftrightarrow\bT\tilde{\bB}_1=\tilde{\bB}_2\\
    \tilde{\bC}_1&=\bT\tilde{\bC}_2,~\bA_1=\bT\bA_2\bT^{-1},~(\bA_1,~\bA_2,~\text{global system's matrices}).
    \end{aligned}
\end{equation}
\end{lemma}
One way to compute the transformation matrix $\bT$ is by solving the first three equations in system \cref{eq:Tsys}. The above problem involves a quadratic matrix equation that can be iteratively solved employing Newton iterations. Moreover, adding linear constraints prevents the convergence of the Newton method to the zero solution. To seek a formal solution, we analytically derive the iterative Newton scheme over the Fr\'echet derivative in what follows (\cref{sec:ConQME}). 

\subsubsection{Computing the equilibrium points}\label{sec:equilibria}
When measurements around equilibria have been collected, the proposed method in \cref{alg:Qmodel} can infer a quadratic system. In particular, the inferred quadratic system is localized around the reachable equilibrium point $q=1,2$; thus, the available matrices are $(\tilde{\bA}_q,\tilde{\bQ}_q,\tilde{\bB}_q,\tilde{\bC}_q),~q=1,2$. The systems $\Sigma_q,~q=1,2$ are not necessarily equivalent to the original (global) one because a bifurcation might have occurred. When the matrix $\bA_q$ is without the tilde notation, it represents the global one, and the subscript "q" indicates the different coordinate system. We can discover the global $\bA_q$ matrix, from \cref{eq:syssq} as $\bA_q=\tilde{\bA}_q-\tilde{\bQ}(\tilde{\bx}_e^{(q)}\otimes\bI)$ after first computing the equilibrium point $\tilde{\bx}_e^{(q)}$ for $q=1,2$. In general, it holds
\begin{equation}\label{eq:equilibriumhattedsystem}
     \bL_{q}(\tilde{\bx}_e^{(q)}):=\mathbf{0}\Leftrightarrow\bA_q\tilde{\bx}_e^{(q)}+\tilde{\bQ}(\tilde{\bx}_e^{(q)}\otimes\tilde{\bx}_e^{(q)})\Leftrightarrow\tilde{\bA}_q\tilde{\bx}_e^{(q)}-\tilde{\bQ}(\tilde{\bx}_e^{(q)}\otimes\tilde{\bx}_e^{(q)})=\mathbf{0}.     
\end{equation}
\cref{eq:equilibriumhattedsystem} suggests that the equilibrium state vector of dimension $n$ can be computed as the unknowns match the equations and the matrices $\tilde{\bA}_q,~\tilde{\bQ}$ are available. \Cref{alg:qve} can then be used directly, but the Newton scheme will converge to the zero solution. Thus, we must constrain our solution to satisfy the observable equilibrium from the DC term. In particular, and for every equilibrium point, additional information comes from the direct current (DC) terms that can be measured from the power spectrum $\alpha_1,~\alpha_2$. For instance, we enforce
\begin{equation}\label{eq:dcterms}
        \tilde{\bC}_q\tilde{\bx}_e^{(q)}=\alpha_q.
\end{equation}
Consequently, we must solve the coupled equations \cref{eq:equilibriumhattedsystem} and \cref{eq:dcterms} for the unknown state equilibrium vector $\tilde{x}_e^{(q)}$. In particular, \cref{eq:dcterms} can give a parametric solution as
\begin{equation}\label{eq:parametricequil}
    \tilde{\bC}_q\tilde{\bx}_e^{(q)}=\alpha_q\Leftrightarrow \tilde{\bx}_e^{(q)}=\tilde{\bx}_{es}^{(q)}+\sum_{i=1}^{r-1}\lambda_i\tilde{\bx}_{ei}^{(q)},~q=1,2.
\end{equation}
By substituting \cref{eq:parametricequil} in \cref{eq:equilibriumhattedsystem}, we result for each system with $q=1,2$ to the following quadratic vector equation
\begin{equation}\label{eq:L0}
    \begin{aligned}      
&\tilde{\bA}_q\left(\tilde{\bx}_{es}^{(q)}+\sum_{i=1}^{r-1}\lambda_i\bx_{ei}^{(q)}\right)-\tilde{\bQ}_q\left(\tilde{\bx}_{es}^{(q)}+\sum_{i=1}^{r-1}\lambda_i\bx_{ei}^{(q)}\right)\otimes\left(\tilde{\bx}_{es}^{(q)}+\sum_{i=1}^{r-1}\lambda_i\tilde{\bx}_{ei}^{(q)}\right)=\mathbf{0}\\
        &\underbrace{\tilde{\bA}_q\tilde{\bx}_{es}^{(q)}-\tilde{\bQ}_{q}\left(\tilde{\bx}_{es}^{(q)}\otimes\tilde{\bx}_{es}^{(q)}\right)}_{\bS\in\IR^{r\times 1}}+\sum_{i=1}^{r-1}\lambda_i\underbrace{\left(\tilde{\bA}_q\bx_{ei}^{(q)}-2\hat{\bQ}_q\bx_{es}^{(q)}\otimes\bx_{ei}^{(q)}\right)}_{\bZ\in\IR^{r\times r}}-\\
        &-\sum_{i=1}^{r-1}\sum_{j=1}^{r-1}\lambda_i\lambda_j\underbrace{\tilde{\bQ}_{q}\left(\bx_{ei}^{(q)}\otimes\bx_{ej}^{(q)}\right)}_{\bW\in\IR^{r\times r^2}}=\mathbf{0}.
    \end{aligned}
\end{equation}
where solution can be obtained after enforcing $\bW(\blambda\otimes\blambda)+\bZ\blambda+\bS=\mathbf{0}$ with the developed \cref{alg:qve} and for each system $q=1,2$. Having the vector $\blambda^\ast=\left[\begin{array}{ccc}
    \lambda_1 & \cdots & \lambda_{r-1} \\
\end{array}\right]^T$ that solves \cref{eq:L0}, we compute the equilibrium states $\tilde{\bx}_e^{(q)}$ from $\cref{eq:parametricequil}$. The original linear operator of the quadratic system can be computed from $\bA_q=\tilde{\bA}_q-\tilde{\bQ}(\tilde{\bx}_e^{(q)}\otimes\bI)$ that lives to a coordinate system where differs from the original (global) system up to a similarity transformation \cref{app:align}.
\section{Numerical results}\label{sec:Results}
We test the proposed method in a more involved introductory example as a next step from the linear one tailored to quadratic structure; afterward, we proceed to different cases of identifying the Lorenz'63 system operating in different parametric regimes. Finally, the viscous Burgers' equation model illustrates the effectiveness of inferring a reduced surrogate model with the proposed data-driven method directly from i/o measurements. 

\subsection{An introductory quadratic toy example}
Consider the following SISO quadratic system of dimension $n=2$ with non-zero initial condition
\begin{equation}\label{eq:smallquadratic}
    \left\{\begin{aligned}
        \dot{\bx}(t)&=\left[\begin{array}{cc}
           1  & 0\\
           0 & -2
        \end{array}\right]\bx(t)+\left[\begin{array}{cccc}
            -1 & 0 & 0 & 0  \\
             0 & 0 & 0 & 0
        \end{array}\right]\bx(t)\otimes\bx(t)+\left[\begin{array}{c}
             1  \\
             1 
        \end{array}\right]u(t),\\
        y(t)&=\left[\begin{array}{cc}
            1 & 1 \end{array}\right]\bx(t),~y(0)=\bC\bx(0)=\left[\begin{array}{cc}
                1 & 1 \\
            \end{array}\right]\left[\begin{array}{c}
                 0.5  \\
                 0 
            \end{array}\right]=0.5,~t\geq 0.
    \end{aligned}\right.
\end{equation}

\begin{table}[h!]
    \centering
    \begin{tabular}{l|cc}
     Input amplitude  &  $\Vert P(0)-\bC\bx_e\Vert $ & $\Vert P(2\pi j10)/a-H_1(2\pi j10)\Vert$ \\
        $2a=20$ & $5.3500e-02$ & $2.7105e-05$\\
        $2a=2$ & $5.0673e-04$ & $2.5666e-07$\\
        $2a=0.002$ & $5.0647e-10$ & $2.1652e-11$\\
    \end{tabular}
    \caption{Input as $u(t)=2a\cos(2\pi 10t)$ to \cref{eq:smallquadratic} and measurement estimation errors of the 1st GFRF at $10$ (Hz) and of the DC term as the observable equilibrium.}
    \label{tab:kernelestimation}
\end{table}

To illustrate the data acquisition in our example with state-space representation as in \cref{eq:smallquadratic}, we simulate with single-tone inputs by sweeping\footnote{This technique of sweeping the amplitude of the multi-harmonic signal helps the Kernel separation.} the amplitude of the signal and compare with the theoretical measurement in \cref{tab:kernelestimation}. As we explained in \cref{sec:dataaquisition}, the propagating harmonics are $(k\omega-p\omega),~k,p\in\IZ$. The higher order symmetric GFRFs scale with the input amplitude, $a^nH_n(\omega,\ldots,\omega)$. Thus, with a small enough amplitude $a$, the higher-order kernels will decay faster to zero than the lower-order kernels, directly allowing good estimations.

In \cref{tab:kernelestimation}, we provide measurements that can be accurately obtained when real inputs have been considered. To avoid the above-known complications, we will consider complex inputs where the indexing of harmonics and kernels is injective. Under this assumption, the rest of the data acquisition will be generated by sampling the derived GFRFs as efficient methods and machinery exist to solve the real case.  

In \cref{eq:smallquadratic}, the linear operator is unstable (due to the eigenvalue $1$), but the quadratic system can operate to stable regimes. The measurements that we achieved in \cref{tab:kernelestimation} are associated with the bifurcated model around its quadratically stable equilibrium $\bx_e=\left[\begin{array}{cc}
    1 & 0 \\
\end{array}\right]^T$ as analyzed in \cref{sec:multiEqulibria}. Thus, the system that we measure subsequently is
\begin{equation}\label{eq:smallquadraticequilibria}
     \left\{\begin{aligned}
        \dot{\bx}(t)&=\underbrace{\left[\begin{array}{cc}
           -1  & 0\\
           0 & -2
        \end{array}\right]}_{\bA+2\bQ(\bx_e\otimes\bI)}\bx(t)+\left[\begin{array}{cccc}
            -1 & 0 & 0 & 0  \\
             0 & 0 & 0 & 0
        \end{array}\right]\bx(t)\otimes\bx(t)+\left[\begin{array}{c}
             1  \\
             1 
        \end{array}\right]u(t),\\
        y(t)&=\left[\begin{array}{cc}
            1 & 1 \end{array}\right]\bx(t)+\underbrace{1}_{\bC\bx_e},~\bx(0)=\bx_0=\left[\begin{array}{c}
                 0.5  \\
                 0 
            \end{array}\right]-x_e,~t\geq 0.
    \end{aligned}\right.
\end{equation}
The linear subsystem \cref{eq:smallquadraticequilibria} is exactly the linear system we identified in our introductory example in \cref{sec:linearintro}. Thus, we can proceed with the same linear identification result. To identify the quadratic system, we must identify the remaining quadratic operator $\hat{\bQ}\in\IR^{2\times 4}$ along with the initial conditions. Following the methodology in \cref{sec:estimationQfromH2}, we compute
\begin{equation}
    \hat{\bQ}_0=\left[\begin{array}{cccc} 0.11335 & -0.0076666 & -0.0076666 & 0.00049333\\ 2.3282 & -0.1619 & -0.1619 & 0.010564 \end{array}\right],
\end{equation}
and the null space dimension of the linear system is $1$. Thus, the parametric solution can be written
\begin{equation}
    \hat{\bQ}=\hat{\bQ}_0+\lambda\hat{\bQ}_1,~\lambda\in\IR.
\end{equation}
In such a case, the \cref{eq:quadraticvector} simplifies to a scalar quadratic as
\begin{equation}
    \hat{H}_3^{(11)}(s_1,s_2,s_3)\lambda^2+(\hat{H}_3^{(01)}(s_1,s_2,s_3)+\hat{H}_3^{(10)}(s_1,s_2,s_3))\lambda+\hat{H}_3^{(00)}(s_1,s_2,s_3)-\underbrace{H_3(s_1,s_2,s_3)}_{\text{measurement}}=0.
\end{equation}
With two measurements of the $H_3$ and by including the conjugate pairs to enforce real entries, we get $\lambda^\ast=-0.090951$. The corrected quadratic operator results to
\begin{equation}
\hat{\bQ}=\hat{\bQ}_0+\lambda^\ast\hat{\bQ}_1=\left[\begin{array}{cccc} 0.11375 & -0.0047759 & -0.0047759 & 0.00020051\\ 2.3371 & -0.098122 & -0.098122 & 0.0041196 \end{array}\right].
\end{equation}
What is the initial condition state vector $\hat{\bx}_0$ to this coordinate system? To address this issue, first, we need to notice that $\hat{\bC}\hat{\bx}_e=1$. Moreover, the equilibrium to the hatted system satisfies $\tilde{\bA}\hat{x}_e+\hat{\bQ}(\hat{\bx}_e\otimes\hat{\bx}_e)=\mathbf{0}$, but we don't know the operator $\tilde{\bA}$. As explained in \cref{sec:equilibria}, we have
\begin{equation}\label{eq:equilsmallquad}
\begin{aligned}
    &\tilde{\bA}\hat{x}_e+\hat{\bQ}(\hat{\bx}_e\otimes\hat{\bx}_e)=\mathbf{0}\Leftrightarrow\hat{\bA}\hat{\bx}_e-\hat{\bQ}(\hat{x}_e\otimes\hat{\bx}_e)=\mathbf{0}.
    \end{aligned}
\end{equation}
Computing the equilibrium $\hat{\bx}_e$ from \cref{eq:equilsmallquad} with \cref{alg:qve} seems feasible, but the solution might converge to the origin that would not explain the nonzero DC term. Thus, we must enforce information from $\hat{\bC}\hat{\bx}_e=1$. Our unknown state equilibrium $\hat{\bx}_e$ has $r=2$ degrees of freedom. Thus, we write
\begin{equation}
\begin{aligned}
    &\hat{\bC}\hat{\bx}_e=1\Leftrightarrow\left[\begin{array}{cc}
        \hat{c}_1 & \hat{c}_2 \\
    \end{array}\right]\left[\begin{array}{c}
         \hat{x}_e^{(1)}  \\
         \hat{x}_e^{(2)}
\end{array}\right]=1\Leftrightarrow\hat{c}_1\hat{x}_e^{(1)}+\hat{c}_2\hat{x}_e^{(2)}=1\overset{\hat{c}_1\neq 0}{\Leftrightarrow}\hat{x}_e^{(1)}=\frac{1}{\hat{c}_1}-\frac{\hat{c}_2}{\hat{c}_1}\hat{x}_e^{(2)}.
\end{aligned}
\end{equation}
 To exploit the affine structure of the system near the equilibrium, we set $p_1=1/\hat{c}_1,~p_2=-\hat{c}_2/\hat{c}_1$ and $\hat{x}_e^{(2)}=\lambda_1\in\IR$ and it results to
\begin{equation}\label{eq:affineequil}
    \hat{\bx}_e=\left[\begin{array}{c}
         \hat{\bx}_e^{(1)}  \\
         \hat{\bx}_e^{(2)} 
    \end{array}\right]=\left[\begin{array}{c}
         p_1+p_2\hat{\bx}_e^{(2)}  \\
         \hat{\bx}_e^{(2)} 
    \end{array}\right]=\underbrace{\left[\begin{array}{c}
         p_1  \\
          0
    \end{array}\right]}_{{\bp}_1}+\underbrace{\left[\begin{array}{c}
         p_2  \\
          1
    \end{array}\right]}_{{\bp}_2}\lambda_1=\bp_1+\bp_2\lambda_1.
\end{equation}
Substituting \cref{eq:affineequil} to \cref{eq:equilsmallquad}, we result to
\begin{equation}\label{eq:quadrequismall}
\begin{aligned}
    &\hat{\bA}(\bp_1+\bp_2\lambda_1)-\hat{\bQ}\left((\bp_1+\bp_2\lambda_1)\otimes(\bp_1+\bp_2\lambda_1)\right)=0\Leftrightarrow\\
    &\hat{\bQ}(\bp_2\otimes\bp_2)\lambda_1^2+\left(2\hat{\bQ}(\bp_2\otimes\bp_1)-\hat{\bA}\bp_2\right)\lambda_1+\hat{\bQ}(\bp_1\otimes\bp_1)-\hat{\bA}\bp_1=\mathbf{0}.
\end{aligned}
\end{equation}
To solve \cref{eq:quadrequismall}, we can use the \cref{alg:qve}, which results in $\lambda_1=-9.5656e+03$.

Using \cref{alg:qve}, we compute the equilibrium state at the hatted system.
\begin{equation}
    \hat{\bx}_e=\left[\begin{array}{c} -465.58\\ -9565.6 \end{array}\right],~\hat{\bC}\hat{\bx}_e=1~~\text{(machine precision accuracy)}.
\end{equation}
The global linear operator can be computed from 
\begin{equation}
    \tilde{\bA}=\hat{\bA}-2\hat{\bQ}(\hat{\bx}_e\otimes\bI)=\left[\begin{array}{cc} 19.832 & -0.91663\\ 448.55 & -20.832 \end{array}\right],~\texttt{eig}(\tilde{\bA})=\{1,-2\}.
\end{equation}
After measuring the local dynamic behavior, we have identified the global quadratic system's dynamics. Finally, we must align the coordinates to use the given observable initial conditions $y_0=0.5$. As in the linear case, we must infer information from transient dynamics to discover the initial condition. First, we switch off the controller by setting $u(t)=0$, and we track the dynamic evolution of the autonomous quadratic system
\begin{equation}\label{eq:quadraticintrou0}
    \left\{\begin{aligned}
        \dot{\bx}(t)&=\bA\bx(t)+\bQ(\bx(t)\otimes\bx(t)),\\
        y(t)&=\bC\bx(t),~y_0=y(0)=\bC\bx_0=0.5~~\text{(the observable initial condition)}.
    \end{aligned}\right.
\end{equation}
Although we can discretize \cref{eq:quadraticintrou0} with higher order schemes (e.g., multi-step Runge Kutta), here to illustrate the theoretical result, we use a forward Euler scheme $\dot{\bx}(t)=(\bx(t+h)-\bx)/h,~h\neq 0$. The discrete identified system to the hatted coordinates is denoted as
\begin{equation}
    \left\{\begin{aligned}
        \hat{\bx}_{k+1}&=\tilde{\bA}_d\hat{\bx}_k+\hat{\bQ}_d(\hat{\bx}_k\otimes\hat{\bx}_k),\\
        \hat{y}_k&=\hat{\bC}_d\hat{\bx}_k,~\hat{y}_0=\hat{y}(0)=\hat{\bC}_d\hat{\bx}_0=0.5~~\text{(the observable initial condition)}.
    \end{aligned}\right.
\end{equation}
The unknown initial condition vector to the hatted coordinate system scales with $r=2$. Therefore, we need to specify two degrees of freedom, which makes two equations. We enforce
\begin{equation}
    \left\{\begin{aligned}
        \hat{y}_0&=y_0,\\
        \hat{y}_1&=y_1.
    \end{aligned}\right.\Rightarrow\left\{\begin{aligned}
        \hat{\bC}_d\hat{\bx}_0&=y_0,\\
        \hat{\bC}_d\tilde{\bA}_d\hat{\bx}_0+\hat{\bC}_d\hat{\bQ}_d(\hat{\bx}_0\otimes\hat{\bx}_0)&=y_1.
    \end{aligned}\right.
\end{equation}
We follow similar arguments as before in our previous methodology on computing the equilibrium state at the hatted coordinate system. The first equation gives a linear dependency among the entry elements of $\hat{\bx}_0$, and we can parameterize the affine solution as $\hat{\bx}_0=\bq_1+\kappa\bq_2$ with $\kappa\in\IR$. using the same methodology and the \cref{alg:qve}, we get $\kappa=-4782.8$. The initial condition state vector to the hatted coordinate system is
\begin{equation}
    \hat{\bx}_0=\bq_1+\kappa\bq_2=\left[\begin{array}{c} -232.79\\ -4782.8 \end{array}\right].
\end{equation}
To illustrate that the original quadratic system and the discovered one are the same, we construct a similarity transform $\Psi$ from the \cref{app:align}, and we conclude to
\begin{equation}
         \left\{\begin{aligned}
        \dot{\bx}_{\text{id}}(t)&=\left[\begin{array}{cc}
           1  & 0\\
           0 & -2
        \end{array}\right]\bx_{\text{id}}(t)+\left[\begin{array}{cccc}
            -1 & 0 & 0 & 0  \\
             0 & 0 & 0 & 0
        \end{array}\right]\bx_{\text{id}}(t)\otimes\bx_{\text{id}}(t)+\left[\begin{array}{c}
             1  \\
             1 
        \end{array}\right]u(t),\\
        y_{\text{id}}(t)&=\left[\begin{array}{cc}
            1 & 1 \end{array}\right]\bx_{\text{id}}(t),~\bx_{\text{id}}(0)=\left[\begin{array}{c}
                 0.5  \\
                 0 
            \end{array}\right],~t\geq 0.
    \end{aligned}\right.
\end{equation}

\begin{remark}[Initial condition discovery for the quadratic system]
    To identify the initial conditions in the quadratic system's arbitrary state dimension $n$, we need transient dynamics that scale with the discovered dimension $r$. Starting with $\hat{\bC}\hat{\bx}_0$, which is the observable initial condition, we can parameterize the unknown $\hat{\bx}_0$ with an affine structure as
    \begin{equation}
        \hat{\bx}_0=\bp_0+\sum_{i=1}^{r-1}\lambda_i\bp_i,~\lambda_i\in\IR,~(\bp_0,~\bp_i~\text{known}).
    \end{equation}
    The recursive formulation below exploits the Markov parameters for the transient dynamics in the absence of the controller and for the quadratic system.
    \begin{equation}
        \left\{\begin{aligned}
            \hat{\bx}_{k+1}&=\hat{\bA}\hat{\bx}_k+\hat{\bQ}(\hat{\bx}_k\otimes\hat{\bx}_k),~\hat{y}_0=\hat{\bC}\hat{\bx}_0,\\
            \hat{y}_{k+1}&=\hat{\bC}\hat{\bA}\hat{\bx}_k+\hat{\bC}\hat{\bQ}(\hat{\bx}_k\otimes\hat{\bx}_k),~k=0,\ldots,r-1.
        \end{aligned}\right.
    \end{equation}
    Thus, to specify the free parameters $\lambda_i,~i=1,\ldots,r-1$, by denoting $\blambda^{\otimes r}=\blambda\overbrace{\otimes\cdots\otimes}^{r}\blambda$, we need to solve the quadratic polynomial vector equation
 $\bW_r\blambda^{\otimes 2r}+\bW_{2r-1}\blambda^{\otimes(2r-1)}+\cdots+\bW_1\blambda+\bW_0=0$.
\end{remark}
The above result explains the difficulty in discovering initial conditions from i/o measurements for nonlinear systems. Another thought is linearization around an operational point where the initial conditions can affect it. In such a case, the only known linearization state point is the system's equilibrium. Thus, we expect a good agreement between the actual response of the nonlinear system and the linearized one close to equilibrium. Unfortunately, the initial conditions that affect the equilibrium through the transient phase become less important without inferring the initial condition with satisfying accuracy. In conclusion, the problem remains challenging, and a higher polynomial vector equation that scales with $2r$ should be solved. Thus, we proceed with zero initial conditions for the larger-scale examples.   
\subsection{The Lorenz system}
We consider the canonical Lorenz '63 model, from \cite{Lorenz}, adding a control-input $u(t)$ in the 1st and 3rd states. The following state-space form describes the quadratic control system
\begin{equation}
    \left\{\begin{aligned}
    \dot{x}(t)&=-\sigma x(t)+\sigma y(t)+u(t),\\
    \dot{y}(t)&=\rho x(t) -y(t)-x(t)z(t),\\
    \dot{z}(t)&=-\beta z(t)+x(t)y(t)+u(t),
    \end{aligned}\right.
\end{equation}
where zero initial condition are assumed e.g., $(x(0),~y(0),~z(0))=(0,~0,~0)$, and the operators are:
\begin{equation}\label{eq:Lorenz}
    \begin{aligned}
        \bA=\left[\begin{array}{ccc}
            -\sigma & \sigma & 0 \\
             \rho & -1 & 0\\
             0 & 0 & -\beta
        \end{array}\right],~\bB=\bC^T=\left[\begin{array}{c}
             1  \\
             0 \\
             1
        \end{array}\right],~\bQ=\left[\begin{array}{ccccccccc} 0 & 0 & 0 & 0 & 0 & 0 & 0 & 0 & 0\\ 0 & 0 & -\frac{1}{2} & 0 & 0 & 0 & -\frac{1}{2} & 0 & 0\\ 0 & \frac{1}{2} & 0 & \frac{1}{2} & 0 & 0 & 0 & 0 & 0 \end{array}\right].
    \end{aligned}
\end{equation}
The input is $u(t)$, and we observe the linear combination of the 1st and 3rd states to enhance the observability condition. Thus, the output is $x(t)+z(t)$. The above quadratic system \cref{eq:Lorenz} gives rise to chaotic dynamics for different choices of the parameters $(\sigma,~\rho,~\beta)$. This study aims to identify the Lorenz system from i/o time-domain data under harmonic excitation in the non-chaotic operational regime. We choose $\sigma=10,~\beta=8/3$, and for the parameter $\rho$, we investigate two cases 1,2 and comment on case 3, which exhibits chaotic dynamics.
\begin{enumerate}
    \item $\rho=0.5$, where the linear subsystem is stable, and the Lorenz attractor has only zero equilibrium.
    \item $\rho=20$, where the linear subsystem is unstable, the Lorenz system has two different steady-states operating around two non-zero stable equilibrium points \cref{fig:align}.
    \item $\rho=28$, where the linear subsystem is unstable, but the Lorenz system is chaotic (steady-state unreachable) with two non-trivial attractors.
\end{enumerate}
\textbf{Case 1 - $\rho=0.5$.}
 Samples of the first three GFRFs over the following frequency grids can be obtained from a physical measurement setup after processing the time-domain evolution of the potentially unknown system over the frequency domain as explained in \cref{sec:dataaquisition}. 
\begin{itemize}
    \item We take $100$ logarithmic distributed measurements $\omega_i,~i=1,\ldots,100$, from $2\pi[10^{-3},10^{3}]$. Therefore, $100$ pairs of measurements $\{j\omega_i,~H_1(j\omega_i)\}$ are collected. Using the Loewner framework \cref{sec:LF}, the order $r=3$ of the minimal linear subsystem can be identified from the singular value decay \cref{fig:fig1}(left), and a linear realization can be constructed:
    \begin{equation}
        \hat{\bA}=\left[\begin{array}{ccc} 17.1 & 4.08 & -4.51\\ 11.3 & 4.26 & -4.12\\ 140.0 & 30.1 & -35.0 \end{array}\right],~\hat{\bB}=\left[\begin{array}{c} -2.56\\ -1.53\\ -21.5 \end{array}\right],~\hat{\bC}^T=\left[\begin{array}{c} 3.12\\ 0.0985\\ -0.473 \end{array}\right].
    \end{equation}

The coordinate system is different from the original, but we can validate that the first transfer function $H_1$, or the Markov parameters, i.e., $\bC\bA\bB=\hat{\bC}\hat{\bA}\hat{\bB}=-12.6667$, and the DC term which is zero in that case remain the same. Note that the eigenvalues are: $\texttt{eig}(\bA)=\texttt{eig}(\hat{\bA})=(-10.52,-0.4751,-2.667)$.
\item We take $10$ logarithmically spaced distributed measurements from a squared grid in each dimension\footnote{Cartesian product: $(2\pi)^2[a,b]^2=2\pi[a,b]\times2\pi[a,b]$ for $a<b$.} as $(2\pi)^2[10^{-3},10^{3}]^2$, and $100$ pairs of measurements $\{(j\omega_1^{k},~j\omega_2^{k}),~H_2(j\omega_1^{k},~j\omega_2^{k})\}
$
are collected. Solving the linear system \cref{eq:H2LS} by minimizing the $2$-norm (least-squares), we estimate $\hat{\bQ}_0$ as
\begin{equation*}
    \footnotesize\hat{\bQ}_0=\left[\begin{array}{ccccccccc} -1.58 & 0.192 & 0.34 & 0.192 & -0.0195 & 0.0192 & 0.34 & 0.0192 & -0.0682\\ 0.183 & 1.0 & 0.042 & 1.0 & -0.312 & 0.0308 & 0.042 & 0.0308 & -0.009\\ 0.509 & 0.0395 & 1.08 & 0.0395 & 0.138 & 0.278 & 1.08 & 0.278 & -0.324 \end{array}\right].
\end{equation*}
\item The least-squares matrix in \cref{eq:H2LS} is rank deficient $\texttt{rank}(\bM)=21<27=3^3$. Therefore, a parameterization as in \cref{eq:Qsk} is introduced. In this particular case, the dimension of the vector $\blambda$ is six. As the proposed method is arbitrary in the number of measurements, we take $10$ logarithmically distributed measurements from the cubic grid $(2\pi)^3[10^{-3},10^3]^3$ in each dimension; $1,000$ triplets of measurements $\{(j\omega_1^{k},~j\omega_2^{k},~j\omega_3^{k}),~H_3(j\omega_1^{k},~j\omega_2^{k},~j\omega_3^{k})\}$ are collected. Solving the quadratic equation with \cref{alg:qve}, and starting with different seeds of $\lambda_0$, as it is depicted on the right of \cref{fig:fig1}, the parameter vector $\blambda\in\IR^6$ can be obtained in that case uniquely. Thus, the updated estimation of the quadratic operator $\hat{\bQ}=\hat{\bQ}_0+\sum_{i=1}^r\lambda\hat{\bQ}_i$, with $\hat{\bQ}_i$ the appropriate reshaped matrices from the null space vectors results to 
\begin{equation*}
\footnotesize\hat{\bQ}=\left[\begin{array}{ccccccccc} -2.38 & -0.184 & 0.391 & -0.184 & -0.00636 & 0.0304 & 0.391 & 0.0304 & -0.0642\\ 26.0 & 0.497 & -3.98 & 0.497 & 0.00948 & -0.0758 & -3.98 & -0.0758 & 0.605\\ -11.3 & -1.05 & 1.89 & -1.05 & -0.0372 & 0.174 & 1.89 & 0.174 & -0.315 \end{array}\right].
\end{equation*}

\begin{figure}[!ht]
    \centering
    \includegraphics[scale=0.14]{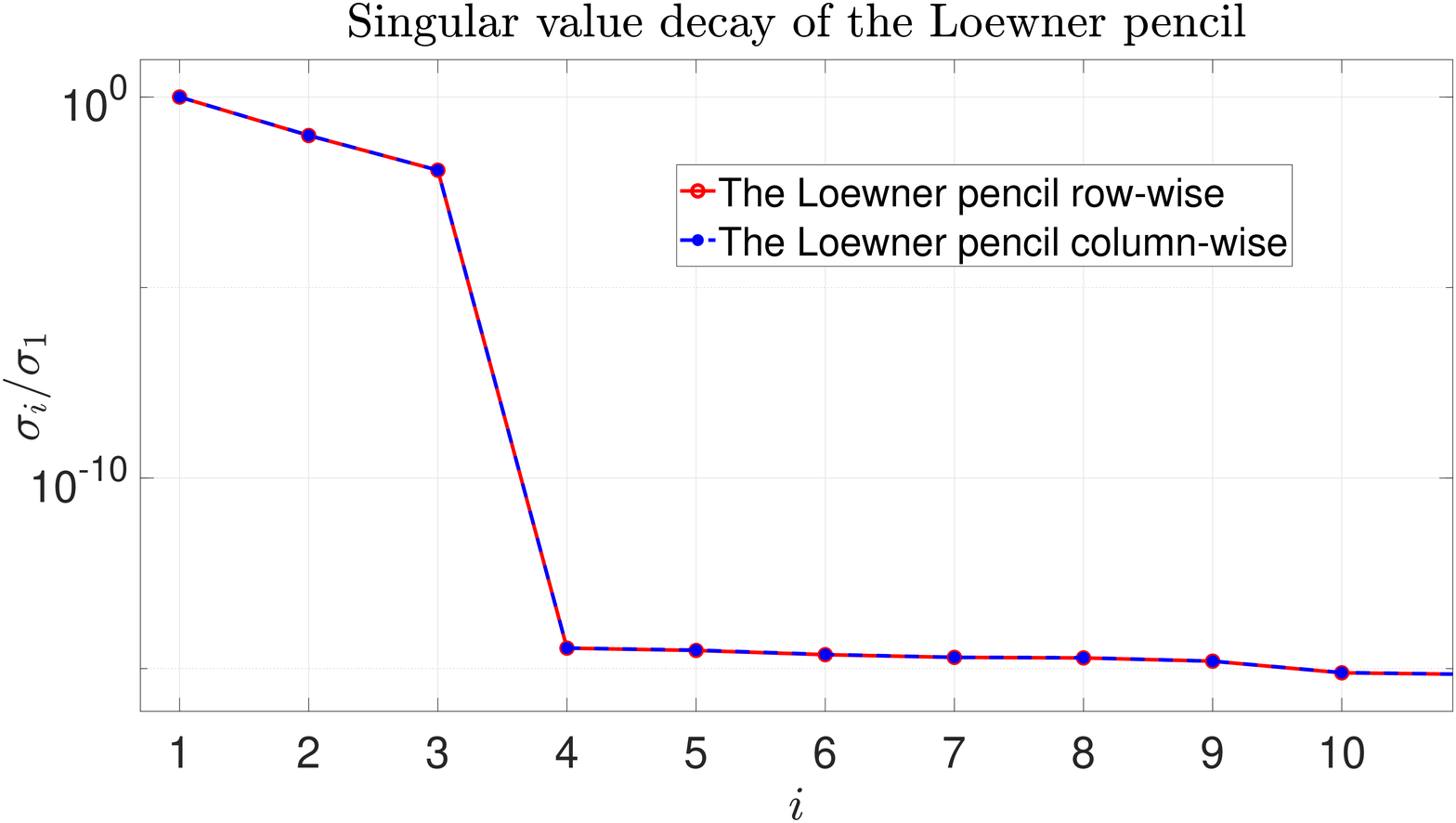}
    \includegraphics[scale=0.14]{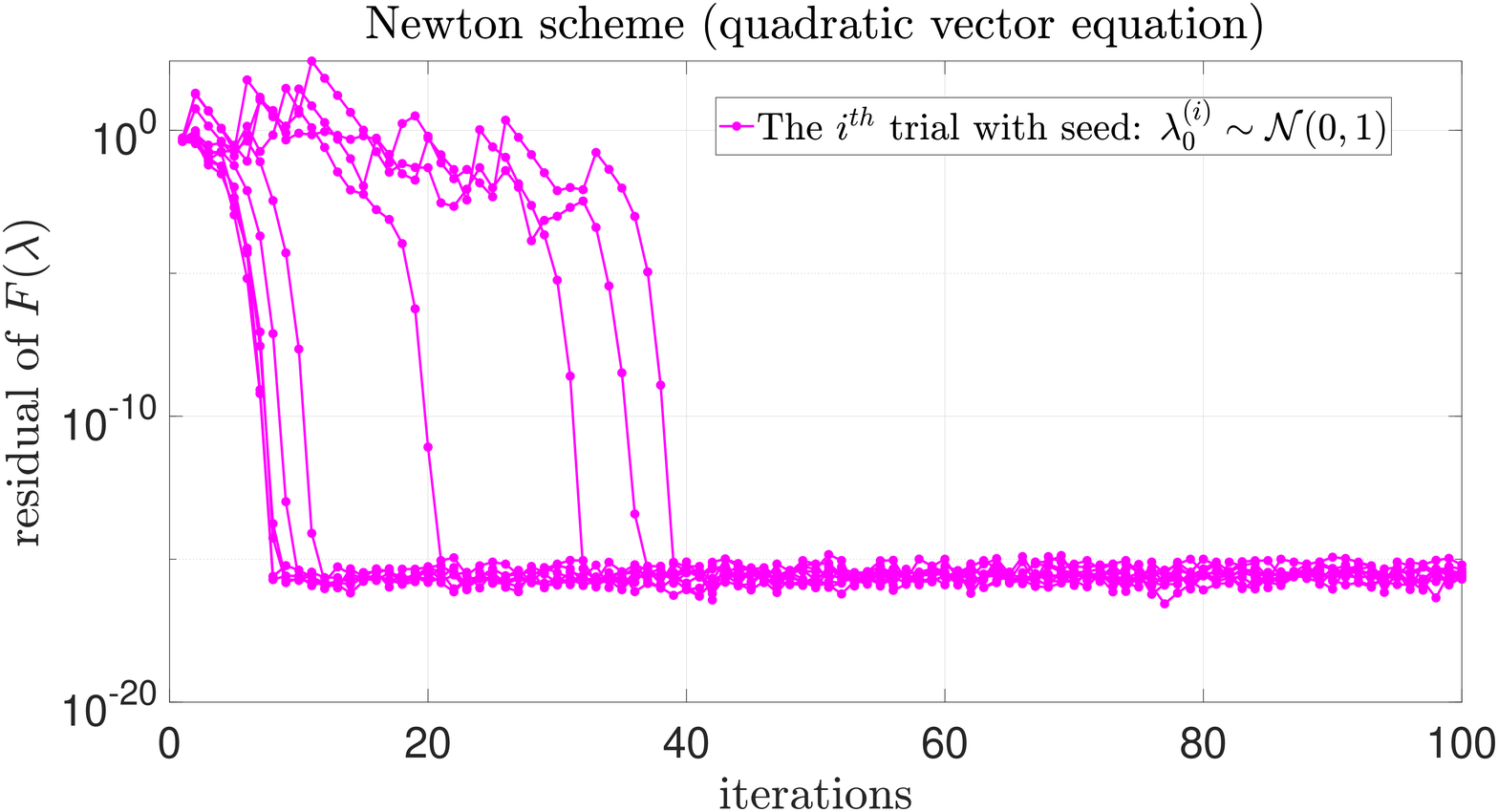}
    \caption{\textbf{Left}: The Loewner singular value decay, discovered minimal order $r=3$ with $\sigma_4/\sigma_1\sim 1e-14$. \textbf{Right:} The Newton convergence scheme \cref{alg:qve}. Solution vector $\blambda^*$ can be obtained uniquely after starting with different random seeds $\blambda_0$.}
    \label{fig:fig1}
\end{figure}

Finally, aligning the identified system to the original coordinates shows that the two systems are identical after computing a transformation $\bPsi$ from the \cref{app:align}.
\small
\begin{equation*}
\begin{aligned}
    \bA&=\bPsi^{-1}\hat{\bA}\bPsi=\left[\begin{array}{ccc} -10.0 & 10.0 & 4.334e-13\\ 0.5 & -1.0 & 5.219e-13\\ -5.254e-11 & 5.448e-11 & -2.667 \end{array}\right],\\
    \bB&=\bPsi^{-1}\hat{\bB}=\left[\begin{array}{c} 1.0\\ -2.207e-11\\ 1.0 \end{array}\right],~\bC=\hat{\bC}\bPsi=\left[\begin{array}{ccc} 1.0 & -1.806e-11 & 1.0 \end{array}\right],\\
    \bQ&=\bPsi^{-1}\hat{\bQ}(\bPsi\otimes\bPsi)=\left[\begin{array}{ccccccccc} 0 & 0 & 0 & 0 & 0 & 0 & 0 & 0 & 0\\ 0 & 0 & -0.5 & 0 & 0 & 0 & -0.5 & 0 & 0\\ 0 & 0.5 & 0 & 0.5 & 0 & 0 & 0 & 0 & 0 \end{array}\right] \pm\epsilon\cdot\mathbf{1},
    \end{aligned}
\end{equation*}

where $\epsilon\in[1e-12,1e-10]$ and $\mathbf{1}\in\IR^{3\times 9}$. The above result certifies that the original and identified systems are equivalent under a coordinate transformation. Since $\bQ\neq \bPsi^{-1}\hat{\bQ}_s(\bPsi\otimes\bPsi)$, quadratic identification with information from the first two kernels $H_1,~H_2$ is impossible. Here, the significant improvement compared with other similar studies \cite{KARACHALIOS20227} is the systematic way to add more information to the constructed model from the higher-order GFRFs, i.e., $H_3$. As a result, the forced Lorenz system was successfully identified when measurements of the first three symmetric GFRFs were considered, as it is illustrated in \cref{fig:fig2}(left) in contrast with the unstable result obtained with information available only from the first two $H_1,~H_2$. Finally, in \cref{fig:fig2}(right), the identified and the original systems are equivalent state-space models after transforming them to the same coordinate system. 
\end{itemize}

\begin{figure}[!ht]
    \centering
    \includegraphics[scale=0.14]{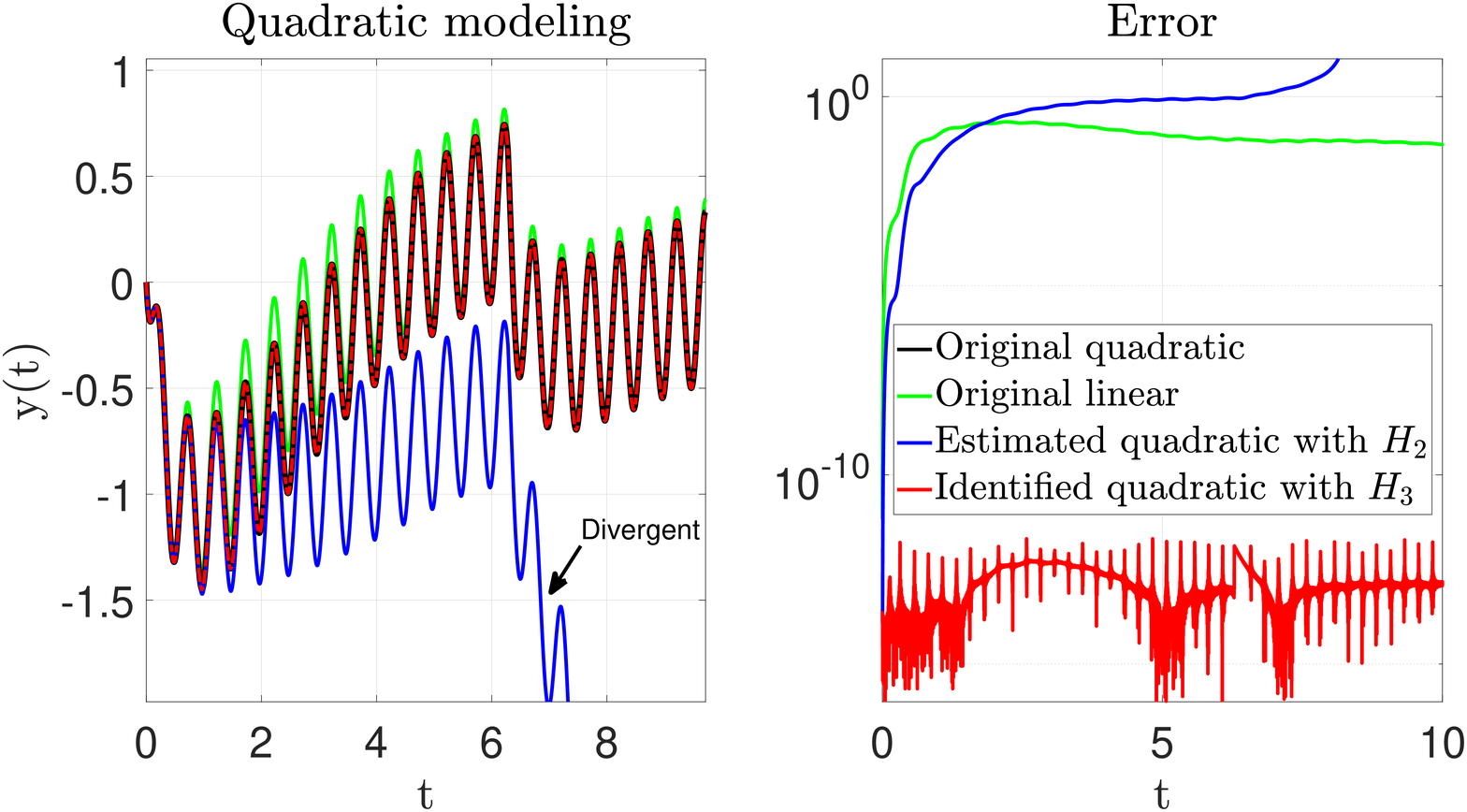}
    \includegraphics[scale=0.14]{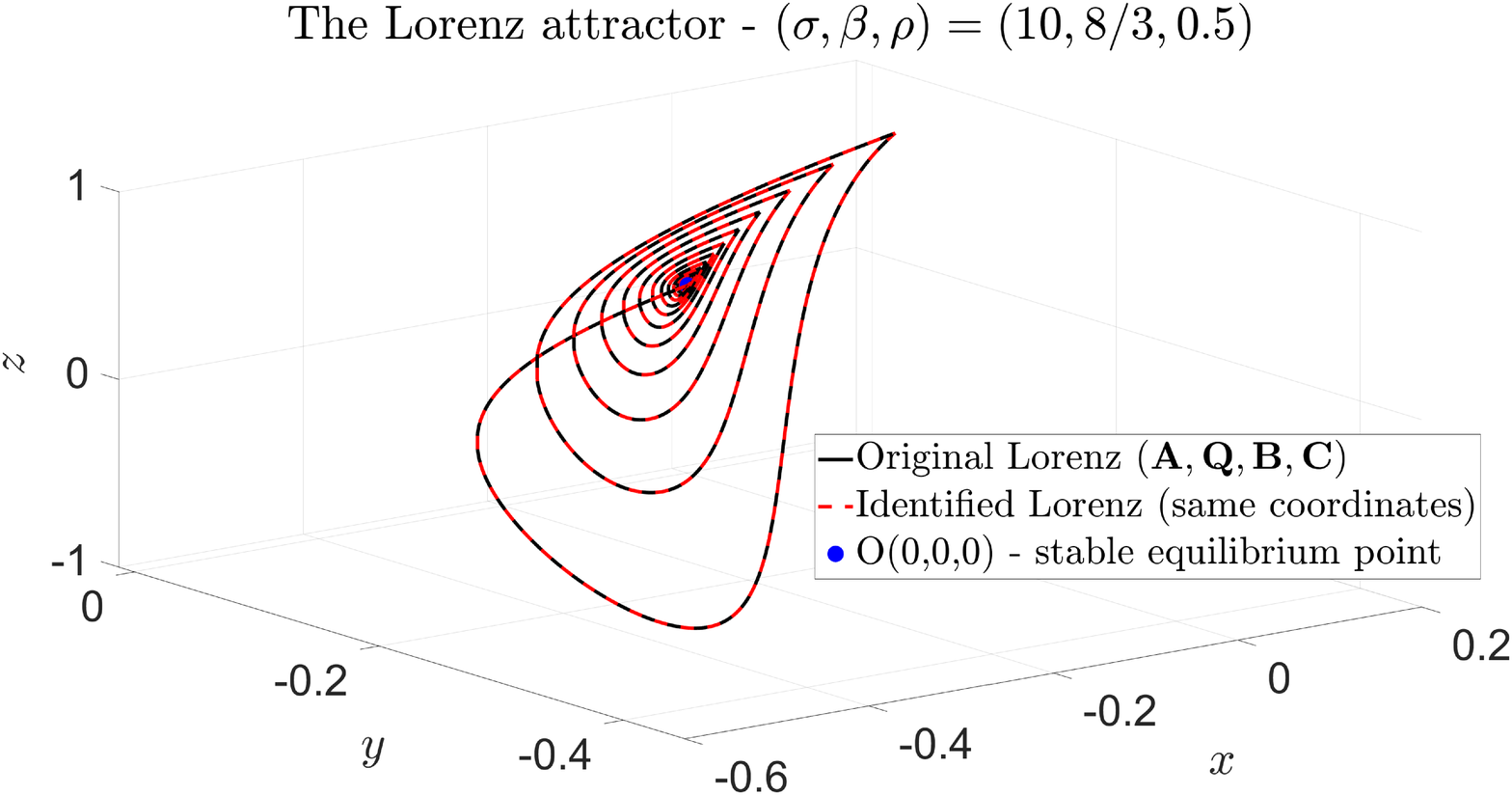}
    \caption{\textbf{Left}: The linear model gives a poor approximation. Also, the $H_2$ does not contribute to a reasonable estimation of the quadratic operator; therefore, numerical instability is observed. After incorporating the information from the 3rd kernel, the Lorenz system was identified with a numerical error near machine precision. \textbf{Right}: The 3D state space is reconstructed from the identified system with the proposed method and is compared to the original one after aligning both systems to the same coordinates.}
    \label{fig:fig2}
\end{figure}
\textbf{Case 2 - $\rho=20>1$.} In that case, the Lorenz attractor has two non-zero equilibrium points, i.e., $\bx_e^{(1)}=\left[\begin{array}{ccc}
-\sqrt{\beta(\rho-1)} & -\sqrt{\beta(\rho-1)} & \rho-1\end{array}\right]^T$ and $\bx_e^{(2)}=\left[\begin{array}{ccc}
\sqrt{\beta(\rho-1)} & \sqrt{\beta(\rho-1)} & \rho-1\end{array}\right]^T$. 

\begin{figure}[!ht]
    \centering
    \includegraphics[scale=0.14]{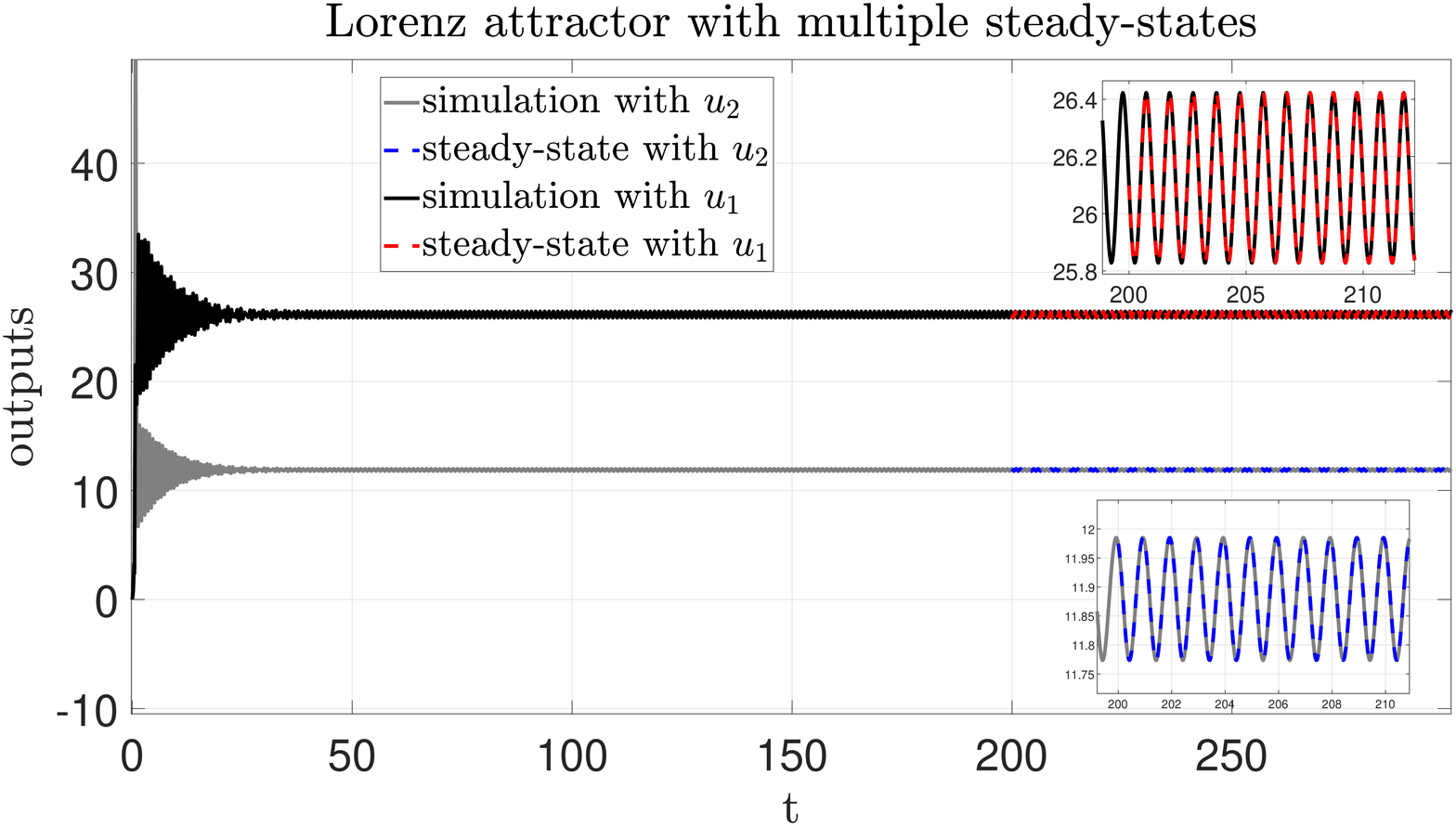}
    \includegraphics[scale=0.14]{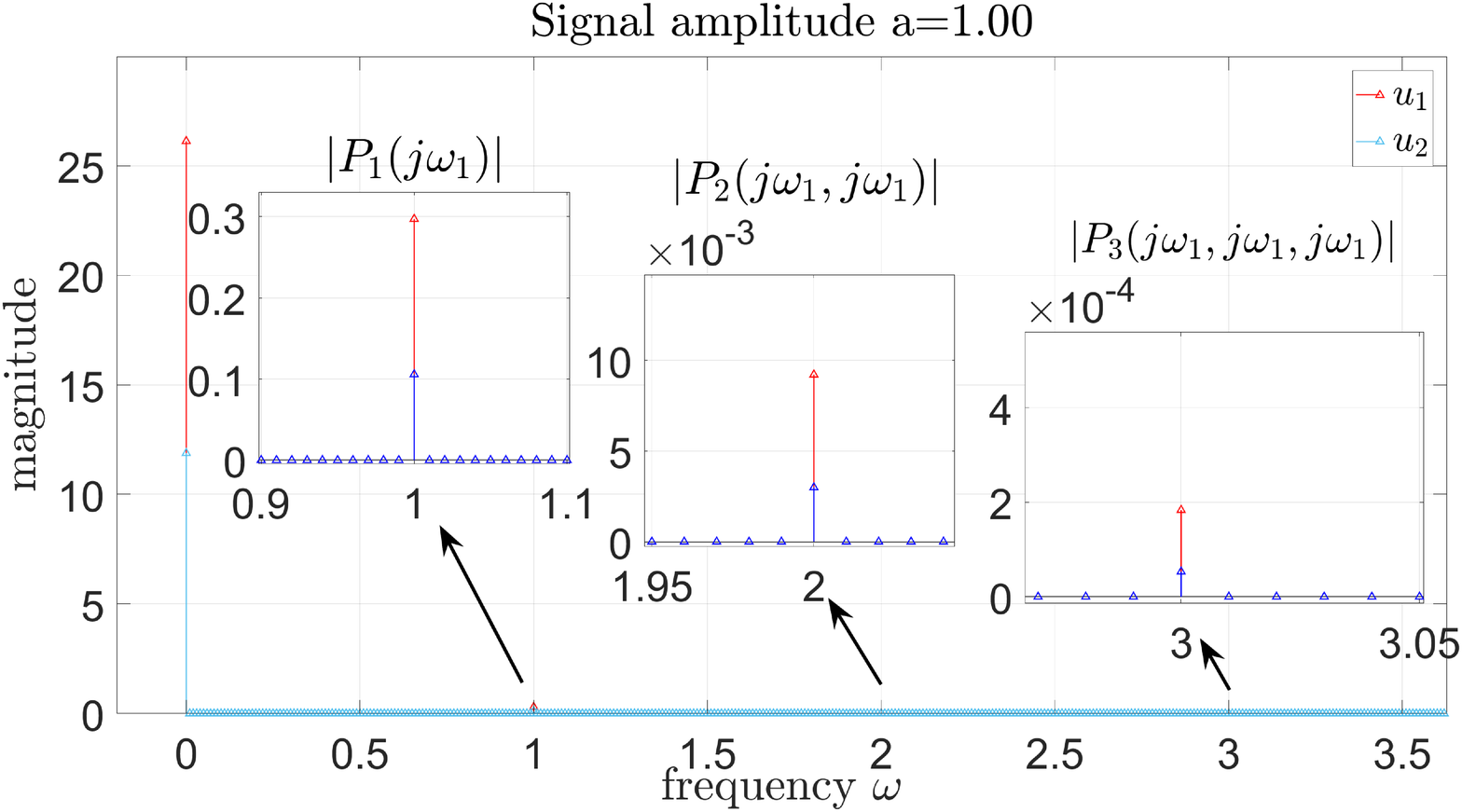}
    \caption{Multiple steady-states and corresponding power spectrum. Numerical details in the \cref{tab:meastab}.}
    \label{fig:figspectrum}
\end{figure}

Under excitation or non-zero initial conditions, the system's trajectories spiral around these two attractors.\\
\textbf{Data assimilation over multiple steady-states.}
In \cref{tab:meastab}, we illustrate how measurements of the higher-order GFRFs can be obtained via harmonic input excitation of a given control system. Here, we illustrate the experimental process by exciting with harmonic inputs the Lorenz attractor with the parameter $\rho=20$. Typically, bifurcation analysis concerns the system's sensitivity to parameters, which leads to structural changes in the system dynamics. In our case, we consider system parameters resulting from the input term, which can potentially drive the system to such structurally different regimes. For instance, with $\alpha=1,~\omega_1=1$, we have two different designed complex inputs\footnote{With complex inputs, e.g., $u(t):\IR_+\rightarrow\IC$, indexing harmonics and estimating kernels are straightforward tasks compared to the real input case, e.g., $u(t):\IR_+\rightarrow\IR$, where additional operations s.a. kernel separation with amplitude shifting, should be addressed \cite{morKarGA21a,Boyd1983MeasuringVK,Xpar2006}.} with identical asymptotic behavior as the perturbation signal decaying to zero by the time the system reaches a quasi-steady state. These inputs are: 
\begin{enumerate}
    \item $\rho=20$ - input 1: $u_1(t)=\underbrace{3 e^{-0.1 t}\texttt{sawtooth}(t)}_{\text{perturbation}}+\alpha e^{2j\pi\omega_1 t}$.
    \item $\rho=20$ - input 2: $u_2(t)=\alpha e^{2j\pi\omega_1 t}$.
\end{enumerate}
As it is depicted in \cref{fig:figspectrum}(left), for the different designed inputs $u_1(t),~u_2(t)$, we obtain two different quasi-steady state solutions with a different power spectrum \cref{fig:figspectrum}(right). Measurements can be obtained for both systems, and the different DC terms will help distinguish the equilibrium operation.  

\begin{table}[h!]
        \centering
        \begin{tabular}{c|cccc}
          \textbf{Data}	& DC & $H_1(j\omega_1)$	 & $H_2(j\omega_1,j\omega_2)$ & $H_3(j\omega_1,j\omega_2,j\omega_3)$\\\hline
          $u_1(t)$ & $11.8819$ & $-0.0148+0.297{}\mathrm{i}$ & $-0.00687-0.00614{}\mathrm{i}$ & $1.74e-4-5.82e-5{}\mathrm{i}$\\
          $u_2(t)$ & $26.1181$ & $0.09303+0.05011{}\mathrm{i}$ & $-3.0e-4-3.0e-3{}\mathrm{i}$ & $6.0e-6+5.3e-5{}\mathrm{i}$\\
        \end{tabular}
        \caption{In this table, and for each system (blue, red), the single-sided Fourier spectrum (magnitude, phase) $P$ provides the following measurements (complex inputs): $H_1(j\omega_1)=P_1(j\omega_1)/a$, $H_2(j\omega_1,j\omega_1)={P_2(j\omega_1,j\omega_1)}/{a^2}$, $H_3(j\omega_1,j\omega_1,j\omega_1)={P_3(j\omega_1,j\omega_1,j\omega_1)}/{a^3}$. For each system, the non-periodic DC term, e.g., \cref{eq:dcterms}, can be computed from the non-periodic value $P(0)$ in the power spectrum. The measurement set-up is generalizable for multi-harmonic real signals (kernel separation) and in the $X$-parameter machinery \cite{Xpar2006} that deals with harmonic distortion.}
        \label{tab:meastab}
    \end{table}

The symmetric GFRFs that can interpret the measurements in \cref{tab:meastab} are related with the corresponding linear operator, i.e., $\tilde{\bA}_{q}=\bA+2\bQ(\bx_e^{(q)}\otimes\bI),~q=1,2$ for each equilibrium respectively. For instance, using the equilibrium $\bx_e^{(1)}$, we compute $\tilde{\bA}_1=\bA+2\bQ(\bx_e^{(1)}\otimes\bI)$. The 1st transfer function $H_1$ (for the equilibrium $\bx_e^{(1)}$) yields the following value at frequency $\omega_1=2\pi$: $H_1(j\omega_1)=\bC(j\omega_1\bI-\tilde{\bA})^{-1}\bB=-0.0148+0.297\mathrm{i}$ which explains the measurements in the 1st row of \cref{tab:meastab}, and similarly, the higher GFRFs explain the subsequent rows of \cref{tab:meastab}. Similar results can be obtained for the 2nd input-$u_2$. Importantly, when the Lorenz system is concerned with $\rho=20>1$, we measure two different local quadratic systems that can be recognized/distinguished from the two different DC terms in \cref{fig:figspectrum}. Thus, the proposed method can identify the two deduced quadratic systems \cref{eq:syss}. 

With the proposed method \cref{alg:Qmodel} and measurement setup as in the previous case 1, we can identify the two quadratic models that interpret the dynamics locally \cref{eq:syss} and after computing the state equilibrium from \cref{eq:L0} the original operators $A_q,~q=1,2$, remain from \cref{eq:syssq}; we obtain the system $(\tilde{\bA}_1,~\tilde{\bQ}_1,~\tilde{\bB}_1,\tilde{\bC}_1)$ around a non-zero equilibrium as
\begin{equation}
\footnotesize
    \begin{aligned}
\tilde{\bA}_1&=\left[\begin{array}{ccc} -0.96 & -8.79 & 5.51\\ 2.12 & 2.17 & 21.0\\ 3.14 & -4.69 & -14.9 \end{array}\right],~\tilde{\bB}_1=\left[\begin{array}{c} 0.719\\ 2.72\\ -2.33 \end{array}\right],\tilde{\bC}_1^T=\left[\begin{array}{c} -0.0479\\ 0.0361\\ -0.831 \end{array}\right],~\bx_e^{(1)}=\left[\begin{array}{c} 6.78\\ 51.7\\ -12.4 \end{array}\right],\\
        \tilde{\bQ}_1&=\left[\begin{array}{ccccccccc} -0.0853 & 0.212 & 0.236 & 0.212 & -0.524 & -0.596 & 0.236 & -0.596 & 0.00706\\ -0.11 & 0.156 & 0.259 & 0.156 & -0.107 & -0.00667 & 0.259 & -0.00667 & 0.525\\ 0.0331 & -0.0303 & -0.0709 & -0.0303 & -0.0478 & -0.105 & -0.0709 & -0.105 & -0.229 \end{array}\right],\\
        \bA_1&=\tilde{\bA}_1-2\tilde{\bQ}_1(\bx_e^{(1)}\otimes\bI)=\left[\begin{array}{ccc} -15.8 & 27.7 & 64.1\\ -6.11 & 10.9 & 31.3\\ 4.06 & -1.94 & -8.78 \end{array}\right],~\texttt{eig}(\bA_1)=\{-20.3408,~-2.6667,~9.3408\}.
    \end{aligned}
\end{equation}
Likewise, for completeness, the identified quadratic model around the second equilibrium is
\begin{equation}
\footnotesize
    \begin{aligned}
        \tilde{\bA}_2&=\left[\begin{array}{ccc} -0.21 & 9.09 & 6.02\\ -7.47 & 2.77 & -19.0\\ -1.07 & 2.52 & -16.2 \end{array}\right],~\tilde{\bB}_2=\left[\begin{array}{c} 0.712\\ -2.61\\ -2.46 \end{array}\right],\tilde{\bC}_2^T=\left[\begin{array}{c} -0.706\\ -0.698\\ -0.278 \end{array}\right],~\bx_e^{(2)}=\left[\begin{array}{c} -35.3\\ -7.22\\ 13.7 \end{array}\right],\\
        \tilde{\bQ}_2&=\left[\begin{array}{ccccccccc} 0.136 & 0.209 & 0.108 & 0.209 & -0.612 & 0.527 & 0.108 & 0.527 & -0.0535\\ -0.041 & -0.0666 & -0.0296 & -0.0666 & 0.156 & -0.146 & -0.0296 & -0.146 & 0.0148\\ -0.16 & -0.313 & -0.0692 & -0.313 & 0.161 & -0.364 & -0.0692 & -0.364 & 0.0378 \end{array}\right],\\
        \bA_2&=\tilde{\bA}_2-2\tilde{\bQ}_2(\bx_e^{(2)}\otimes\bI)=\left[\begin{array}{ccc} 9.42 & 0.531 & 22.7\\ -10.5 & 4.32 & -23.6\\ -14.9 & -7.21 & -27.4 \end{array}\right],~\texttt{eig}(\bA_2)=\{-20.3408,~-2.6667,~9.3408\}.
    \end{aligned}
\end{equation}

\begin{figure}[h!]
    \centering
    \includegraphics[scale=0.14]{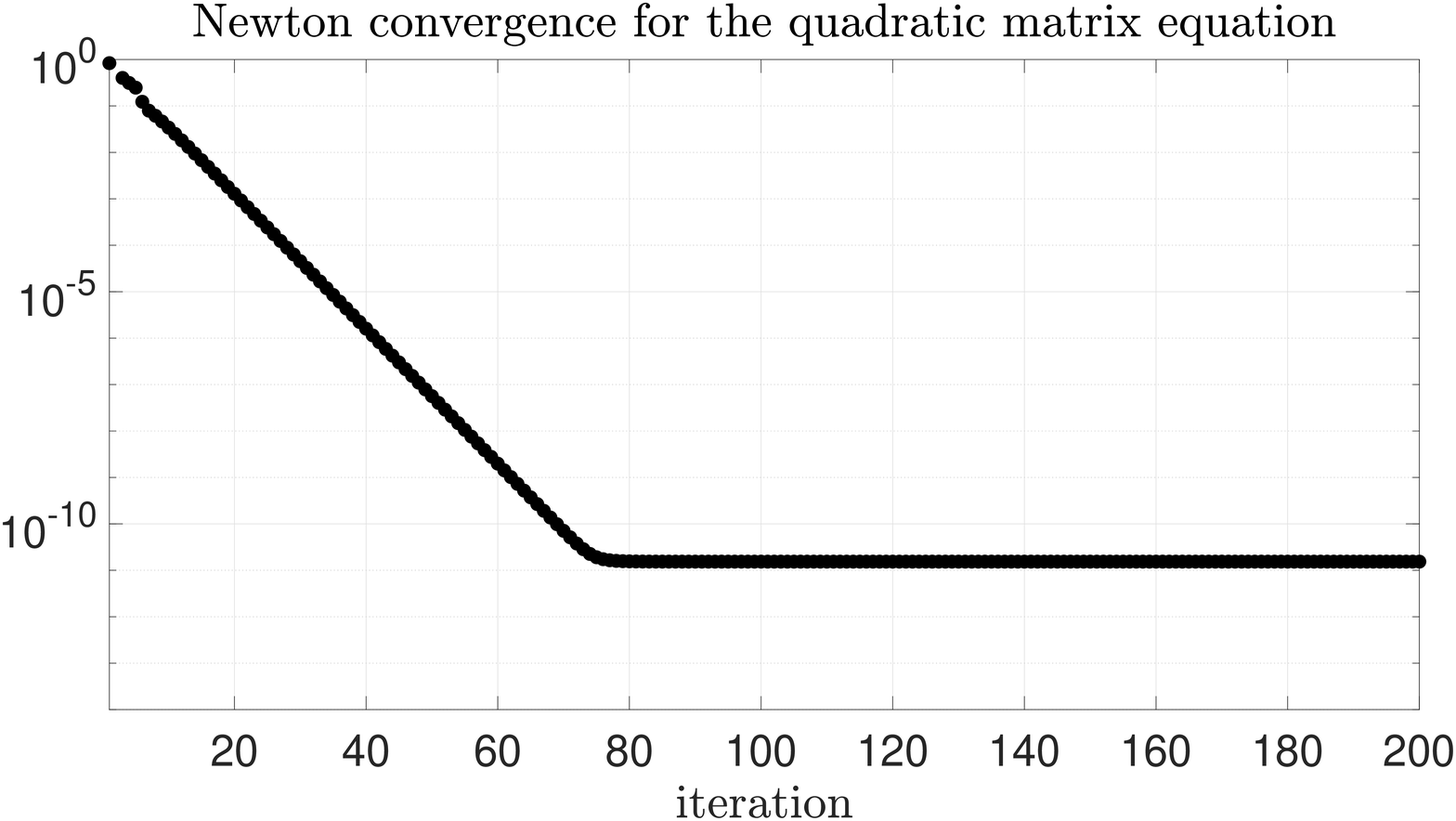}
    \includegraphics[scale=0.14]{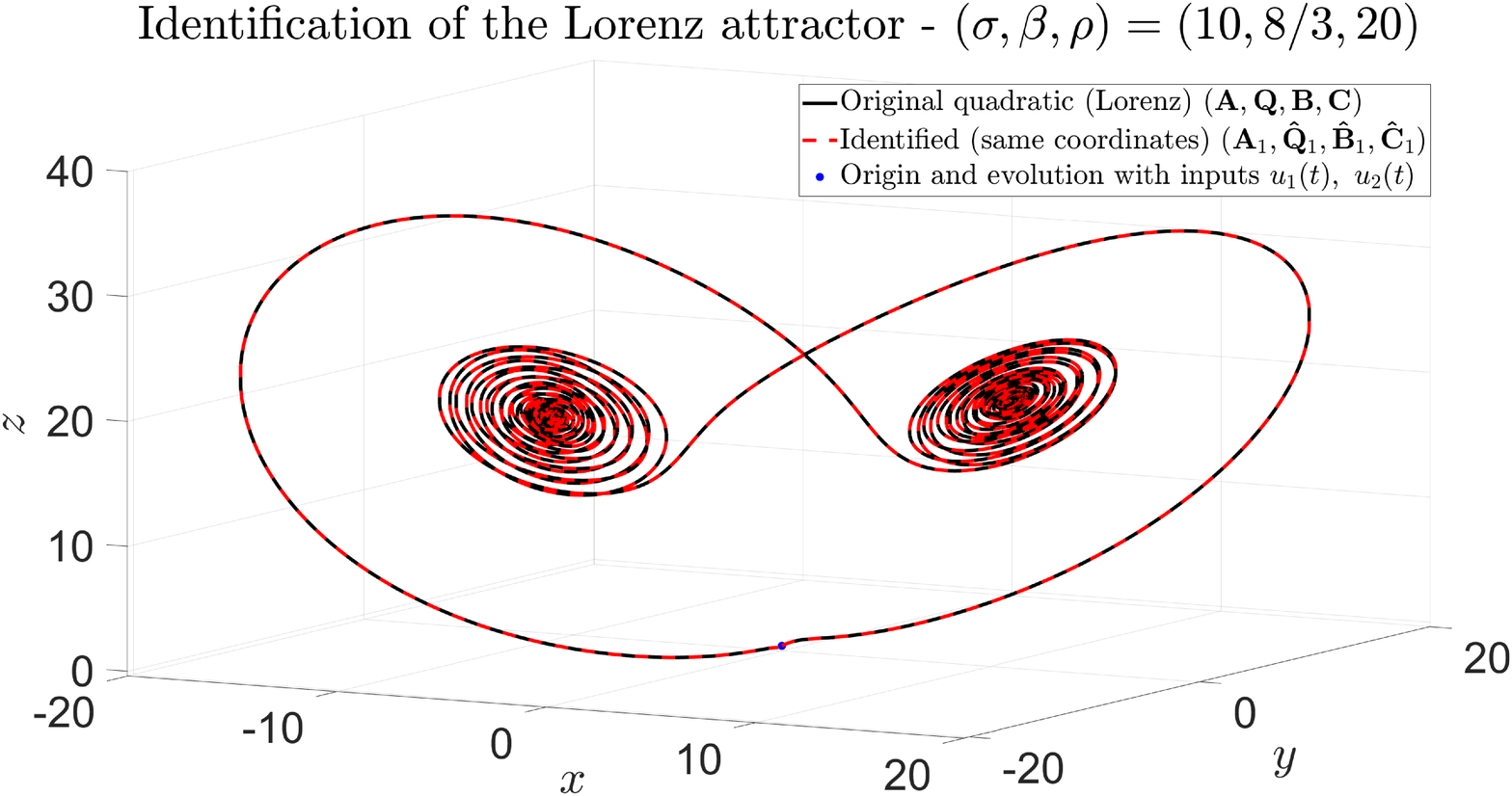}
    \caption{\textbf{Left}: Convergence of the Newton scheme in \cref{alg:qme}. \textbf{Right}: The original Lorenz system is identified with the two non-trivial equilibrium points. The comparison between the original system and the identified one is depicted. The constructed state space evolution for both systems and at the same coordinates remains the same with zero numerical error.}
    \label{fig:fig4}
\end{figure}

We discovered the original global quadratic model, i.e., the Lorenz system, independently from the local inference around the two equilibria. Thus, to illustrate the identification result of the global quadratic dynamics, we choose one of the two models, and for various inputs, we observe that the dynamics evolve around two non-trivial attractors. To validate our result further, we align one of the above models with the original Lorenz system and compare them in \cref{fig:fig4}.

As this study is devoted to quadratic modeling, another way to discover the original model and the two equilibrium points simultaneously can be obtained after first utilizing \cref{lem:2.6} for aligning the system's operators to the same coordinates without using information from the corresponding global linear matrices and then proceeding to find the equilibrium points as a coupled quadratic vector equation explained in \cref{app:align}. For the alignment, we start with a random seed for the Newton method \cref{alg:qme} e.g., $\bT_0\sim\mathcal{N}(\bmu,\bsigma)$, and applying \cref{alg:qme}. The following convergence is depicted in \cref{fig:fig4}, and the transformation is
\begin{equation}\label{eq:transformQBC}
    \bT=\left[\begin{array}{ccc} -0.613 & 1.98 & -2.57\\ -1.79 & 0.375 & -2.02\\ 0.807 & 0.742 & 0.395 \end{array}\right],~\bT^{-1}=\left[\begin{array}{ccc} 1.22 & -1.99 & -2.25\\ -0.686 & 1.36 & 2.49\\ -1.21 & 1.52 & 2.46 \end{array}\right]
\end{equation}
The two quadratic systems have been aligned with the transformation $\bT$ without the linear matrices. 


\begin{remark}[Global system identification from local measurements]
    After measuring the local dynamics at least over an achievable stable equilibrium that provides a quasi-steady state under harmonic excitation, we identified the original global Lorenz system with the unstable linear operator with the proposed method. After transforming the operators $(\bA_1,\tilde{\bQ}_1,\tilde{\bB}_1,\tilde{\bC}_1)$ to the original coordinates utilizing \cref{app:align}, the two systems; original \cref{eq:Lorenz} and identified are the same \cref{fig:fig4}(right) to any operational regime. Moreover, after appropriate input excitation of the discovered system, the dynamics can switch across the different attractors by realizing the global dynamical behavior as illustrated in \cref{fig:fig4}(right).
\end{remark}

\textbf{ Case 3 - $\rho=28$.}  Since the dynamics are chaotic for this parameter range (when $\rho>24.74$), a steady state cannot be achieved, and measurements cannot be obtained for estimating GFRFs with the proposed method. Identifying such systems is beyond this study's scope, and the Volterra series would need a new formulation to address such as chaotic dynamics or dynamics with a continuous spectrum. Generally, operators can be identified for any identification method to finite numerical precision (e.g., IEEE machine precision $\epsilon\approx 2.22e-16$). Since this is already an approximation with a nonzero numerical error, this slight numerical discrepancy will not allow any safe prediction of such a sensitive system as the (deterministic) chaotic Lorenz attractor. Other approaches, mainly stochastic toward dealing with dynamics in the chaotic regime, with continuous spectrum, can be found in \cite{LuschBruntonKutz,Ting2023,Pathak18,PhysRevE.91.032915,Chorin2015}.

\subsection{Data-driven quadratic inference of the viscous Burgers' model}
In this example, we want to illustrate the capability of the proposed method in discovering reduced models from a potentially unknown system (large-scale) from i/o time-domain measurements. In particular, the aim is to construct robust surrogate models of reduced order from i/o measurements that can be obtained physically without any prior model (i.e., samples of the symmetric GFRFs obtained from time-domain simulations). These measurements could be obtained from a real engineering process (e.g., X-parameters \cite{Xpar2006}). We consider the Burgers' equation with Robin boundary conditions from \cite{morAntGH19}. Given viscosity $\nu>0$ and parameters $\sigma_0\leq 0,~\sigma_1\geq 0$, consider
\begin{equation}\label{eq:BurgersPDE}
    \begin{aligned}
        &\frac{\partial}{\partial t}\upsilon(x,t)-\nu\frac{\partial^2}{\partial x^2}\upsilon(x,t)+\upsilon(x,t)\frac{\partial}{\partial x}\upsilon(x,t)=0,~x\in(0,1),~t\in(0,T),\\
        &\nu\frac{\partial}{\partial x}\upsilon(x,0)+\sigma_0\upsilon(0,t)=u_0(t),~t\in(0,T),\\
        &\nu\frac{\partial}{\partial x}\upsilon(x,1)+\sigma_1\upsilon(1,t)=0,~t\in(0,T),\\
        &\upsilon(x,0)=0,~x\in(0,1).
    \end{aligned}
\end{equation} 
The input is the function $u_0(t)$ and the output is $y(t)=\int_0^1\upsilon(x,t)dx$.

The Burgers model is one challenging benchmark due to its mathematical nature, with nonlinear advection phenomena dominating for low viscosity values. Methods in dealing with this problem's stochastic version and predicting viscous shocks exist as in \cite{Chen22}. We keep the same model setup in \cite{morAntGH19} with a different viscosity parameter $\nu$. Here, we consider the flow's mean velocity field as an output. As illustrated in the study \cite{morAntGH19}, the Loewner models for small viscosity coefficients $\nu$ may produce unstable results. As the current study relies on the Volterra series representation, analysis of convergence with arbitrary viscosity and input amplitude remains an open issue. Therefore, we illustrate a more conservative case with a higher viscosity.

We use the problem data $\nu=0.5,~\sigma_0=0,~\sigma_1=0.1$ that represent the same physical quantities as in \cite{morAntGH19}. The full-order model (FOM) is the linear semi-discretization of \cref{eq:BurgersPDE} with finite elements of dimension $n=257$. The semi-discretized system can result in \cref{eq:qsys} after inverting the well-conditioned mass matrix $\bE$. The system solution \cref{eq:qsys} is approximated with the Runge-Kutta multistep integration method with a uniform time-discretization step $dt=1/1000$ with the Matlab command \texttt{ode45}.
Similarly, as in the Lorenz example, we take the following measurements:
\begin{itemize}
    \item $100$ logarithmic distributed measurements from the interval $2\pi[10^{-2},10^{1}]$,
    \item $400$ logarithmic distributed measurements from the square grid $(2\pi)^2[10^{-2},10^{1}]^2$,
    \item $216$ logarithmic distributed measurements from the cubic grid $(2\pi)^3[10^{-2},10^{1}]^3$.
\end{itemize} 

\begin{figure}[h!]
    \centering
    \includegraphics[scale=0.25]{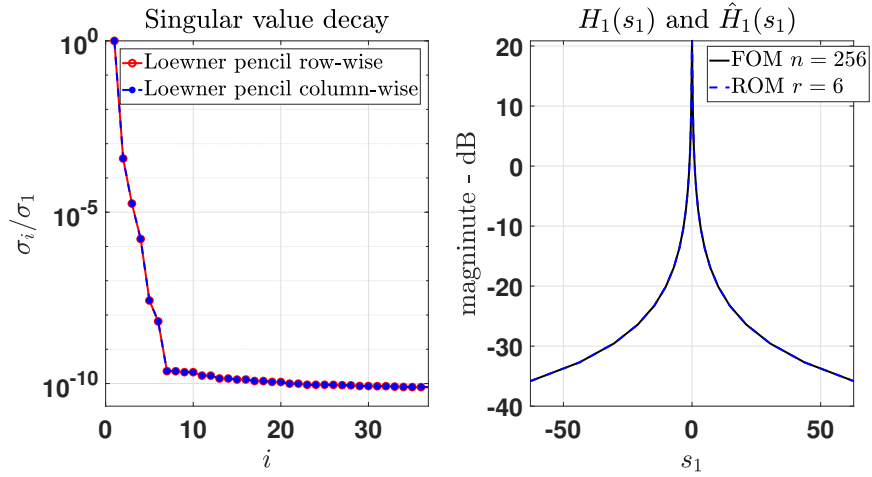}
    \caption{Left: the singular value decay of the Loewner matrix. Right: Comparison between the FOM and ROM kernels of the first level with error $\lVert\bH_1(\bs_1)-\hat{\bH}_1(\bs_1)\rVert=9.14e-08$.}
    \label{fig:BurgersLFdecay}
\end{figure}

In \cref{fig:BurgersLFdecay} (left), the singular value decay of the Loewner matrix pencil is presented. Therefore, we choose the minimal linear order $r=6$ with the first normalized truncated singular value of magnitude $\sigma_7/\sigma_1=2.3069\cdot 1e-10$. The recovery of the 1st GFRF $H_1$ that results from the FOM of dimension $n=257$ is compared with the reduced $\hat{H}_1$ of dimension $r=6$ in \cref{fig:BurgersLFdecay}(right).

\begin{figure}[h!]
    \centering
    \includegraphics[scale=0.14]{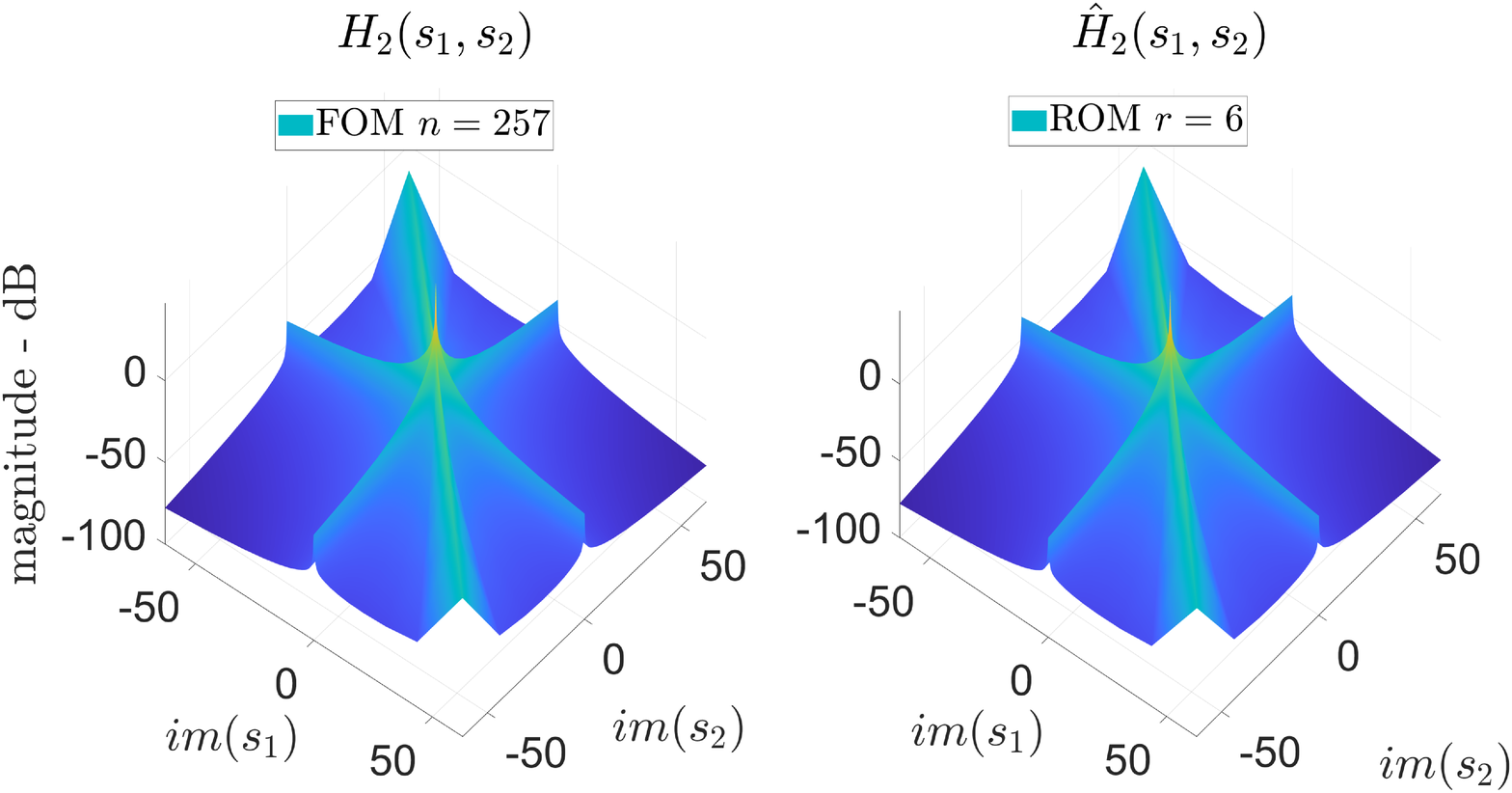}
     \includegraphics[scale=0.14]{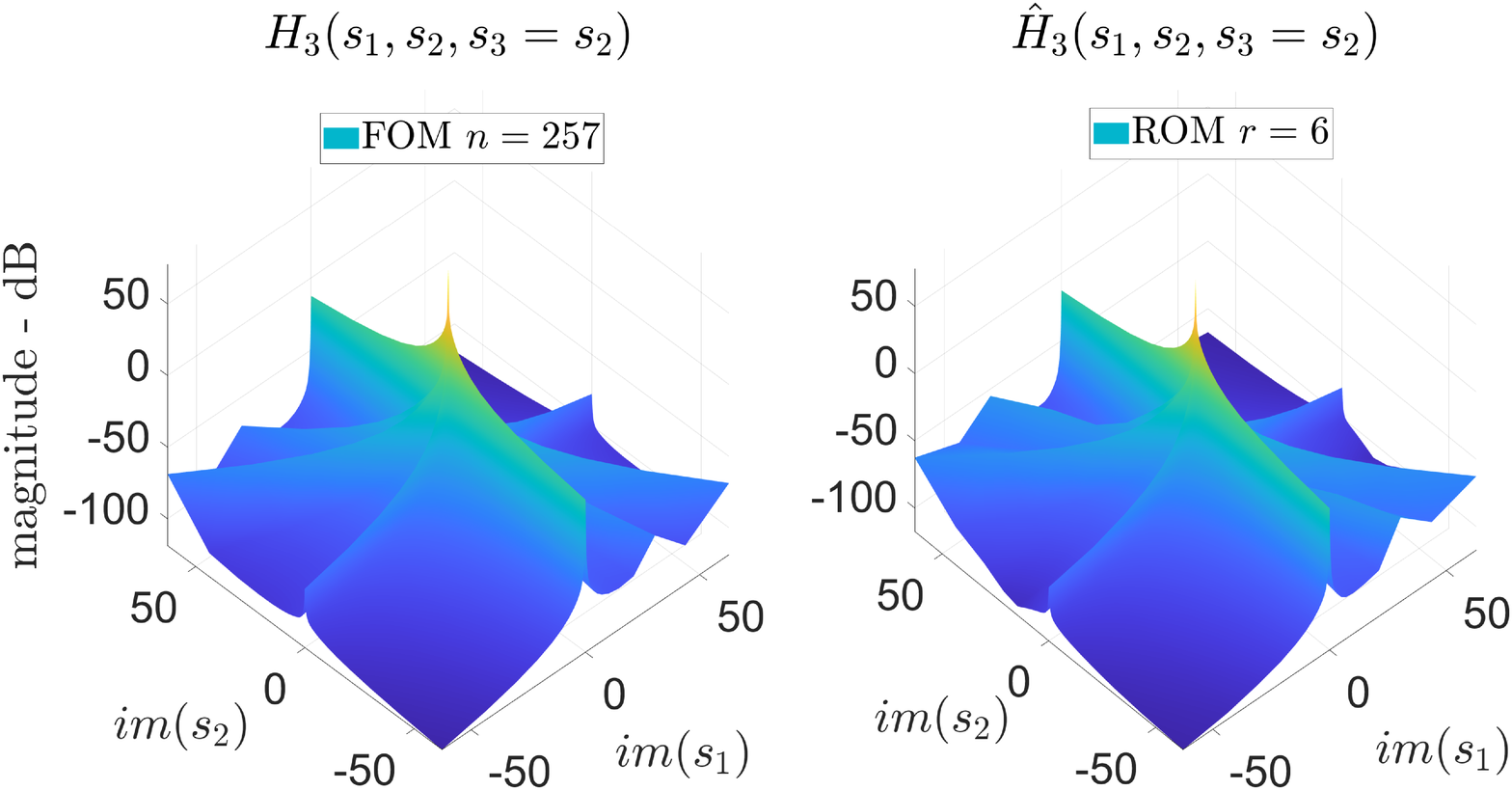}
    \caption{Comparison between the FOM and ROM; On the left, for the 2nd GFRF $\lVert \bH_2(\bs_1,\bs_2)-\hat{\bH_2}(\bs_1,\bs_2)\rVert_2=2.21e-02$, and on the right, for the 3rd GFRF. The error over the 2D plain-domain after fixing the 3rd dimension as $s_2=s_3$ is $\lVert\bH_3(\bs_1,\bs_2,\bs_2=\bs_3)-\hat{\bH}_3(\bs_1,\bs_2,\bs_2=\bs_3)\rVert_2=4.0e-03$.}
    \label{fig:BurgersH2H3}
\end{figure}

After solving \cref{eq:H2LS} with a threshold $\eta=1e-9$, we obtain the quadratic operator $\bQ_s\in\IR^{6\times6^2}$. In a classical regularization sense, the hyperparameter is tuned to leverage the error with the mathematical norm of $\lVert\bQ\rVert_{2}$. Different ways to find the optimal regularization parameter $\eta$, e.g., Tikhonov regularization, L-curve, work similarly to the thresholding SVD. Furthermore, the choice of $\eta$ affects the null space size. In particular, we have $r^3=6^3=216$ degrees of freedom, and after enforcing the symmetries of the quadratic operator, the maximum rank could be $\texttt{rank}=141<216$ when $\eta$ is close to machine precision. Therefore, by inverting with the threshold $\eta=1e-9$, the rank is $117$, and the resulting null space has length $216-117=99$. These extra $99$ free parameters will also be estimated to interpolate the third GFRF $(i.e., H_3)$. The accuracy of fitting the higher GFRFs for the FOM and ROM systems is depicted in \cref{fig:BurgersH2H3}. 

We estimate the quadratic operator denoted $\bQ_0$ with information from the first two kernels $H_1,~H_2$. Using parameterization with $\lambda_i,~i=1,\ldots,99$ from \cref{eq:Qsk}, we enforce interpolation with the third kernel $H_3$ to estimate the remaining parameters $m=99$. When the data matrices in \cref{eq:QuadVec} are formed, and the conjugate counterparts are included to enforce the real symmetry, the dimensions remain: $\bW\in\IR^{432\times 99^2},~\bZ\in\IR^{432\times 99},~\bS\in\IR^{432\times 1},~\blambda\in\IR^{99\times 1}$.

To solve \cref{eq:QuadVec} with the above dimensionality, as the system is (potentially) over-determined with $99$ unknowns from the $432$ equations, we have to use the projection scheme in \cref{rem:ProjectQVE}. That leads to 
\begin{equation}\label{eq:327}
    \hat{\bW}(\blambda\otimes\blambda)+\hat{\bZ}\blambda+\hat{\bS}=\mathbf{0},~\hat{\bW}\in\IR^{99\times 99^2},~\hat{\bZ}\in\IR^{99\times 99},~\hat{\bS}\in\IR^{99\times 1},~\blambda\in\IR^{99\times 1}.
\end{equation}

We seek a solution of \cref{eq:327} with \cref{alg:qve} by setting the hyper-parameter tunned as $\eta=1e-8,~\gamma_0=1e-5$. The residue of Newton iterations stagnates at $1.9343e-06$. In \cref{fig:BurgersH2H3}(right), a comparison between the 3rd level kernels FOM($n=257$) and ROM($r=6$) is depicted. The updated quadratic matrix $\hat{\bQ}_r$ is obtained from \cref{eq:Qsk}.

\begin{figure}[h!]
    \centering
    \includegraphics[scale=0.3]{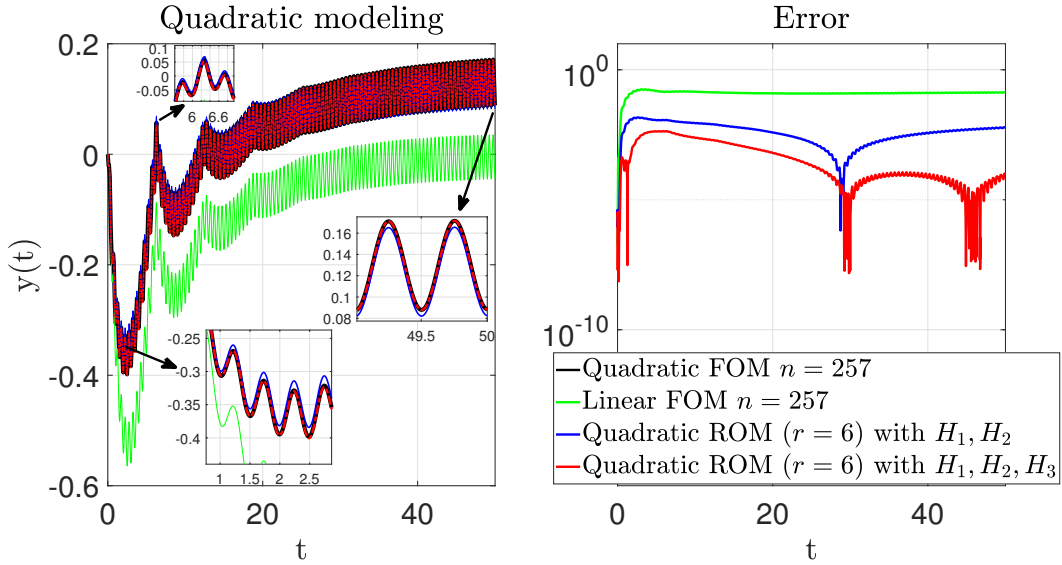}
    \caption{Reduction of a percentage level $98\%$ has been achieved with the error between the  FOM$(n=257)$ and ROM$(r=6)$ to be $\lVert\hat{y}(t)-y(t)\rVert_2\approx 1e-1$ and $\lVert\hat{y}(t)-y(t)\rVert_{max}\approx 1e-3$.}
    \label{fig:Burgers}
\end{figure}

In \cref{fig:Burgers}, both systems FOM($n=257$) and ROM($r=6$) are compared under a nontrivial test input i.e., $u(t)=0.5\exp(-0.2t)\texttt{sawtooth}(t)+0.5\sin(4\pi t)$, and for sufficient time, encompassing the "strong" transient behavior before finally reaching the quasi-steady state response. By considering more measurements, i.e., from $H_3$, the fit accuracy improved significantly for the complete simulation. The proposed method achieves an approximate interpolation to all measurement data sets, and the accuracy is improved when using the three GFRFs \cref{fig:Burgers}(right). 

Finally, the updated quadratic matrix $\hat{\bQ}$ of dimension $6\times6^2$ interpolates approximately the $H_3$ without affecting the interpolation achieved in $H_2$ from $\hat{\bQ}_0$. To illustrate this result, in \cref{tab:inter}, we choose a random point, and we test the interpolation error for both estimations of the quadratic matrix; first, from the $H_1,~H_2$ as $\hat{\bQ}_0$; secondly, from the $H_1,~H_2,~H_3$ as $\hat{\bQ}_r$.

\begin{table}
    \centering
    \small
    \begin{tabular}{c|c|c}
       GFRFs & Evaluation at $(s_1,s_2,s_3)=(1j,2j,3j)$ & Interpolation with the FOM\\[2mm]\hline
        FOM $H_2(s_1,s_2)$ & $-0.2511 + 0.0004{}\mathrm{i}$ & theoretical value \\ 
          $\hat{H}_2(s_1,s_2,\hat{\bQ}_r)$ & $-0.2511 + 0.0004{}\mathrm{i}$ & \checkmark \\
         $\hat{H}_2(s_1,s_2,\hat{\bQ}_0)$ & $-0.2511 + 0.0004{}\mathrm{i}$ & \checkmark\\[2mm]\hline
        FOM  $H_3(s_1,s_2,s_3)$ & $\mathbf{0.0042 + 0.0448{}\mathrm{i}}$ & theoretical value\\ 
          $\hat{H}_3(s_1,s_2,s_3,\hat{\bQ}_r,\hat{\bQ}_r)$ & $\mathbf{0.0042 + 0.0448{}\mathrm{i}}$ & \checkmark \\
         $\hat{H}_3(s_1,s_2,s_3,\hat{\bQ}_0,\hat{\bQ}_0)$ & $0.0122 + 0.0249{}\mathrm{i}$ & $\times$
    \end{tabular}
    \normalsize
    \caption{Symmetric Volterra kernel interpolation at a test point. The updated estimation $\hat{\bQ}_r$ from the three kernels enforces the interpolation to the third kernel. As a result, the overall fitting performance has improved significantly \cref{fig:Burgers}.}
    \label{tab:inter}
\end{table}

\section{Discussion and concluding remarks}\label{sec:conclusions}
In this study, we were concerned with inferring quadratic state-space models from i/o time-domain data from measurements collected over local stable operational quasi-steady state profiles. The proposed method successfully inferred quadratic models after estimating the first three symmetric GFRFs within the Volterra framework from i/o time-domain data in conjunction with the developed rational interpolation method known as the Loewner framework and nonlinear optimization methods. The inferred quadratic models have the potential to achieve global identification of low-order systems that could operate and be measured only locally, as we detailed with the Lorenz example. Another aspect of this method towards model inference enables the construction of accurate reduced quadratic models as surrogates of a potentially unknown large/complex model where only i/o measurements could be collected, as illustrated with the viscous Burgers' model that reduction achieved with a factor of $0.98$ and normalized error $\approx 1e-3$. Such efficient reduction can assess real-time implementation suitable for multi-query tasks such as optimization and control.

The current method could be extended to generalize its features to incorporate higher-order GFRFs, increasing the computational complexity in both algorithmic and measurement aspects. In contrast to this theoretical open challenge, the practicality of the derived method is connected strongly with the common engineering regularity in real applications. In particular, for the models with Volterra representation that exhibit weakly nonlinear dynamics, higher kernel information that comes with higher harmonics can be considered negligible in practice. The Volterra series representation already excludes chaotic dynamics and systems with continuous spectra, where the extension to handle such limitations and theoretical establishment of convergence properties tailored with hyperbolic PDEs would be interesting.

Moreover, the inference task via the proposed method could be challenging in cases where the linear subsystem, compared to the nonlinear, has a different minimal representation with dimension mismatch. Fortunately, when a quadratic model operates around a nonzero equilibrium, the linear minimality is affected by the quadratic operator (e.g., \cref{eq:syssq}), which can significantly increase and much better the quadratic minimality. The tools used in this study are robust to noise (such as most spectral transforms); a more involved analysis of the noise impact is left for future studies.

Finally, machine learning techniques can be advantageous to methods such as the proposed one due to their power to learn nonlinear i/o maps. For instance, when solely one input-output sequence of measurements is accessible, a neural network (NN) can be used as a surrogate black-box model to transfer the whole measurement process to more efficient cheap simulations. Finally, connecting data science and computational science tools, e.g., NNs, with variants of the proposed method will increase interpretability in ML approaches, enhancing the safety of utilizing such predictive models.

\appendix
\section{Coordinate tranformation}\label{app:align}
Dynamical systems that form an equivalent modulo can be aligned using a similarity transform $\bPsi\in\IR^{n\times n}$ as:
\begin{equation}
    \begin{aligned}
        \bPhi_1&=\left[\begin{array}{cccc}
        (\bC_1\bA_1)^T & (\bC_1\bA_1^2)^T & \cdots & (\bC_1\bA_1^n)^T  \end{array}\right],~\bPhi_2=\left[\begin{array}{cccc}
        (\bC_2\bA_2)^T & (\bC_2\bA_2^2)^T & \cdots & (\bC_2\bA_2^n)^T  \end{array}\right]^T,\\
        \bPsi&=\bPhi_2^{-1}\bPhi_1,~\text{and for the operators of the quadratic system holds:}\\
        &(\bA_2,~\bQ_2,~\bB_2,~\bC_2)=(\bPsi^{-1}\bA_1\bPsi,~\bPsi^{-1}\bQ_1(\bPsi\otimes\bPsi),~\bPsi^{-1}\bB_1,~\bC_1\bPsi).
    \end{aligned}
\end{equation}

\section{Generalized Frequency Response Functions (GFRFs)}\label{app:GFRFs}
To determine the $R$th order generalized frequency response function $H_{R}(s_{1},\ldots,s_{R})$, the probing input $u(t)=\sum_{i=1}^{R}e^{s_i t}$ must be applied, with at least $R$ harmonics. We introduce the input-state GFRFs $\bG_i(s_1,\ldots,s_i),~i=1,...,n$ to simplify the next derivations. These simply result in the input-output GFRFs by multiplying the $\bH_i$'s from the left by the output vector $\bC$ (in the SISO case), that is, $H_n(s_1,\ldots,s_n)=\bC\bG_n(s_1,\ldots,s_n)$. Note further that the transfer function $\bG_i$ is a vector of length equal to the state dimension $n$ (this latter is identical to $\bH_i$, when $\bC = \bI_n$, i.e., when all the state elements are individually observed).\\
$\bullet$ \underline{\textbf{$R=1$ - 1st order GFRF $H_{1}(s_1)$:}} With input $u(t)=e^{s_{1}t}$ the state solution $\bx(t)$  and the time derivative $\dot{\bx}(t)$ are respectively:
\begin{equation}
\bx(t)=\sum_{n=1}^{\infty}\sum_{m(n)}\tilde{\bG}_{m_{1}(n)}(s_{1})e^{(m_{1}(n)s_1)t},~\dot{\bx}(t)=\sum_{n=1}^{\infty}\sum_{m(n)}\tilde{\bG}_{m_{1}(n)}(s_{1})m_1(n)s_1e^{(m_{1}(n)s_1)t}
\end{equation}
By substituting to the differential equation of the system \cref{eq:qsys}, we have
\begin{equation}
    \begin{aligned}
    &\dot{\bx}(t)-\bA\bx(t)-\bB u(t)=\sum_{n=1}^{\infty}\sum_{m(n)}(m_1(n)s_1\bI-\bA)\tilde{\bG}_{m_{1}(n)}(s_{1})e^{(m_{1}(n)s_1)t}-\bB e^{s_1 t}\\
    &=\bQ\left(\sum_{n=1}^{\infty}\sum_{m(n)}\tilde{\bG}_{m_{1}(n)}(s_{1})e^{(m_{1}(n)s_1)t}\otimes\sum_{n=1}^{\infty}\sum_{m(n)}\tilde{\bG}_{m_{1}(n)}(s_{1})e^{(m_{1}(n)s_1)t}\right).
    \end{aligned}
\end{equation}
By collecting the terms with $e^{s_1 t}$, we result to
\begin{equation}
    (s_1\bI-\bA)\tilde{\bG}_1(s_1)=\bB\Rightarrow\tilde{\bG}_1(s_1)=(s_1\bI-\bA)^{-1}\bB.
\end{equation}
We adjust the weighted $\bG_1(s_1)=\frac{1!}{m_1(1)}\tilde{\bG}_1(s_1)=(s_1\bI-\bA)^{-1}\bB$ in a similar way to $H_n$. Multiplication by the vector $\bC$ from the left gives the first-order GFRF consistent with the linear subsystem. Using the resolvent notation $(\bE=\bI)$,
\begin{equation}
    H_1(s_1)=\bC(s_1\bI-\bA)^{-1}\bB=\bC\underbrace{\bPhi(s_1)\overbrace{\bB}^{\bR_1}}_{\bG_1(s_1)}.
\end{equation}
Higher-order kernels ($H_2, H_3,\ldots$) can be derived at this point. Still, these can be evaluated only on the diagonal of the hyperplane that spans the domain of definition, e.g., $H_2(s_1,s_1)$, and $H_3(s_1,s_1,s_1)$ which is not enough for achieving the identification goal as $H_2$ has a 2D domain support where a single harmonic input will always give information on the univariate diagonal (NFR method). Therefore, the next step is to excite the system with more complex harmonic inputs to identify the structure of the higher kernels.\\
$\bullet$ \underline{\textbf{$R=2$ - 2nd order GFRF $H_{2}(s_1,s_2)$:}} With input $u(t)=e^{s_{1}t}+e^{s_{2}t}$ the state solution is:
\begin{equation}
    \bx(t)=\sum_{n=1}^{\infty}\sum_{m(n)}\tilde{\bG}_{m_{1}(n)m_{2}(n)}(s_{1},s_{2})e^{(m_{1}(n)s_1+m_{2}(n)s_{2})t},
\end{equation}
and the time derivative results to
\begin{equation}
    \dot{\bx}(t)=\sum_{n=1}^{\infty}\sum_{m(n)}(m_1(n)s_1+m_2(n)s_2)\tilde{\bG}_{m_{1}(n)m_{2}(n)}(s_{1},s_{2})e^{(m_{1}(n)s_1+m_{2}(n)s_{2})t}.
\end{equation}
By substituting into the differential equation of the system \eqref{eq:qsys}, we obtain
\begin{equation}\label{eq:G2probingsystem}
    \begin{aligned}
    &\dot{\bx}(t)-\bA\bx(t)=\bQ\left(\bx(t)\otimes\bx(t)\right)+\bB u(t)\Leftrightarrow\\
    &\sum_{n=1}^{\infty}\sum_{m(n)}((m_1(n)s_1+m_2(n)s_2)\bI-\bA)\tilde{\bG}_{m_{1}(n)m_2(n)}(s_{1},s_2)e^{(m_{1}(n)s_1+m_{2}(n)s_2)t}\\
    &=\bQ\left(\sum_{n=1}^{\infty}\sum_{m(n)}\tilde{\bG}_{m_{1}(n)m_{2}(n)}(s_{1},s_{2})e^{(m_{1}(n)s_1+m_{2}(n)s_{2})t}\right.\otimes...\\
    &...\left.\otimes\sum_{n=1}^{\infty}\sum_{m(n)}\tilde{\bG}_{m_{1}(n)m_{2}(n)}(s_{1},s_{2})e^{(m_{1}(n)s_1+m_{2}(n)s_{2})t}\right)+\bB(e^{s_1 t}+e^{s_2 t}).
    \end{aligned}
\end{equation}
By collecting the terms $e^{s_1 t+s_2 t}$ with $(n=2,~m_1(2)=1,~m_2(2)=1)$, we result to
\begin{equation}\label{eq:G2derivation}
\begin{aligned}
&\left((s_1+s_2)\bI-\bA)\right)\tilde{\bG}_{11}(s_1,s_2)e^{(s_1+s_2)t}=\\
&=\bQ\left[\tilde{\bG}_{10}(s_{1})e^{s_1 t}\otimes\tilde{\bG}_{01}(s_{2})e^{s_2 t}+\tilde{\bG}_{01}(s_{2})e^{s_2 t}\otimes\tilde{\bG}_{10}(s_{1})e^{s_1 t}\right]\Rightarrow\\
&\left((s_1+s_2)\bI-\bA)\right)\tilde{\bG}_{11}(s_1,s_2)=\bQ\left[\tilde{\bG}_{10}(s_{1})\otimes\tilde{\bG}_{01}(s_{2})+\tilde{\bG}_{01}(s_{2})\otimes\tilde{\bG}_{10}(s_{1})\right]\Rightarrow\\
&\tilde{\bG}_{11}(s_{1},s_{2})=\left((s_1+s_2)\bI-\bA\right)^{-1}\cdot\\
&\cdot\bQ\left[(s_1\bI-\bA)^{-1}\bB\otimes(s_2\bI-\bA)^{-1}\bB+(s_2\bI-\bA)^{-1}\bB\otimes(s_1\bI-\bA)^{-1}\bB\right]
\end{aligned}
\end{equation}
Finally, by adjusting the weighted $\tilde{\bG}_{11}(s_1,s_2)=\frac{2!}{1!1!}\bG_{2}(s_1,s_2)=2\bG_{2}(s_1,s_2)$, and multiplying from the left with $\bC$, we can define the 2nd order GFRF after using the resolvent notation as
\begin{equation}
\begin{aligned}
    H_2(s_1,s_2)&=\frac{1}{2}\bC\bPhi(s_1+s_2)\bQ\underbrace{\left[\bPhi(s_1)\bB\otimes\bPhi(s_2)\bB+\bPhi(s_2)\bB\otimes\bPhi(s_1)\bB\right]}_{\bR_2(s_1,s_2)}\\
    &=\bC\underbrace{\frac{1}{2}\bPhi(s_1+s_2)\bQ\left[\bG_1(s_1)\otimes\bG_1(s_2)+\bG_1(s_2)\otimes\bG_1(s_1)\right]}_{\bG_2(s_1,s_2)}
    \end{aligned}
\end{equation}

$\bullet$ \underline{\textbf{$R=3$ - 3rd order GFRF $H_{3}(s_1,s_2,s_3)$:}} With input $u(t)=e^{s_{1}t}+e^{s_{2}t}+e^{s_{3}t}$, and  similar arguments, we can derive 
\begin{equation}
  \begin{aligned}
    H_3(s_1,s_2,s_3)&=\bC\underbrace{\frac{1}{6}\bPhi(s_1+s_2+s_3)\bQ\bR_3(s_1,s_2,s_3)}_{\bG_3(s_1,s_2,s_3)},~\text{with}\\
    \bR_3(s_1,s_2,s_3)=&\bG_1(s_1)\otimes\bG_2(s_2,s_3)+\bG_2(s_2,s_3)\otimes\bG_1(s_1)+\\
    &\bG_1(s_2)\otimes\bG_2(s_1,s_3)+\bG_2(s_1,s_3)\otimes\bG_1(s_2)+\\
    &\bG_1(s_3)\otimes\bG_2(s_1,s_2)+\bG_2(s_1,s_2)\otimes\bG_1(s_3).
    \end{aligned}
\end{equation}

\section{Solution of the constrained quadratic matrix equation}\label{sec:ConQME}
The analysis starts with the quadratic matrix equation. We define the following operator: $\mathcal{F}:\IR^{n\times n}\rightarrow\IR^{n\times n}$ with $\mathcal{F}(\bX):=\bX\bU-\bQ\left(\bX\otimes\bX\right)$. For known $\bU,~\bQ\in\IR^{n\times n^2}$, we seek $\bX\neq\mathbf{0}\in\IR^{n\times n}$ such that $\mathcal{F}(\bX)=\mathbf{0}$. Furthermore, $\bX$ should be invertible $(\exists~\bX^{-1})$. The idea is to differentiate w.r.t. the Fr\'echet derivative and solve a linear matrix equation for every Newton step similar to the Newton-Kleinmann algorithm for the solution of the Ricatti matrix equation \cite{Kleinman1968}. Therefore, we introduce a small perturbation to the matrix $\bX$ with $\bN\in\IR^{n\times n}$ and with $h$ a small real number. We define
\begin{equation}
    (\mathcal{F}'(\bX))(\bN)=\lim_{h\to 0}\frac{1}{h}\left(\mathcal{F}(\bX+h\bN)-\mathcal{F}(\bX)\right)=\bN\bU-\bQ(\bX\otimes\bN+\bN\otimes\bX).
\end{equation}
Since $\bQ$ is symmetric, we can write equivalently
\begin{equation}
    (\mathcal{F}'(\bX))(\bN)=\bN\bU-2\bQ(\bX\otimes\bN).
\end{equation}
The Newton iteration is given by
\begin{equation}
    (\mathcal{F}'(\bX_{j-1}))(\bN_{j-1})=-\mathcal{F}(\bX_{j-1}),~\bX_j=\bX_{j-1}+\bN_{j-1}.
\end{equation}
We compute
\begin{equation}
\begin{aligned}
    \bN_{j-1}\bU-2\bQ(\bX_{j-1}\otimes\bN_{j-1})&=-\bX_{j-1}\bU+\bQ(\bX_{j-1}\otimes\bX_{j-1})\Rightarrow\\
    (\bX_{j}-\bX_{j-1})\bU-2\bQ(\bX_{j-1}\otimes(\bX_{j}-\bX_{j-1}))&=-\bX_{j-1}\bU+\bQ(\bX_{j-1}\otimes\bX_{j-1})\Rightarrow\\
    \bX_{j}\bU-2\bQ(\bX_{j-1}\otimes\bX_{j})+2\bQ(\bX_{j-1}\otimes\bX_{j-1})&=\bQ(\bX_{j-1}\otimes\bX_{j-1})
     \end{aligned}
\end{equation}
which results to the following linear matrix equation \cref{eq_Xj_Xj-1} w.r.t the forward step solution $\bX_j$:
\begin{equation}\label{eq_Xj_Xj-1}
\begin{aligned}
     \bX_{j}\bU-2\bQ(\bX_{j-1}\otimes\bX_{j})+\bQ(\bX_{j-1}\otimes\bX_{j-1})&=0.
     \end{aligned}
\end{equation}

\begin{remark}
In \cref{eq_Xj_Xj-1}, it should be observed that in step $j$, the matrix equation is linear in $\bX_{j}$, provided that $\bX_{j-1}$ is explicitly known, which is assumed (from the Newton iteration).
 \end{remark}
 \begin{remark}
\cref{eq_Xj_Xj-1} is linear w.r.t $\bX_{j} \in \IR^{n \times n}$; since $\bU,\bQ \in \IR^{n \times n^2}$, there are $n^3$ linear scalar equations to solve and only $n^2$ unknowns. Hence, we have to deal with an over-determined linear system of equations with a possibly non-empty null space.
 \end{remark}

In what follows, we show how to isolate the $\bX_j$ term from the rest and how to re-write this equation more conventionally. More specifically, based on the previous remark, we show that equation (\ref{eq_Xj_Xj-1}) can equivalently be written as $n$ classical Sylvester equations, each characterized by $n^2$ scalar equations in $n^2$ unknowns. From (\ref{eq_Xj_Xj-1}), it follows that
\small
 \begin{equation}
 \begin{aligned}
      \bX_{j}\bU-2\bQ(\bX_{j-1}\otimes\bX_{j})+\bQ(\bX_{j-1}\otimes\bX_{j-1})&=0\Rightarrow\\
       \bX_{j}\bU-\underbrace{2\bQ(\bX_{j-1}\otimes \bI_n)}_{:= \bV_{j-1}} (\bI_n \otimes \bX_{j})&=  \underbrace{-\bQ(\bX_{j-1}\otimes\bX_{j-1})}_{:= \bZ_{j-1}} \Rightarrow
      \end{aligned}
 \end{equation}
\small
\vspace{-5mm}
  \begin{equation}
      \begin{aligned}  
      \bX_j \underbrace{\begin{bmatrix} \bU^{(1)}  & \cdots \bU^{(n)} \end{bmatrix}}_{\bU} &- \underbrace{\begin{bmatrix} \bV_{j-1}^{(1)}  & \cdots \bV_{j-1}^{(n)} \end{bmatrix}}_{\bV_{j-1}}
      \begin{bmatrix} \bX_j & \bfz & \cdots & \bfz \\ \bfz & \bX_j & \cdots & \bfz \\ \vdots & \vdots & \ddots & \vdots \\ \bfz & \bfz & \cdots & \bX_j \end{bmatrix} 
= \underbrace{\begin{bmatrix} \bZ_{j-1}^{(1)} & \bZ_{j-1}^{(2)} & \cdots \bZ_{j-1}^{(n)} \end{bmatrix}}_{\bZ_{j-1}}.
      \end{aligned}
 \end{equation}
 \normalsize
 Above, we have that $\bU^{(k)}, \bV_{j-1}^{(k)}, \bZ_{j-1}^{(k)}$ are known $n \times n$ real-valued matrices at step $j$, for all $1 \leq k \leq n$. These are the building blocks of the following matrices:
 \begin{equation}
     \bV_{j-1} := 2\bQ(\bX_{j-1}\otimes \bI_n) \in \mathbb{R}^{n \times n^2}, \ \ \bZ_{j-1} :=- \bQ(\bX_{j-1}\otimes\bX_{j-1}) \in \mathbb{R}^{n \times n^2}.
 \end{equation}

We can hence write this equation equivalently as follows:
 \begin{align}
          \begin{bmatrix}  \bX_j \bU^{(1)} &  \cdots  \bX_j \bU^{(n)} \end{bmatrix} &- \begin{bmatrix} \bV_{j-1}^{(1)} & \cdots \bV_{j-1}^{(n)}  \bX_j \end{bmatrix}=\begin{bmatrix} \bZ_{j-1}^{(1)} & \cdots \bZ_{j-1}^{(n)} \end{bmatrix}.
 \end{align}
 Then, for all $1 \leq k \leq n$, solving (\ref{eq_Xj_Xj-1}) boils down to solving $n$ (linear) Sylvester equations as:
 \begin{align}
     \bX_j \bU^{(k)} - \bV_{j-1}^{(k)} \bX_j = \bZ_{j-1}^{(k)}.
 \end{align}
The solution $\bX_j \in \IR^{n \times n}$, after vectorization, becomes $\text{vec}(\bX_j) \in \IR^{n^2 \times 1}$. Putting together the $n$ Sylvester equations in vectorized form by using the identity $\text{vec}(\bT \bO \bR) = (\bR^T \otimes \bT) \text{vec}(\bO)$,
will yield the following system of $n^3$ scalar equations in $n^2$ unknowns:
\small
\begin{equation}\label{eq:LinSysBig}
\underbrace{
  \begin{bmatrix}
        \left(\bU^{(1)} \right)^T \otimes \bI_n - \bI_n \otimes \bV_{j-1}^{(1)} \\
              \left(\bU^{(2} \right)^T \otimes \bI_n - \bI_n \otimes \bV_{j-1}^{(2)}\\
              \vdots \\
                    \left(\bU^{(n)} \right)^T \otimes \bI_n - \bI_n \otimes \bV_{j-1}^{(n)}
  \end{bmatrix}}_{\in \IR^{n^3 \times n^2}} \text{vec}(\bX_j) = \underbrace{\begin{bmatrix}
      \text{vec}(\bZ_{j-1}^{(1)}) \\[2mm]
      \text{vec}(\bZ_{j-1}^{(2)}) \\
      \vdots \\
      \text{vec}(\bZ_{j-1}^{(n)})
  \end{bmatrix}}_{\in \IR^{n^3 \times 1}}
\end{equation}
\normalsize
For low values of $n$, such a procedure is indeed feasible. However, for moderate to large values of $n$, i.e., $n>50$, it is quite challenging or even impossible to find the next value $\bX_j$, by explicitly forming the matrix $n^3 \times n^2$ in (\ref{eq:LinSysBig}). In what follows, we are concerned with low-order systems, as we emphasize quadratic identification in a reduced-order sense.

\begin{lemma}
The square matrix $\bT^{-1}$ that aligns the operators $(\hat{\bQ}_1,\hat{\bB}_1,\hat{\bC}_1)$ and $(\breve{\bQ}_2,\breve{\bB}_2,\breve{\bC}_2)$ from \cref{lem:2.6}, can be computed, on the convergence of Newton's method \cref{alg:qme}, as the iterative solution of the following constrained linear system of equations \cref{eq:LinSysBig2} with $\bT^{-1}:= \lim_{j\to\infty}\bX_{j}$ which gives $\mathcal{F}(\bT^{-1})\approx\mathbf{0}$.
\begin{equation}\label{eq:LinSysBig2}
\underbrace{
  \begin{bmatrix}
        \left(\bU^{(1)} \right)^T \otimes \bI_n - \bI_n \otimes \bV_{j-1}^{(1)} \\
              \left(\bU^{(2} \right)^T \otimes \bI_n - \bI_n \otimes \bV_{j-1}^{(2)}\\
              \vdots \\
                    \left(\bU^{(n)} \right)^T \otimes \bI_n - \bI_n \otimes \bV_{j-1}^{(n)} \\ \hat{\bB}_1^T \otimes \bI_n \\[1mm] \bI_n\otimes\breve{\bC}_2
  \end{bmatrix}}_{\in \IR^{(n^3+2n) \times n^2}} \underbrace{\text{vec}(\bX_j)}_{\in \IR^{n^2 \times 1}} = \underbrace{\begin{bmatrix}
      \text{vec}(\bZ_{j-1}^{(1)}) \\[2mm]
      \text{vec}(\bZ_{j-1}^{(2)}) \\
      \vdots \\
      \underbrace{\text{vec}(\bZ_{j-1}^{(n)})}_{\in \IR^{n^2 \times 1}} \\[1mm]
      \breve{\bB}_2 \\[1mm]
      \hat{\bC}_1^T
  \end{bmatrix}}_{\in\IR^{(n^3+2n) \times 1}}
\end{equation} 
\end{lemma}

\normalsize

\begin{algorithm}
\caption{Solution of the constrained quadratic matrix equation with Newton method}
\label{alg:qme}
\begin{algorithmic}
\STATE{Seek: $\bX$ s.t. $\mathcal{F}(\bX):=\bX\bU-\bQ(\bX\otimes\bX)=\mathbf{0}$ and satisfies the  constrains (two last rows) in \cref{eq:LinSysBig2}.} 
\STATE{Choose an initial random seed: $\bX_{j=0}\in\IR^{n\times n}$.}
\WHILE{$\gamma>\gamma_0$}
\STATE{Update: $j\leftarrow j+1$.}
\STATE{Compute $\bX_j$ by solving the linear system of equations
\cref{eq:LinSysBig2}.}
\STATE{Compute the residue $\lVert\mathcal{F}(\bX_j)\rVert=\gamma$.}
\ENDWHILE
\RETURN $\bX$
\end{algorithmic}
\end{algorithm}

We can align the two systems in \cref{eq:syss} explained in \cref{lem:2.6} to the same coordinates with the solution $\bT^{-1}$ obtained from \cref{alg:qme}. 

\normalsize

\section*{Acknowledgments} 
The first author would like to acknowledge that this work was supported mainly by the Max Planck Institute in Magdeburg during his Ph.D. project.

\normalsize
\bibliographystyle{siamplain}
\bibliography{qpaper}

\begin{thebibliography}{10}

\bibitem{aa90}
{\textsc B.~D.~O. Anderson and A.~C. Antoulas}, {\em Rational interpolation and
  state-variable realizations}, Linear Algebra and Its Applications, 137/138
  (1990), pp.~479--509.

\bibitem{ACA05}
{\textsc A.~C. Antoulas}, {\em Approximation of large-scale dynamical systems},
  SIAM, Philadelphia, 2005.

\bibitem{ABG20}
{\textsc A.~C. Antoulas, C.~A. Beattie, and S.~Gugercin}, {\em Interpolatory
  Methods for Model Reduction}, SIAM, Philadelphia, 2020.

\bibitem{morAntGH19}
{\textsc A.~C. Antoulas, I.~V. Gosea, and M.~Heinkenschloss}, {\em On the
  {L}oewner framework for model reduction of {B}urgers' equation}, in Active
  Flow and Combustion Control, R.~King, ed., Springer, Cham, Switzerland, 2019,
  pp.~255--270.

\bibitem{AGI16}
{\textsc A.~C. Antoulas, I.~V. Gosea, and A.~C. Ionita}, {\em Model reduction
  of bilinear systems in the {L}oewner framework}, SIAM Journal on Scientific
  Computing, 38(5) (2016), pp.~B889--B916.

\bibitem{birkjour}
{\textsc A.~C. Antoulas, S.~Lefteriu, and A.~C. Ionita}, {\em A tutorial
  introduction to the {L}oewner Framework for Model Reduction}, Model Reduction
  and Approximation for Complex Systems, Edited by P. Benner, A. Cohen, M.
  Ohlberger, and K. Willcox,Series: Computational Science $\&$ Engineering,
  \url{https://doi.org/10.1137/1.9781611974829}, SIAM, 2017, pp.~335--376.

\bibitem{baur2014model}
{\textsc U.~Baur, P.~Benner, and L.~Feng}, {\em Model order reduction for
  linear and nonlinear systems: a system-theoretic perspective}, Archives of
  Computational Methods in Engineering, 21 (2014), pp.~331--358.

\bibitem{benner2021model}
{\textsc P.~Benner and L.~Feng}, {\em Model order reduction based on
  moment-matching}, in Model Order Reduction: Volume 1: System-and Data-Driven
  Methods and Algorithms, De Gruyter, 2021, pp.~57--96.

\bibitem{PSGoyBen}
{\textsc P.~Benner and P.~Goyal}, {\em Interpolation-based model order
  reduction for polynomial systems}, SIAM Journal on Scientific Computing, 43
  (2021), pp.~A84--A108, \url{https://doi.org/10.1137/19M1259171}.

\bibitem{morBenGKetal20}
{\textsc P.~Benner, P.~Goyal, B.~Kramer, B.~Peherstorfer, and K.~Willcox}, {\em
  Operator inference for non-intrusive model reduction of systems with
  non-polynomial nonlinear terms}, CompMethAppMechEng, 372 (2020), p.~113433,
  \url{https://doi.org/10.1016/j.cma.2020.113433}.

\bibitem{PhysRevE.91.032915}
{\textsc T.~Berry, D.~Giannakis, and J.~Harlim}, {\em Nonparametric forecasting
  of low-dimensional dynamical systems}, Phys. Rev. E, 91 (2015), p.~032915,
  \url{https://doi.org/10.1103/PhysRevE.91.032915},
  \url{https://link.aps.org/doi/10.1103/PhysRevE.91.032915}.

\bibitem{Billings2013NonlinearSI}
{\textsc S.~A. Billings}, {\em Nonlinear system identification: Narmax methods
  in the time, frequency, and spatio-temporal domains}, 2013,
  \url{https://api.semanticscholar.org/CorpusID:117105230}.

\bibitem{Billings2003Subharmonic}
{\textsc O.~Boaghe and S.~Billings}, {\em Subharmonic oscillation modeling and
  miso volterra series}, IEEE Transactions on Circuits and Systems I:
  Fundamental Theory and Applications, 50 (2003), pp.~877--884,
  \url{https://doi.org/10.1109/TCSI.2003.813965}.

\bibitem{BoydChua85}
{\textsc S.~Boyd and L.~S. Chua}, {\em Fading memory and the problem of
  approximating nonlinear operators with volterra series}, IEEE Transactions on
  Circuits and Systems, 30 (1985), pp.~1150--1161.

\bibitem{Boyd1983MeasuringVK}
{\textsc S.~Boyd, Y.~Tang, and L.~Chua}, {\em Measuring volterra kernels}, IEEE
  Transactions on Circuits and Systems, 30 (1983), pp.~571--577,
  \url{https://doi.org/10.1109/TCS.1983.1085391}.

\bibitem{breiten20212}
{\textsc T.~Breiten and T.~Stykel}, {\em 2 balancing-related model reduction
  methods}, System-and Data-Driven Methods and Algorithms,  (2021), pp.~15--56.

\bibitem{ArnoldiNaka}
{\textsc P.~D. Brubeck, Y.~Nakatsukasa, and L.~N. Trefethen}, {\em Vandermonde
  with arnoldi}, SIAM Review, 63 (2021), pp.~405--415,
  \url{https://doi.org/10.1137/19M130100X},
  \url{https://doi.org/10.1137/19M130100X},
  \url{https://arxiv.org/abs/https://doi.org/10.1137/19M130100X}.

\bibitem{brunton2022data}
{\textsc S.~Brunton and J.~Kutz}, {\em Data-Driven Science and Engineering:
  Machine Learning, Dynamical Systems, and Control}, Cambridge University
  Press, 2022, \url{https://books.google.de/books?id=rxNkEAAAQBAJ}.

\bibitem{SINDY}
{\textsc S.~L. Brunton, J.~L. Proctor, and J.~N. Kutz}, {\em Discovering
  governing equations from data by sparse identification of nonlinear dynamical
  systems}, Proceedings of the National Academy of Sciences, 113 (2016),
  pp.~3932--3937, \url{https://doi.org/10.1073/pnas.1517384113}.

\bibitem{CHEN2000}
{\textsc H.~Chen and J.~Maciejowski}, {\em Subspace identification method for
  combined deterministic-stochastic bilinear systems}, IFAC Proceedings
  Volumes, 33 (2000), pp.~229--234,
  \url{https://doi.org/https://doi.org/10.1016/S1474-6670(17)39755-0}.
\newblock 12th IFAC Symposium on System Identification (SYSID 2000), Santa
  Barbara, CA, USA, 21-23 June 2000.

\bibitem{Chen22}
{\textsc N.~Chen, H.~Liu, and F.~Lu}, {\em {Shock trace prediction by reduced
  models for a viscous stochastic Burgers equation}}, Chaos: An
  Interdisciplinary Journal of Nonlinear Science, 32 (2022), p.~043109,
  \url{https://doi.org/10.1063/5.0084955},
  \url{https://doi.org/10.1063/5.0084955},
  \url{https://arxiv.org/abs/https://pubs.aip.org/aip/cha/article-pdf/doi/10.1063/5.0084955/16452495/043109\_1\_online.pdf}.

\bibitem{BillingsNARMAX}
{\textsc S.~Chen and S.~A. Billings}, {\em Representations of non-linear
  systems: the narmax model}, International Journal of Control, 49 (1989),
  pp.~1013--1032, \url{https://doi.org/10.1080/00207178908559683}.

\bibitem{Chorin2015}
{\textsc A.~J. Chorin and F.~Lu}, {\em Discrete approach to stochastic
  parametrization and dimension reduction in nonlinear dynamics}, Proceedings
  of the National Academy of Sciences, 112 (2015), pp.~9804--9809,
  \url{https://doi.org/10.1073/pnas.1512080112},
  \url{https://www.pnas.org/doi/abs/10.1073/pnas.1512080112},
  \url{https://arxiv.org/abs/https://www.pnas.org/doi/pdf/10.1073/pnas.1512080112}.

\bibitem{HandbookVol2}
{\textsc P.~B. et~al.}, ed., {\em Model Order Reduction Volume 2:
  Snapshot-Based Methods and Algorithms}, De Gruyter, Berlin, Boston, 2020,
  \url{https://doi.org/doi:10.1515/9783110671490}.

\bibitem{HandbookVol3}
{\textsc P.~B. et~al.}, ed., {\em Model Order Reduction Volume 3:
  Applications}, De Gruyter, Berlin, Boston, 2020,
  \url{https://doi.org/doi:10.1515/9783110499001}.

\bibitem{HandbookVol1}
{\textsc P.~B. et~al.}, ed., {\em Model Order Reduction Volume 1: System- and
  Data-Driven Methods and Algorithms}, De Gruyter, Berlin, Boston, 2021,
  \url{https://doi.org/doi:10.1515/9783110498967}.

\bibitem{Favoreel}
{\textsc W.~Favoreel, B.~De~Moor, and P.~Van~Overschee}, {\em Subspace
  identification of bilinear systems subject to white inputs}, IEEE
  Transactions on Automatic Control, 44 (1999), pp.~1157--1165,
  \url{https://doi.org/10.1109/9.769370}.

\bibitem{Gosea22mdpi}
{\textsc I.~V. Gosea}, {\em Exact and inexact lifting transformations of
  nonlinear dynamical systems: Transfer functions, equivalence, and complexity
  reduction}, Applied Sciences, 12 (2022),
  \url{https://doi.org/10.3390/app12052333}.

\bibitem{morGA15}
{\textsc I.~V. Gosea and A.~C. Antoulas}, {\em Model reduction of linear and
  nonlinear systems in the loewner framework: A summary}, in 2015 European
  Control Conference (ECC), 2015, pp.~345--349,
  \url{https://doi.org/10.1109/ECC.2015.7330568}.

\bibitem{morGosA18}
{\textsc I.~V. Gosea and A.~C. Antoulas}, {\em Data-driven model order
  reduction of quadratic-bilinear systems}, Numerical Linear Algebra with
  Applications, 25 (2018), p.~e2200, \url{https://doi.org/10.1002/nla.2200}.

\bibitem{gosea2021learning}
{\textsc I.~V. Gosea, D.~S. Karachalios, and A.~C. Antoulas}, {\em Learning
  reduced-order models of quadratic dynamical systems from input-output data},
  in 2021 European Control Conference (ECC), IEEE, 2021, pp.~1426--1431.

\bibitem{GPAswitch}
{\textsc I.~V. Gosea, M.~Petreczky, and A.~C. Antoulas}, {\em Data-driven model
  order reduction of linear switched systems in the loewner framework}, SIAM
  Journal on Scientific Computing, 40 (2018), pp.~B572--B610,
  \url{https://doi.org/10.1137/17M1120233}.

\bibitem{GPA21CDC}
{\textsc I.~V. Gosea, M.~Petreczky, and A.~C. Antoulas}, {\em Reduced-order
  modeling of {LPV} systems in the {L}oewner framework}, in 2021 60th IEEE
  Conference on Decision and Control (CDC), 2021, pp.~3299--3305,
  \url{https://doi.org/10.1109/CDC45484.2021.9683742}.

\bibitem{Gu11}
{\textsc C.~Gu}, {\em A projection-based nonlinear model order reduction
  approach using quadratic-linear representation of nonlinear systems}, IEEE
  Transactions on Computer-Aided Design of Integrated Circuits and Systems, 30
  (2011), pp.~1307--1320.

\bibitem{VF}
{\textsc B.~Gustavsen and A.~Semlyen}, {\em Rational approximation of frequency
  domain responses by vector fitting}, IEEE Transactions on Power Delivery, 14
  (1999), pp.~1052--1061, \url{https://doi.org/10.1109/61.772353}.

\bibitem{morHirHKetal20}
{\textsc S.~M. Hirsch, K.~D. Harris, J.~N. Kutz, and B.~W. Brunton}, {\em
  Centering data improves the dynamic mode decomposition}, {SIAM} J. Appl. Dyn.
  Syst., 19 (2020), pp.~1920--1955, \url{https://doi.org/10.1137/19M1289881}.

\bibitem{HOKALMAN}
{\textsc B.~L. Ho and R.~E. Kalman}, {\em Effective construction of linear
  state-variable models from input/output functions}, at -
  Automatisierungstechnik, 14 (1966), pp.~545--548,
  \url{https://doi.org/10.1524/auto.1966.14.112.545}.

\bibitem{Isidori1973}
{\textsc A.~Isidori}, {\em Direct construction of minimal bilinear realizations
  from nonlinear input-output maps}, IEEE Transactions on Automatic Control, 18
  (1973), pp.~626--631, \url{https://doi.org/10.1109/TAC.1973.1100424}.

\bibitem{JuangPappa}
{\textsc J.-N. Juang and R.~S. Pappa}, {\em An eigensystem realization
  algorithm for modal parameter identification and model reduction}, Journal of
  Guidance, Control, and Dynamics, 8 (1985), pp.~620--627,
  \url{https://doi.org/10.2514/3.20031}.

\bibitem{morKarGA19}
{\textsc D.~S. Karachalios, I.~V. Gosea, and A.~C. Antoulas}, {\em A bilinear
  identification-modeling framework from time domain data}, Proc. Appl. Math.
  Mech., 19 (2019), p.~e201900246,
  \url{https://doi.org/10.1002/pamm.201900246}.

\bibitem{morKarGA19a}
{\textsc D.~S. Karachalios, I.~V. Gosea, and A.~C. Antoulas}, {\em The Loewner
  framework for system identification and reduction}, De Gruyter, Berlin,
  Boston, 2021, pp.~181--228,
  \url{https://doi.org/doi:10.1515/9783110498967-006}.

\bibitem{morKarGA21a}
{\textsc D.~S. Karachalios, I.~V. Gosea, and A.~C. Antoulas}, {\em On bilinear
  time-domain identification and reduction in the {L}oewner framework}, in
  Model Reduction of Complex Dynamical Systems, vol.~171 of International
  Series of Numerical Mathematics, Birkh{\"a}user, Cham, 2021, pp.~3--30,
  \url{https://doi.org/10.1007/978-3-030-72983-7_1}.

\bibitem{KARACHALIOS20227}
{\textsc D.~S. Karachalios, I.~V. Gosea, and A.~C. Antoulas}, {\em A framework
  for fitting quadratic-bilinear systems with applications to models of
  electrical circuits}, IFAC-PapersOnLine, 55 (2022), pp.~7--12,
  \url{https://doi.org/https://doi.org/10.1016/j.ifacol.2022.09.064}.
\newblock 10th Vienna International Conference on Mathematical Modelling
  MATHMOD 2022.

\bibitem{KHODABAKHSHI2022114296}
{\textsc P.~Khodabakhshi and K.~E. Willcox}, {\em Non-intrusive data-driven
  model reduction for differential–algebraic equations derived from lifting
  transformations}, Computer Methods in Applied Mechanics and Engineering, 389
  (2022), p.~114296,
  \url{https://doi.org/https://doi.org/10.1016/j.cma.2021.114296}.

\bibitem{Kleinman1968}
{\textsc D.~Kleinman}, {\em On an iterative technique for {R}iccati equation
  computations}, IEEE Transactions on Automatic Control, 13 (1968),
  pp.~114--115, \url{https://doi.org/10.1109/TAC.1968.1098829}.

\bibitem{li2021fourier}
{\textsc Z.~Li, N.~Kovachki, K.~Azizzadenesheli, B.~Liu, K.~Bhattacharya,
  A.~Stuart, and A.~Anandkumar}, {\em Fourier neural operator for parametric
  partial differential equations}, 2021,
  \url{https://arxiv.org/abs/2010.08895}.

\bibitem{LIN2021109864}
{\textsc K.~K. Lin and F.~Lu}, {\em Data-driven model reduction, {W}iener
  projections, and the {K}oopman-{M}ori-{Z}wanzig formalism}, Journal of
  Computational Physics, 424 (2021), p.~109864,
  \url{https://doi.org/https://doi.org/10.1016/j.jcp.2020.109864}.

\bibitem{Ting2023}
{\textsc Y.~T. Lin, Y.~Tian, D.~Perez, and D.~Livescu}, {\em Regression-based
  projection for learning mori–zwanzig operators}, SIAM Journal on Applied
  Dynamical Systems, 22 (2023), pp.~2890--2926,
  \url{https://doi.org/10.1137/22M1506146},
  \url{https://doi.org/10.1137/22M1506146},
  \url{https://arxiv.org/abs/https://doi.org/10.1137/22M1506146}.

\bibitem{Ramos2009}
{\textsc P.~Lopes~dos Santos, J.~A. Ramos, and J.~L. Martins~de Carvalho}, {\em
  Identification of bilinear systems with white noise inputs: An iterative
  deterministic-stochastic subspace approach}, IEEE Transactions on Control
  Systems Technology, 17 (2009), pp.~1145--1153,
  \url{https://doi.org/10.1109/TCST.2008.2002041}.

\bibitem{Lorenz}
{\textsc E.~N. Lorenz}, {\em Deterministic nonperiodic flow}, Journal of
  Atmospheric Sciences, 20 (1963), pp.~130 -- 141,
  \url{https://doi.org/10.1175/1520-0469(1963)020<0130:DNF>2.0.CO;2}.

\bibitem{LuschBruntonKutz}
{\textsc B.~Lusch, J.~Kutz, and S.~Brunton}, {\em Deep learning for universal
  linear embeddings of nonlinear dynamics}, Nat Commun, 9 (2018),
  \url{https://doi.org/https://doi.org/10.1038/s41467-018-07210-0}.

\bibitem{ajm07}
{\textsc A.~J. Mayo and A.~C. Antoulas}, {\em A framework for the generalized
  realization problem}, Linear Algebra and Its Applications, 426 (2007),
  pp.~634--662.

\bibitem{NST18}
{\textsc Y.~Nakatsukasa, O.~S\`{e}te, and L.~N. Trefethen}, {\em The aaa
  algorithm for rational approximation}, SIAM Journal on Scientific Computing,
  40 (2018), pp.~A1494--A1522, \url{https://doi.org/10.1137/16M1106122}.

\bibitem{Pathak18}
{\textsc J.~Pathak, A.~Wikner, R.~Fussell, S.~Chandra, B.~R. Hunt, M.~Girvan,
  and E.~Ott}, {\em {Hybrid forecasting of chaotic processes: Using machine
  learning in conjunction with a knowledge-based model}}, Chaos: An
  Interdisciplinary Journal of Nonlinear Science, 28 (2018), p.~041101,
  \url{https://doi.org/10.1063/1.5028373},
  \url{https://doi.org/10.1063/1.5028373},
  \url{https://arxiv.org/abs/https://pubs.aip.org/aip/cha/article-pdf/doi/10.1063/1.5028373/19792240/041101\_1\_online.pdf}.

\bibitem{PGW17}
{\textsc B.~Peherstorfer, S.~Gugercin, and K.~Willcox}, {\em Data-driven
  reduced model construction with time-domain loewner models}, SIAM Journal on
  Scientific Computing, 39 (2017), pp.~A2152--A2178,
  \url{https://doi.org/10.1137/16M1094750}.

\bibitem{morPehW16}
{\textsc B.~Peherstorfer and K.~Willcox}, {\em Data-driven operator inference
  for nonintrusive projection-based model reduction}, Computer Methods in
  Applied Mechanics and Engineering, 306 (2016), pp.~196--215,
  \url{https://doi.org/10.1016/j.cma.2016.03.025}.

\bibitem{PEHERSTORFER2016196}
{\textsc B.~Peherstorfer and K.~Willcox}, {\em Data-driven operator inference
  for nonintrusive projection-based model reduction}, Computer Methods in
  Applied Mechanics and Engineering, 306 (2016), pp.~196 -- 215.

\bibitem{ProBruKut2016}
{\textsc J.~L. Proctor, S.~L. Brunton, and J.~N. Kutz}, {\em Dynamic mode
  decomposition with control}, SIAM Journal on Applied Dynamical Systems, 15
  (2016), pp.~142--161, \url{https://doi.org/10.1137/15M1013857}.

\bibitem{RAISSI2019686}
{\textsc M.~Raissi, P.~Perdikaris, and G.~Karniadakis}, {\em Physics-informed
  neural networks: A deep learning framework for solving forward and inverse
  problems involving nonlinear partial differential equations}, Journal of
  Computational Physics, 378 (2019), pp.~686--707,
  \url{https://doi.org/https://doi.org/10.1016/j.jcp.2018.10.045},
  \url{https://www.sciencedirect.com/science/article/pii/S0021999118307125}.

\bibitem{Ru82}
{\textsc W.~J. Rugh}, {\em Nonlinear system theory - The {V}olterra/{W}iener
  Approach}, The Johns Hopkins University Press, 1981.

\bibitem{schmid_2010}
{\textsc P.~J. SCHMID}, {\em Dynamic mode decomposition of numerical and
  experimental data}, Journal of Fluid Mechanics, 656 (2010), p.~5–28,
  \url{https://doi.org/10.1017/S0022112010001217}.

\bibitem{TothHossam}
{\textsc R.~Toth, H.~S. Abbas, and H.~Werner}, {\em On the state-space
  realization of lpv input-output models: Practical approaches}, IEEE
  Transactions on Control Systems Technology, 20 (2012), pp.~139--153,
  \url{https://doi.org/10.1109/TCST.2011.2109040}.

\bibitem{Xpar2006}
{\textsc Verspecht and D.~Root}, {\em Polyharmonic distortion modeling}, IEEE
  Microwave Magazine, 7 (2006), pp.~44--57,
  \url{https://doi.org/10.1109/MMW.2006.1638289}.

\bibitem{PODWeissJulien}
{\textsc J.~Weiss}, {\em A Tutorial on the Proper Orthogonal Decomposition},
  Aerospace Research Center, 2019, \url{https://doi.org/10.2514/6.2019-3333},
  \url{https://arc.aiaa.org/doi/abs/10.2514/6.2019-3333},
  \url{https://arxiv.org/abs/https://arc.aiaa.org/doi/pdf/10.2514/6.2019-3333}.

\bibitem{VoltMacrXpar}
{\textsc X.~Y.~Z. Xiong, L.~J. Jiang, J.~E. Schutt-Ainé, and W.~C. Chew}, {\em
  Volterra series-based time-domain macromodeling of nonlinear circuits}, IEEE
  Transactions on Components, Packaging and Manufacturing Technology, 7 (2017),
  pp.~39--49, \url{https://doi.org/10.1109/TCPMT.2016.2627601}.

\end{thebibliography}
\end{document}